\documentclass[12pt,a4paper,twoside]{article}

\newcommand{\pageformat}[6]{\setlength{\hoffset}{-1in}
                  \setlength{\voffset}{-1in}
                  \addtolength{\hoffset}{#5}
                            \addtolength{\voffset}{#6}
                            \setlength{\oddsidemargin}{#1}
                            \setlength{\evensidemargin}{#2}
                            \setlength{\textwidth}{\paperwidth}
                  \addtolength{\textwidth}{-\oddsidemargin}
                  \addtolength{\textwidth}{-\evensidemargin}
                  \addtolength{\textwidth}{-\marginparsep}
                  \addtolength{\textwidth}{-\marginparwidth}
                            \setlength{\topmargin}{#3}
                            \setlength{\textheight}{\paperheight}
                  \addtolength{\textheight}{-\topmargin}
                  \addtolength{\textheight}{-\headheight}
                  \addtolength{\textheight}{-\headsep}
                  \addtolength{\textheight}{-\footskip}
                  \addtolength{\textheight}{-#4}}
\pageformat{2cm}{3cm}{25mm}{25mm}{1pt}{0pt}

\usepackage{ifthen}
\newboolean{article}
    \setboolean{article}{true}
\newboolean{report}
\newboolean{book}
\newboolean{letter}
\newboolean{german}
\newboolean{italian}
\newboolean{nobaselinestretch}
\newboolean{nosectionappendix}
\newboolean{oldtoc}
\newboolean{nosectionequn}
\newboolean{notheorem}
\ifthenelse{\boolean{german}}{
    \usepackage{german}}{}

\usepackage[latin1]{inputenc}
\usepackage{amsmath}
\usepackage{amssymb}
\usepackage[mathscr]{eucal}

\ifthenelse{\boolean{notheorem}}{}{
    \usepackage{theorem}}

\ifthenelse{\boolean{nobaselinestretch}}{}{
    \renewcommand{\baselinestretch}{1.25}}

\newenvironment{env}[2]{\begin{#1}#2\end{#1}}{}
    \newcommand{\beq}[1]{\begin{env}{equation}{#1}}
    \newcommand{\beqn}[1]{\begin{env}{equation*}{#1}}
    \newcommand{\bal}[1]{\begin{env}{align}{#1}}
    \newcommand{\baln}[1]{\begin{env}{align*}{#1}}
    \newcommand{\bga}[1]{\begin{env}{gather}{#1}}
    \newcommand{\bgan}[1]{\begin{env}{gather*}{#1}}
    \newcommand{\bflal}[1]{\begin{env}{flalign}{#1}}
    \newcommand{\bflaln}[1]{\begin{env}{flalign*}{#1}}
    \newcommand{\bmu}[1]{\begin{env}{multline}{#1}}
    \newcommand{\bmun}[1]{\begin{env}{multline*}{#1}}
    \newcommand{\bsp}[1]{\begin{env}{split}{#1}}

    \newcommand{\eeq}{\end{env}}
    \newcommand{\eeqn}{\end{env}}
    \newcommand{\eal}{\end{env}}
    \newcommand{\ealn}{\end{env}}
    \newcommand{\ega}{\end{env}}
    \newcommand{\egan}{\end{env}}
    \newcommand{\eflal}{\end{env}}
    \newcommand{\eflaln}{\end{env}}
    \newcommand{\emu}{\end{env}}
    \newcommand{\emun}{\end{env}}
    \newcommand{\esp}{\end{env}}

\newcommand{\lf}{\vspace{2ex}}
\newcommand{\bulletline}[1][]{\lf\noindent~\hfill$\bullet\bullet\bullet$\hfill~

\lf\noindent\bf{#1}}

\renewcommand{\bf}[1]{\textbf{#1}}
\renewcommand{\it}[1]{\textit{#1}}

\renewcommand{\sf}[1]{\textsf{#1}}

\renewcommand{\tt}[1]{\texttt{#1}}
\newcommand{\hl}[1]{\bf{\it{#1}}}
\newcommand{\mrm}[1]{\mathrm{#1}}
\newcommand{\mbf}[1]{\mathbf{#1}}
\newcommand{\msf}[1]{\text{\small$\sf{#1}$}}

\newcommand{\cmc}[1]{\mathcal{#1}}
\newcommand{\eus}[1]{\mathscr{#1}}
\newcommand{\euf}[1]{\mathfrak{#1}}
\newcommand{\bb}[1]{\mathbb{#1}}

\newcommand{\mtiny}[1]{{\setlength{\arraycolsep}{.3ex}\text{\tiny$#1$}}}
\newcommand{\nbd}[1]{$#1$\nobreakdash--}
\newcommand{\ol}[1]{\overline{#1}}
\newcommand{\ul}[1]{\underline{#1}}
\newcommand{\wt}[1]{\widetilde{#1}}

\newcommand{\ve}{\varepsilon}

\newcommand{\vk}{\varkappa}
\newcommand{\vp}{\varphi}
\newcommand{\om}{\omega}
\newcommand{\Om}{\Omega}
\newcommand{\abs}[1]{\left\lvert#1\right\rvert}
\newcommand{\norm}[1]{\left\lVert#1\right\rVert}

\newcommand{\snorm}[1]{\norm{\smash{#1}}}
\newcommand{\sabs}[1]{\abs{\smash{#1}}}

\newcommand{\bfam}[1]{\bigl(#1\bigr)}

\newcommand{\AB}[1]{\langle#1\rangle}

\newcommand{\BAB}[1]{\Bigl\langle#1\Bigr\rangle}
\newcommand{\CB}[1]{\{#1\}}
\newcommand{\bCB}[1]{\bigl\{#1\bigr\}}
\newcommand{\BCB}[1]{\Bigl\{#1\Bigr\}}
\newcommand{\SB}[1]{[#1]}

\newcommand{\Matrix}[1]{\begin{pmatrix}#1\end{pmatrix}}

\newcommand{\tMatrix}[1]{\mtiny{\Matrix{#1}}}
\newcommand{\rtMatrix}[1]{\raisebox{.3ex}{\tMatrix{#1}}}

\newcommand{\sbar}[1]{\:\bar{#1}\:}
\newcommand{\sodot}{\sbar{\odot}}

\newcommand{\set}[2][]{
    \ifthenelse{\equal{#1}{}}{
        \CB{#2}}{
        \CB{#1~|~#2}}}
\newcommand{\bset}[2][]{
    \ifthenelse{\equal{#1}{}}{
        \bCB{#2}}{
        \bCB{#1~|~#2}}}
\newcommand{\Bset}[2][]{
    \ifthenelse{\equal{#1}{}}{
        \BCB{#2}}{
        \BCB{#1~\big|~#2}}}
\newcommand{\zero}{\CB{0}}

\DeclareMathOperator{\ls}{\normalfont\msf{span}}
\DeclareMathOperator{\cls}{\ol{\ls}}

\DeclareMathOperator*{\coplus}{\ol{\bigoplus}}

\DeclareMathOperator{\id}{\normalfont\msf{id}}

\DeclareMathOperator{\alg}{\normalfont\msf{alg}}
\renewcommand{\ker}{\operatorname{\msf{ker}}}
\renewcommand{\dim}{\operatorname{\msf{dim}}}

\newcommand{\C}{\bb{C}}

\newcommand{\E}{\bb{E}}

\newcommand{\N}{\bb{N}}

\newcommand{\R}{\bb{R}}

\newcommand{\Z}{\bb{Z}}
\newcommand{\cA}{\cmc{A}}
\newcommand{\cB}{\cmc{B}}
\newcommand{\cC}{\cmc{C}}
\newcommand{\cD}{\cmc{D}}
\newcommand{\cF}{\cmc{F}}
\newcommand{\cG}{\cmc{G}}
\newcommand{\cI}{\cmc{I}}
\newcommand{\cJ}{\cmc{J}}

\newcommand{\sB}{\eus{B}}

\newcommand{\sF}{\eus{F}}

\newcommand{\sK}{\eus{K}}
\newcommand{\sL}{\eus{L}}

\newcommand{\sN}{\eus{N}}

\newcommand{\U}{\mbf{1}}

\newcommand{\I}{{I\!\!\!\;I}}
\newcommand{\f}{\text{\scriptsize$\sF$}}

\ifthenelse{\boolean{nosectionequn}}{}{
    \numberwithin{equation}{section}
    }

\ifthenelse{\boolean{article}\or\boolean{letter}\or\boolean{nosectionequn}}{
    \setboolean{nosectionappendix}{true}}{}
\ifthenelse{\boolean{nosectionappendix}}{}{
    \renewcommand{\appendix}{
        \chapter*{\appendixname}
        \addcontentsline{toc}{chapter}{\appendixname}
        \renewcommand{\thesection}{\Alph{section}}
        \setcounter{section}{0}}}
   
\ifthenelse{\boolean{report}\or\boolean{book}}{
    }{}

\ifthenelse{\boolean{notheorem}}{}{
        \newcommand{\notename}{Note.}
        \newcommand{\mnname}{Mathematical note.}
        \newcommand{\enname}{End of the note.}
        \newcommand{\definame}{Definition.}
        \newcommand{\propname}{Proposition.}
        \newcommand{\lemname}{Lemma.}
        \newcommand{\exname}{Example.}
        \newcommand{\exername}{Exercise.}
        \newcommand{\remname}{Remark.}
        \newcommand{\obname}{Observation.}
        \newcommand{\thmname}{Theorem.}
        \newcommand{\corname}{Corollary.}
        \newcommand{\proofname}{Proof.}
    \ifthenelse{\boolean{german}}{
        \renewcommand{\mnname}{Mathematische Notiz.}
        \renewcommand{\enname}{Ende der Notiz.}
        \renewcommand{\exname}{Beispiel.}
        \renewcommand{\exername}{Übung.}
        \renewcommand{\remname}{Bemerkung.}
        \renewcommand{\obname}{Beobachtung.}
        \renewcommand{\thmname}{Satz.}
        \renewcommand{\corname}{Korollar.}
        \renewcommand{\proofname}{Beweis.}}{}
    \ifthenelse{\boolean{italian}}{
        \renewcommand{\mnname}{Nota matematica.}
        \renewcommand{\enname}{Fina della nota.}
        \renewcommand{\definame}{Definizione.}
        \renewcommand{\propname}{Proposizione.}
        \renewcommand{\exname}{Esempio.}
        \renewcommand{\exername}{Esercizio.}
        \renewcommand{\remname}{Nota.}
        \renewcommand{\obname}{Osservazione.}
        \renewcommand{\thmname}{Teorema.}
        \renewcommand{\corname}{Corollario.}
        \renewcommand{\proofname}{Dimostrazione.}

       \renewcommand{\appendixname}{Appendice}

       }{}
    \theoremheaderfont{\normalfont\bfseries}
    \theoremstyle{change}
        \theorembodyfont{\rmfamily}
            \newtheorem{emp}{}[section]
                \newcommand{\bemp}[1][]{
                    \begin{emp}\hskip-\labelsep\bf{#1}\hskip\labelsep}
                \newcommand{\eemp}{\end{emp}}
\newtheorem{itemp}[emp]{}
                \newcommand{\bitemp}[1][]{
                    \begin{itemp}\hskip-\labelsep\bf{#1}\hskip\labelsep\normalfont\itshape}
                \newcommand{\eitemp}{\end{itemp}}
            \newtheorem{note}[emp]{\notename}
                \newcommand{\bnote}{\begin{note}}
                \newcommand{\enote}{\end{note}}
            \newtheorem{mn}[emp]{\mnname}
                \newcommand{\bnm}{\begin{mn}~\begin{quotation}\renewcommand{\baselinestretch}{1}\small\noindent\ignorespaces}
                \newcommand{\enm}{\end{quotation}\hfill\bf{\enname}\end{mn}}
            \newtheorem{ex}[emp]{\exname}
                \newcommand{\bex}{\begin{ex}}
                \newcommand{\eex}{\end{ex}}
            \newtheorem{exer}[emp]{\exername}
                \newcommand{\bexer}{\begin{exer}}
                \newcommand{\eexer}{\end{exer}}
            \newtheorem{defi}[emp]{\definame}
                \newcommand{\bdefi}{\begin{defi}}
                \newcommand{\edefi}{\end{defi}}
            \newtheorem{rem}[emp]{\remname}
                \newcommand{\brem}{\begin{rem}}
                \newcommand{\erem}{\end{rem}}
            \newtheorem{ob}[emp]{\obname}
                \newcommand{\bob}{\begin{ob}}
                \newcommand{\eob}{\end{ob}}
        \theorembodyfont{\normalfont\itshape}
            \newtheorem{thm}[emp]{\thmname}
                \newcommand{\bthm}{\begin{thm}}
                \newcommand{\ethm}{\end{thm}}
            \newtheorem{prop}[emp]{\propname}
                \newcommand{\bprop}{\begin{prop}}
                \newcommand{\eprop}{\end{prop}}
            \newtheorem{cor}[emp]{\corname}
                \newcommand{\bcor}{\begin{cor}}
                \newcommand{\ecor}{\end{cor}}
            \newtheorem{lem}[emp]{\lemname}
                \newcommand{\blem}{\begin{lem}}
                \newcommand{\elem}{\end{lem}}
\newenvironment{empn}[1]{\lf\noindent\bf{#1}\ignorespaces\hskip\labelsep}{\lf}
		\newcommand{\bempn}[1]{\begin{empn}{#1}}
		\newcommand{\eempn}{\end{empn}}
		\newcommand{\bitempn}[1]{\begin{empn}{#1}\normalfont\itshape}
		\newcommand{\eitempn}{\end{empn}}
                \newcommand{\bnmn}{\begin{empn}{\mnname}~\begin{quotation}\renewcommand{\baselinestretch}{1}\small\noindent\ignorespaces}
                \newcommand{\enmn}{\end{quotation}\hfill\bf{\enname}\end{empn}}
		\newcommand{\bexn}{\begin{empn}{\exname}}
		\newcommand{\eexn}{\end{empn}}
		\newcommand{\bexern}{\begin{empn}{\exername}}
		\newcommand{\eexern}{\end{empn}}
		\newcommand{\bdefin}{\begin{empn}{\definame}}
		\newcommand{\edefin}{\end{empn}}
		\newcommand{\bremn}{\begin{empn}{\remname}}
		\newcommand{\eremn}{\end{empn}}
		\newcommand{\bobn}{\begin{empn}{\obname}}
		\newcommand{\eobn}{\end{empn}}

\newcommand{\qedsymbol}{~\rule[-0.35mm]{2mm}{2mm}}
    \newcounter{proof}[emp]
    \newenvironment{Proof}[1]{
        \vspace{1ex}
        \renewcommand{\item}[1][\stepcounter{proof}(\roman{proof})]%
            {##1\hskip\labelsep}
        \noindent\textsc{#1\hskip\labelsep}}{
        \nolinebreak\qedsymbol}
    \newcommand{\proof}[1][\proofname]{
        \begin{Proof}{#1}\ignorespaces}
    \newcommand{\qed}{\end{Proof}}
    \newcommand{\noqed}{
        \renewcommand{\qedsymbol}{}
        \end{Proof}}}
    \ifthenelse{\boolean{italian}}{
        \renewcommand{\proofname}{Dimostrazione.}}{}

\newcounter{OP}
		\newcommand{\bOP}{\stepcounter{OP}\begin{empn}{Open Problem \arabic{OP}:}}
		\newcommand{\eOP}{\end{empn}}

\usepackage[varg]{txfonts}

\newcommand{\ee}{\euf{e}}

\usepackage{stmaryrd}

\newcommand{\botimes}{\varogreaterthan}

\usepackage[matrix,arrow,curve]{xy}

\renewcommand{\thefootnote}{[\arabic{footnote}]}

\usepackage{hyperref}

\begin{document}

\bibliographystyle{amsalpha}

\title{Interacting Fock Spaces and Subproduct Systems}

\author{Malte Gerhold{\renewcommand{\thefootnote}{}{~}% \footnote{Keywords: Operator algebras; tensor and Cuntz-Pimsner-Toeplitz algebras; interacting Fock spaces; subproduct systems}
  }\!\!\!\;\thanks{MG acknowledges funding from the German Research Foundation (DFG) through the project "Non-Commutative Stochastic Independence: Algebraic and Analytic Aspects", project number 397960675.}  ~and Michael Skeide}

\date{~}

\maketitle

%\vspace{-3ex}
\begin{abstract}
\noindent
Interacting Fock spaces are the most general \nbd{\N_0}graded (pre-)Hilbert spaces with creation operators that have degree $1$ and generate everything out of a single vacuum vector $\Om$. It is the creators alone that generate the space out of the vacuum; so the same is true for the non-selfadjoint operator algebra generated by the creators. A formal definition has been given by Accardi, Lu, and Volovich (1997). Forthcoming work by Accardi and Skeide (2008), gave a different but equivalent definition, and also several desirable properties (embeddability, and what we are going to call here regularity, but also embeddability in Cuntz-Pimsner -Toeplitz type algebras) have been pointed out there.

In this paper we show that every interacting Fock space is embeddable, provided we ask the question the right way. This requires and motivates a new more flexible definition. (The definition does not allow for more interacting Fock spaces, but for more freedom how to capture their structure in a more useful way.) We show that the same statement for regularity must fail: There are irregular interacting Fock spaces -- irregular beyond repair; and they are quite natural. Embeddability allows to recover an interacting Fock space as a so-called \nbd{\vk}interacting Fock space. ($\vk$ is an operator on a usual full Fock space that allows to write the `interacting' creators $a^*(x)$ in terms of the usual creators $\ell^*(x)$ as $\vk\ell^*(x)$.) We show that interacting Fock spaces are classified by the $\vk$. We give criteria for when the creators of an interacting Fock space are bounded in general and under regularity. If all creators are bounded, then the Banach algebra and the \nbd{C^*}algebra generated by them, embed into the tensor algebras and the Cuntz-Pimsner-Toeplitz algebras, respectively, associated with several suitably chosen \nbd{C^*}correspondences.

We illustrate all this in the case of interacting Fock spaces coming from so-called subproduct systems, and determine for which $\vk$ the \nbd{\vk}interacting Fock space comes from a subproduct system. In the concluding multi-part section, we pose a number of problems for future work; for several of them we also illustrate why they do not possess uniform solutions, but solutions that strongly depend on the case.
\end{abstract}

\setcounter{footnote}{0}
\newpage
\section{Introduction} \label{intro}

After the detailed abstract, let us start right with the gist of the definition of \it{interacting Fock space} by Accardi, Lu, and Volovich \cite[Definition 18.1]{ALV97}:

\bdefi \label{ALVdefi}
Let $H$ be a (complex) vector space and form the \hl{tensor algebra} $\sF(H):=\Om\C\oplus\bigoplus_{n\in\N}H^{\otimes n}$ over $H$, where $\Om$ is some nonzero reference vector, the \hl{vacuum}. For each $x\in H$, define the \hl{creation operator} $\ell^*(x)$ on $\sF(H)$ by
\baln{
\ell^*(x)X_n
&
~:=~
x\otimes X_n
~~~
(X_n\in H^{\otimes n}, n\ge1),
&
\ell^*(x)\Om
&
~:=~x
}\ealn
(that is, $\ell^*(x)X= xX$, the product in the tensor algebra with unit $\Om$, for all $X\in\sF(H)$). Put $H^{\otimes 0}:=\Om\C$. Suppose on each $H^{\otimes n}$ ($n\in\N$) we have a semiinner product $(\bullet, \bullet)_n$ with kernel $\sN_n$ and put $(\Om,\Om)_0:=1$, so that $(\bullet,\bullet):=\bigoplus_{n\in\N_0}(\bullet, \bullet)_n$ is a semiinner product on $\sF(H)$ with kernel $\sN=\bigoplus_{n\in\N_0}\sN_n$. Put $H_n:=H^{\otimes n}/\sN_n$ and $\cI:=\sF(H)/\sN$. (Note that $\sN_0=\zero$, hence, $H_0=\Om\C$.) Then
\beqn{
\cI
~=~
\bigoplus_{n\in\N_0}H_n
}\eeqn
(we omit the simple proof; essentially $H_n\ni X_n+\sN_n=X_n+\sN\in\cI$ for $X_n\in H^{\otimes n}$). We say the pre-Hilbert space $\cI$ is an \hl{ALV-interacting Fock space} (denoting this situation by $\cI=(H,\bfam{(\bullet, \bullet)_n}_{n\in\N_0})$) if
\beqn{ \tag{$*$}\label{*}
H\otimes\sN_n
~\subset~
\sN_{n+1}
}\eeqn
(that is, $H\otimes\sN\subset\sN$), so that $a^*(x)\colon X+\sN\mapsto\ell^*(x)X+\sN$ well-defines the \hl{creation operators} $a^*(x)$ on $\cI$.
\edefi

\brem \label{ALVrem}
We collect some notes that should be mentioned but, otherwise (like all our remarks), should not interrupt the flow of reading.
\begin{enumerate}
\item \label{ALVr1}
The notion of \it{interacting Fock space} was motivated by an example due to Accardi and Lu \cite{AcLu92p,AcLu96}, emerging from QED. In this example, actually, the semiinner product is on the tensor algebra over a \nbd{\cB}bimodule (all tensor products over $\cB$), turning the quotient into a pre-Hilbert \nbd{\cB}module. In fact, one might study also these more general \hl{interacting Fock modules}. We emphasize that, here, we are concerned only with the scalar case. The full Fock module does occur, however, in its ``unperturbed'' form, when we discuss that the algebras generated by the creators on an interacting Fock space embed into Cuntz-Pimsner-Toeplitz algebras. (Actually, Skeide \cite{Ske98} shows that the interacting Fock module from \cite{AcLu96} \bf{is} a usual full Fock module, provided we choose the ``correct'' left module operation.) Still, it might be noteworthy that the papers \cite{AcLu92p,AcLu96} are likely to mark the first occurrence of full Fock modules even before Pimsner \cite{Pim97} and Speicher \cite{Spe98}, and that the three contexts are entirely different. 

\item \label{ALVr2}
The only true difference between Definition \ref{ALVdefi} and \cite[Definition 18.1]{ALV97} is (apart from some unnecessary requirements in the latter which are fulfilled automatically) that, here, we do not require that the creators possess a (formal) adjoint (in which case they are well-defined, automatically), but that we produce well-definedness by the kernel condition in \eqref{*}. In fact, as a minor side effect, in this paper we also free a number of results from Accardi and Skeide \cite{AcSk08} from the requirement that creators have adjoints.

\item
As the reader will have noticed, by the construction in Definition \ref{ALVdefi}, we are concerned with pre-Hilbert spaces, and both tensor products and direct sums are algebraic. Despite the fact that in the end we are interested basically in the case when the creators are bounded and that, therefore, we may and will complete the pre-Hilbert spaces in this case, in order not to lose any of the bounded(\bf{!}) examples it is indispensable to wait with this step until the last moment. In fact, even when all creators are bounded, certain operators that parametrize interacting Fock spaces, will remain unbounded  (see \ref{kappaunbex} and \ref{bA*unbLex}).

\item
The scope of the notion of \it{interacting Fock space} is to capture, in some sense, the most general situation of a Fock type pre-Hilbert space. What \it{in some sense} means, becomes clearer in a moment when we discuss the definition from \cite{AcSk08}. We do not claim that all spaces that are somehow related to Fock spaces are captured. (The GNS-spaces of \it{temperature states} on the \it{CCR-algebras} are not. Also Fock spaces from \it{species} discussed by Guta and Maassen \cite{GuMa02} are not. Actually, the latter would fit quite nicely into a description by Fock modules.) But we would not like to dispense with the properties that interacting Fock spaces possess.
\end{enumerate}
\erem

\noindent
The, in a sense, simplest class of interacting Fock spaces possible is captured in the following example by Accardi and Bozejko \cite{AcBo98}. Despite its striking simplicity, it captures to a surprisingly large extent rudimentary forms of the most important structure results on interacting Fock spaces; for this reason we repeat it here once more.

\bex \label{1-mex}
We consider the case $H=\C$, a so-called \hl{one-mode interacting Fock space}. So, $\sF(\C)=\bigoplus_{n\in\N_0}\C^{\otimes n}$ and we denote $e_0:=\Om$ and $e_n:=1^{\otimes n}$. A family of semiinner products is determined by the numbers $\ell_n=(e_n,e_n)_n\ge0$. For that the $(\bullet,\bullet)_n$ determine an interacting Fock space, we must have $\ell_0=1$ and $\ell_n=0$ $\Rightarrow$ $\ell_{n+1}=0$. These conditions are also sufficient. It follows that there are (unique, if $k_n=0$ $\Rightarrow$ $k_{n+1}=0$) numbers $k_n$ such that $\ell_n=k_n\ldots k_1$.

Suppose $\mu$ is a (nonzero) symmetric measure on the real line with finite moments of all orders. Then the orthogonal polynomials $P_n$ of $\mu$ satisfy and are determined by the following recursion
\baln{
P_0(t)
&
~=~
1,
&
P_1(t)
&
~=~
t,
&
tP_n(t)
&
~=~
P_{n+1}(t)+k_nP_{n-1}(t)
~~~
(n\ge1),
}\ealn
for unique (positive) numbers $k_n$. (If $\mu$ is not symmetric, then on the right-hand side of the recursion there is also a term proportional to $P_n$. \cite{AcBo98} take into account also this case; here, we ignore it.)

Since $\int P_m(t)P_n(t)\mu(dt)=\delta_{m,n}\ell_n$ and since the $P_n$ are real, it follows that $e_n+\sN_n\mapsto P_n$ defines an isometry from $\cI$ onto $\ls P_{\N_0}\subset L^2(\mu)$. The creation operator $a^*:=a^*(1)\colon e_n+\sN_n\mapsto e_{n+1}+\sN_{n+1}$ has an adjoint $(a^*)^*=:a\colon  e_n+\sN_n\mapsto (e_{n-1}+\sN_{n-1})k_n$ (with $e_{-1}:=0$), and the crucial observation in \cite{AcBo98} is that $a^*+a$ (for symmetric $\mu$), under the isomorphism $\cI\rightarrow \ls P_{\N_0}$, acts as
\beqn{
(a^*+a)(e_n+\sN_n)
~\longmapsto~
tP_n,
}\eeqn
that is, $a^*+a$, on $\ls P_{\N_0}$, acts as multiplication with the function $t$. In the context of this paper, we are more interested in the following fact. Suppose we equip $\sF(\C)$ with the canonical inner product where the $e_n$ are orthonormal. Then we may embed $\cI$ into $\sF(\C)$ via the (adjointable) isometry $\xi\colon e_n+\sN_n\mapsto e_n\sqrt{\ell_n}$ and we find
\beqn{ \tag{$**$}\label{**}
\xi a^*(x)\xi^*
~=~
\vk\ell^*(x),
}\eeqn
where $\vk$ is some square root of the operator $k\colon e_n\mapsto e_nk_n$ on $\sF(\C)$. In fact, one of the main results of this paper is (see Theorem \ref{**thm}) that every interacting Fock space $\cI$ can be recovered as $\xi\cI\subset\ol{\sF(H)}$ for an isometry $\xi$ from $\cI$ to the completion of the usual algebraic Fock space $\sF(H)$ over a suitable pre-Hilbert space $H$, $\ol{\sF(H)}$, in such a way that the creators have the form as in \eqref{**}. Moreover, the occurring $\vk$ suitably parametrize interacting Fock spaces.
\eex

After this example, where $H$ is a pre-Hilbert space, we return (really only for a moment) to the situation in Definition \ref{ALVdefi} where $H$ is just a vector space. An ALV-interacting Fock space $\cI=(H,\bfam{(\bullet, \bullet)_n}_{n\in\N_0})$ comes shipped with the \hl{creator map} $a^*\colon H\rightarrow\sL(\cI)$ from $H$ into the linear operators on $\cI$, which is linear and satisfies
\beqn{ \tag{$*\!*\!*$}\label{***}
\ls a^*(H)H_n
~=~
H_{n+1}.
}\eeqn
This means, in particular, that everything in $\cI$ is created out of the vacuum $\Om$ by successive application of creation operators $a^*(x)$. In the definition by Accardi and Skeide \cite{AcSk08}, emphasis is put on the family of pre-Hilbert spaces $H_n$ and the creator map $a^*$. A formulation of \cite[Definition 2.2]{AcSk08} that matches the situation of Definition \ref{ALVdefi} is:

\bdefi \label{ASdefi}
Let $\bfam{H_n}_{n\in\N_0}$ be a family of pre-Hilbert spaces where $H_0=\Om\C$ for a fixed unit vector $\Om$, the \hl{vacuum}, and put $\cI:=\bigoplus_{n\in\N_0}H_n$. Let $H$ be a vector space and suppose $a^*\colon H\rightarrow\sL(\cI)$, the \hl{creator map}, is a linear map satisfying \eqref{***}. Then $\cI$ is an \hl{interacting Fock space} based on $H$ (denoted as $\cI=(\bfam{H_n}_{n\in\N_0},a^*)$).
\edefi

Let us convince ourselves that Definitions \ref{ALVdefi} and \ref{ASdefi} speak about ``the same'' thing. (This has not been clarified that explicitly in \cite{AcSk08}.)
\begin{itemize}
\item
We know already that every ALV-interacting Fock space $\cI=(H,\bfam{(\bullet, \bullet)_n}_{n\in\N_0})$ is an interacting Fock space based on $H$ via $H_n=H^{\otimes n}/\sN_n$ and $a^*\colon x\mapsto a^*(x)$ (obviously, by definition, having the same creators $a^*(x)$).

\item
Every interacting Fock space $\cI=(\bfam{H_n}_{n\in\N_0},a^*)$ based on $H$, comes along with a linear surjective map $\Lambda:=\bigoplus_{n\in\N_0}\Lambda_n\colon\sF(H)\rightarrow\cI$ where $\Lambda_n\in\sL(H^{\otimes n},H_n)$ is defined by
\beqn{ \tag{$*\!*\!*\hspace{.6pt}*$}\label{****}
\Lambda_n
\colon
x_n\otimes\ldots\otimes x_1
~\longmapsto~
a^*(x_n)\ldots a^*(x_1)\Om
}\eeqn
and $\Lambda_0\colon\Om\mapsto\Om$. Then for the semiinner products $(\bullet,\bullet)_n:=\AB{\Lambda_n\bullet,\Lambda_n\bullet}$ on $H^{\otimes n}$, the map $\Lambda_n(x_n\otimes\ldots\otimes x_1)\mapsto x_n\otimes\ldots\otimes x_1+\sN_n$ establishes a unitary $H_n\rightarrow H^{\otimes n}+\sN_n $. Moreover, from
\beqn{
\Lambda_{n+1}(\ell^*(x)X_n)
~=~
a^*(x)\Lambda_nX_n,
}\eeqn
it follows that the semiinner products fulfill \eqref{*} (and $(\Om,\Om)_0=1$) and that, under the stated isomorphism, the ALV-interacting Fock space $(H,\bfam{(\bullet, \bullet)_n}_{n\in\N_0})$ has the same creators as $\cI=(\bfam{H_n}_{n\in\N_0},a^*)$.
\end{itemize}
Note that an ALV-interacting Fock space, with the structures defined in the first part of this discussion, \bf{is} an interacting Fock space based on $H$; and, as a \hl{convention}, we will always consider it as such. On the other hand, the isomorphism that identifies in the second part an interacting Fock space as an ALV-interacting Fock space cannot be discussed away; for an interacting Fock space being ALV is an extra information that tells how the interacting Fock space has been obtained; if the latter is not there, it cannot be recovered better than up to (canonical) isomorphism.

We stated Definition \ref{ASdefi} for $H$ just a vector space in order to be compatible with Definition \ref{ALVdefi}. As compared with \cite[Definition 2.2]{AcSk08} (where $H$ is required to be a pre-Hilbert space), in Definition \ref{ASdefi} (like in Definition \ref{ALVdefi}; see Remark \ref{ALVrem}\eqref{ALVr2}), we also have removed the condition that the $a^*(x)$ be adjointable. In either case, we will speak of an \hl{adjointable} interacting Fock space if all creators have an adjoint. However, in applications $H$ is (almost) always a pre-Hilbert space; so, from now on, as a \hl{convention}, we shall always assume (adding to Definitions \ref{ALVdefi} and \ref{ASdefi}) that $H$ is a pre-Hilbert space. This means, $\sF(H)$ does already possess an inner product $\AB{\bullet,\bullet}$ arising from tensor product and direct sum of pre-Hilbert spaces. In other words, $\sF(H)$ is not just the tensor algebra over the vector space $H$, but the (\hl{algebraic}) \hl{full Fock space} over the pre-Hilbert space $H$.

As already illustrated in Example \ref{1-mex}, the interplay between the semiinner product $(\bullet,\bullet)$ on $\sF(H)$ and the inner product $\AB{\bullet,\bullet}$ on $\sF(H)$ plays a very important role.

\bulletline
Let us briefly describe what we are up to in the following sections, adding also more motivation.

\newpage
% \lf\lf
When interacting Fock spaces are obtained by introducing a semiinner product $(\bullet,\bullet)$ on $\sF(H)$, then in almost all examples of this type the new semiinner product is obtained from the original inner product $\AB{\bullet,\bullet}$ of the pre-Hilbert space $\sF(H)$ with the help of a positive operator $L=\bigoplus_{n\in\N_0}L_n$ (with $L_0\Om=\Om$) as $(\bullet,\bullet):=\AB{\bullet,L\bullet}$. (If we would speak about interacting Fock \it{modules}, then the QED-example from \cite{AcLu92p,AcLu96} mentioned in Remark \ref{ALVrem}\eqref{ALVr1} would be a prominent exception; here the new semiinner product is computed directly and we do not know if it can be induced by positive type operator on the full Fock module.) We discuss such \it{positive operator induced} or \it{POI-interacting Fock spaces} in Section \ref{poiSEC}. In particular, we push forward to the more general Definitions \ref{ALVdefi} and \ref{ASdefi} the result from \cite{AcSk08} that POI-interacting Fock spaces are precisely those ALV-interacting Fock spaces that are \it{regular} interacting Fock spaces based on a pre-Hilbert space $H$ in the sense that the canonical surjection $\Lambda\colon\sF(H)\rightarrow\cI$ defined by \eqref{****} has an adjoint $\cI\rightarrow\ol{\sF(H)}$.

In Theorem \ref{nonregthm}, we show that a large class of interacting Fock spaces is non-regular. While the \it{non-nilpotent full} interacting Fock spaces occurring in Theorem \ref{nonregthm} have the defect to be not \it{vacuum- separated} (see Section \ref{classSS}), in the course of this paper we meet several examples of vacuum-separated interacting Fock spaces that are not regular. A large class of (regular) examples that arises from so-called \it{subproduct systems} directly as interacting Fock spaces based on a Hilbert space, we discuss in Section \ref{piSEC}. Despite the late occurrence, this class motivated large parts of this paper starting from Section \ref{boundSEC}. In Section \ref{pivesSS}, we generalize this scheme to \it{nondegenerate productive systems}.

\lf
Both definitions, ALV-interacting Fock spaces (with its subclass of POI-interacting Fock spaces) and interacting Fock spaces based on a pre-Hilbert space (with its subclass of regular interacting Fock spaces), are relative to a chosen pre-Hilbert space $H$. For several reasons it is indispensable to come up with yet another (new) definition of (\it{abstract}) \it{interacting Fock space} (see Definition \ref{IFSdefi}) that abandons this dependence on $H$. For instance, already \cite{AcSk08} pointed out \it{embeddability}, that is, existence of an \it{even}, \it{vacuum-preserving} isometry $\xi\colon\cI\rightarrow\ol{\sF(H)}$, as a crucial property, which an interacting Fock space based on $H$ may possess or not. (We discuss the important consequences of embeddability in Section \ref{classSEC}.) Definition \ref{IFSdefi} allows to \bf{choose} the dependence on $H$ appropriately so that the resulting interacting Fock space based on some $H$ \bf{is} embeddable; Theorem \ref{aIFSembthm}. Moreover, Theorem \ref{nonregthm} shows that the same attempt to ``repair'' non-regulatirity does not work; there are (abstract) interacting Fock spaces that can in no way be \it{regularly based} on a pre-Hilbert space.

Just to give the idea: The pre-Hilbert space $H$ parametrizes the set $A^*:=a^*(H)$ of all creation operators by means of the creator function $a^*$. Definition \ref{IFSdefi} substitutes in Definition \ref{ASdefi} the pre-Hilbert space $H$ and the creator map $a^*$ by the vector subspace $A^*$ of $\sL(\cI)$, and the cyclicity condition in \eqref{***} with $\ls A^*H_n=H_{n+1}$. We \it{may} base such $\cI$ and $A^*$ by \it{choosing} a linear surjective map $a^*\colon H\rightarrow A^*$. This definition also meets perfectly the frequent situation in quantum probability and operator algebras, where the main object of interest is the \hl{algebra} $\cA:=\alg A^*$ (or the \hl{\nbd{*}algebra} $\cA^*:=\alg^*A^*$) generated by $A^*$ (plus, in some contexts, the \hl{vacuum state} $\AB{\Om,\bullet\Om}$). Two paragraphs below we briefly mention bounded example classes.

\vspace{1ex}
In Section \ref{classSEC}, we examine the important consequences of embeddability. Already \cite{AcSk08} observed that under embeddability, creators can be written as in \eqref{**} for a suitable (even, vacuum-preserving) operator $\vk$ (that goes into the dense subspace $\xi\cI\oplus(\xi\cI)^\perp$ of $\ol{\sF(H)}$ with another dense subspace $(H\otimes(\xi\cI\oplus(\xi\cI)^\perp))\oplus\Om\C$ as domain). Apart from the (minor) effort to push this forward to the situation of Definition \ref{ASdefi} (dropping adjointability), we show two major results. Firstly, the operator $\vk$ is uniquely determined by \eqref{**} (and $\xi$) and satisfies two extra conditions; we shall call $\vk$ satisfying these conditions a \it{squeezing}. Moreover, varying the squeezing $\vk$ and the subspaces corresponding to $\xi\cI\cong\cI$, we get a parametrization (up to suitable isomorphism) of all interacting Fock spaces based on $H$ in terms of $\vk$ (Theorem \ref{kIFSthm}), so-called \it{\nbd{\vk}interacting Fock spaces} (Definition \ref{kdefi}). Secondly, taking into account Theorem \ref{aIFSembthm}, which asserts that all interacting Fock spaces \bf{can} be embeddably based, we get, in Theorem \ref{kIFSisothm}, that every interacting Fock space is (isomorphic to) a \nbd{\vk}interacting Fock space. Whether $\vk$ is bounded or not, is a property intrinsic to $\cI$; see Proposition \ref{kpartprop}.

\vspace{1ex}
While in Section \ref{classSEC} we examine the consequences of embeddability, in Section \ref{boundSEC} we examine the consequences of \it{boundedness}, whereas, in Section \ref{bcritSEC} we give criteria when an interacting Fock space, actually, is bounded. We say an interacting Fock space is \it{bounded} if $A^*$ has only bounded elements. In this case, $A^*$ generates a (\nbd{*})algebra $\cA^{(*)}$ of bounded operators that may be completed to obtain a Banach (a \nbd{C^*})algebra of bounded operators acting on the completion $\,\ol{\!\cI}$ of $\cI$.  (For instance, Davidson, Ramsey, and Shalit \cite{DRS11} showed that  the interacting Fock spaces associated with commutative finite-dimensional subproduct systems are classified up to isomorphism by the isomorphism classes of the (non-selfadjoint) operator algebras $\ol{\cA}$ and that this fails for the selfadjoint operator algebras $\ol{\cA^*}$. Kakariadis and Shalit \cite{KaSha15p} have done the general case. See Sections \ref{piSEC} and \ref{boundSS}.)  In Section \ref{boundSEC}, we basically explain in the bounded case how these operator (\nbd{*})algebras embed into \it{tensor algebras} in the sense of Muhly and Solel \cite{MuSo98} (\it{Cuntz-Pimsner-Toeplitz algebras} in the sense of Pimsner \cite{Pim97}); actually, we present three (potentially) different ways to do that. In the adjointable but not necessarily bounded case, this has been done in \cite{AcSk08}, while \cite{KaSha15p} have done it for subproduct systems. Here, it is really important to free, in the bounded case, the situation from the hypothesis of adjointability. (See the end of this section for the notions of adjointability we are using in this paper.)

\vspace{1ex}
In Section \ref{bcritSEC} we give criteria for when an interacting Fock space is bounded. For an abstract one, which is characterized by the set $A^*$ of creators, there is not really more that can be done, other than just look at the elements of $A^*$ and check if they are bounded. So, Section \ref{bcritSEC} is on interacting Fock spaces that are based or -- better -- embeddably or even regularly based on some pre-Hilbert space. In this case, we have available one or more of the parameters $\Lambda,\lambda,\vk,L$, and the criteria we give are in terms of these parameters. Among the results there are: Boundedness of all $a^*(x)$ (even boundedness of the creator map $a^*$) does not imply regularity; boundedness of $\vk$ is sufficient, but not necessary; boundedness of $L$ is neither sufficient, nor necessary. The necessary and sufficient criterion that all $a^*(x)$ are bounded, given by the following (unbounded operator) inequality
\vspace{-1ex}
\beqn{
\ell(x)L\ell^*(x)
~\le~
M_xL
}\eeqn
for all $x\in H$ (together with an analogous criterion for boundedness of the creator map), answers the long standing question when POI-interacting Fock spaces have bounded creators (creator maps).

\lf
In Section \ref{piSEC}, we finally pass to \it{subproduct systems} (over discrete time $n\in\N_0$, to be precise), and show that their Fock spaces (Shalit and Solel \cite{ShaSo09}) are (completions of) interacting Fock spaces. We show that they are \nbd{\vk}interacting Fock spaces whose squeezings $\vk$ are projections fulfilling an extra condition. As a by-product we determine the structure of all squeezings that are projections.

\lf
In Sections \ref{poiSEC} -- \ref{piSEC}, we put up a framework in which to deal with interacting Fock spaces. We present a general theory with many new results and an even larger number of examples and counter examples that illustrate the difficulties of the general theory. Only rarely, we formulate the results presented in counter examples in the form of no-go theorems, but evidently they make up a prominent part -- maybe, even more important than the positive results -- of this paper. In Section \ref{olSEC}, we concisely formulate and illustrate a number of problems and the questions they trigger for future work, putting into evidence the richness of the theory. Each of the problems is -- in some cases illustrated by tricky preliminary results and no-go statements or toy examples -- too ambitious to be tackled in this, already quite lengthy, paper.  Section \ref{olSEC} is subdivided into six parts.

In \ref{regSS}, we ask for criteria for regularity of interacting Fock spaces. Emphasis (though not exclusively) is put on \nbd{\vk}interacting Fock spaces (that is, embedded ones). We show that there is no direct relation between regularity of a \nbd{\vk}interacting Fock space and adjointability of $\vk$. (Corollary \ref{padcor}.)

In \ref{boundSS}, we explain very briefly what is relevant to us from the works \cite{ShaSo09,DRS11,KaSha15p}. We formulate which classification problems this suggests for interacting Fock spaces or for distinguished subclasses of interacting Fock spaces different from those considered in Section \ref{piSEC}.

In \ref{pivesSS}, based on the notion of \it{productive system} from Shalit and Skeide \cite[Section 6]{ShaSk10p}, we generalize the relation between subproduct systems and \nbd{\vk}interacting Fock spaces to (contractive) \it{non-degenerate productive systems}. In Theorem \ref{prodSkthm}, we prove a criterion for when a squeezing $\vk$ gives rise to the interacting Fock space of such a productive system. The criterion is much more complicated than that for subproduct systems in Section \ref{piSEC}. We leave entirely open the question how to apply this criterion to a given $\vk$. Calling our interacting Fock spaces \it{left}, we put into evidence the intimate relation of productive systems to interacting Fock spaces that are \it{left} and \it{right}. (We wonder, how Voiculescu's \it{bi-freeness} might fit into this context.)

In \ref{classSS}, we give in to the idea that further progress in classification of interacting Fock spaces will depend on finding ``good'' extra conditions they should fulfill. Examples for such conditions are being \it{injectively based} (the creator function $a^*$ is injective, hence, bijective onto $A^*$) or \it{vacuum-separated} ($a^*\Om=0$ implies $a^*=0$ for all $a^*\in A^*$; a property that does not depend on a possible basing). Emphasis is put on the question whether this might help to decide on regularity. For instance, in Proposition \ref{reg-inregprop} we show that restricting to injective basings does not change the answer to the question of regularity. The following case study in \ref{2FockSS} shows that the members of a very simple class of vacuum-separated interacting Fock spaces possess different answers to the question of regularity.

In \ref{2FockSS}, we propose the class of \it{interacting Fock \nbd{n}spaces} (distinguished by being vacuum-separated and having one-dimensional $H_n$ and $H_{n+1}=\zero$), and we completely determine their structure for $n=2$. So, we look at spaces $\cI=\Om\C\oplus H_1\oplus\Om_2\C$. We show: An interacting Fock \nbd{2}space may be irregular (Example \ref{Phinobex}); it is always regular provided $H_1$ has a countable Hamel basis (Theorem \ref{Phicbthm}).

In \ref{autoSS}, we augment the number of notions of isomorphism from $2$ to $4$, and we phrase the natural questions regarding the corresponding automorphism groups of interacting Fock spaces and how the former classify the latter. This adds to the questions about classification in terms of the associated operators algebras from Sections \ref{boundSEC} and \ref{boundSS}.

\lf\noindent
\bf{Notation.~}
For the discussion of interacting Fock spaces we need to work with pre-Hilbert spaces. Direct sums and tensor products are understood algebraically.%
\footnote{
For at least two reasons, this is not exaggerated generality, but necessary and unavoidable flexibility. Firstly, it actually lightens notation quite a bit when we discuss spaces where, like the Boson Fock space (this is Example \ref{POIex}\eqref{POI1} for $q=1$), the creation operators are unbounded; and we, surely, would not be happy to exclude classical examples like the Boson Fock space from the discussion. Secondly, yes, in the end we are interested in the case of bounded creation operators as they occur, for instance, from subproduct systems, and will complete the interacting Fock spaces; but, as our Theorems \ref{a*iffthm} and \ref{a*Lthm} show, it is not possible to characterize efficiently interacting Fock spaces with bounded creators by just bounded parameters $\vk$ or $L$. (This resembles a bit the characterization of morphisms of time ordered product systems from Barreto, Bhat, Liebscher, and Skeide \cite[Secion 5.2]{BBLS04}, which is quite a bit easier than the characterization of \it{bounded} morphisms; see Bhat \cite[Section 6]{Bha01}.) And we do not wish to lose these cases.
}
Consequently, we need the following spaces of operators. The space $\sL(H,H')$ of linear maps from the pre-Hilbert space $H$ to the pre-Hilbert space $H'$. Of course, here and in a similar way for all other spaces of operators, in the case $H'=H$ we will write $\sL(H)$. The space $\sL^a(H,H')$ contains those elements of $\sL(H,H')$ that have an adjoint in $\sL(H',H)$. We do \bf{not} assume that an adjoint has maximal domain. (For $a\in\sL^a(H,H')$, the domain of $a^*$ is $H'$ and $H'$ is mapped by $a^*$ into its codomain $H$.) We use the letter $\sB$ to indicate the bounded parts of these spaces. For fixed $H,H'$, an operator $a\in\sL(H,\ol{H'})$ is called \hl{weakly adjointable} if there exists $a^*\in\sL^a(H',\ol{H})$ such that $\AB{ah,h'}=\AB{h,a^*h'}$ for all $h\in H,h'\in H'$. (Note: \it{Weakly adjointable} is a notion relative to two chosen pre-Hilbert spaces $H$ and $H'$. So if $G=\ol{H'}$, being weakly adjointable as element of $\sL(H,G)$ is a different thing from being weakly adjointable as element of $\sL(H,\ol{H'})$. Obviously, the notion also applies to an operator $a\in\sL(H,H')$, considering it a map $H\rightarrow H'\subset\ol{H'}$.) If $a$ is weakly adjointable, then $(a^*)^*=a$ considered as an element in  $\sL(H,\ol{H'})$ is a weak adjoint of $a^*$. If $a$ is weakly adjointable, then it is a closeable densely defined operator $\ol{H}\rightarrow\ol{H'}$ in the usual sense (see also Section \ref{regSS}), with core $H$. In particular, if $H$ is a Hilbert space, then a weakly adjointable $a$ is bounded. A (weakly) adjointable operator $a$ is (\hl{weakly}) \hl{self-adjoint} if $a^*=a$. (Note: Weakly self-adjoint coincides with the usual definition of \it{symmetric} in functional analysis.)

\lf\noindent
\bf{Pre-Fock notation.~} All Fock-type spaces -- in this paper and elsewhere -- are in the first place \it{graded} vector spaces. This makes available the notion of linear maps with a \it{degree} in $\Z$. It is useful to do this once for all, and introduce a unified way to address these structures. However, Fock-type spaces are more than just graded vector spaces, but have an important specialty about them: The vacuum; that is, a grade-zero space of a particular form. A \hl{pre-Fock space} $\cI$ is a(n \nbd{\N_0})graded vector space, that is, $\cI=\bigoplus_{n\in\N_0}H_n$ for vector spaces $H_n$, with a distinguished non-zero vector $0\ne\Om\in\cI$, the \hl{vacuum}, such that $H_0=\Om\C$. We sometimes write $\cI=(\bfam{H_n}_{n\in\N_0},\Om)$.

For every pair $\cI=\bigoplus_{n\in\N_0}H_n$ and $\cJ=\bigoplus_{n\in\N_0}G_n$ of graded vector spaces, a linear map $a\in\sL(\cI,\cJ)$ has \hl{degree} $n\in\Z$ if $aH_m\subset G_{m+n}$ for all $m\in\N_0$ (where we use the conventions that $G_k=H_k=\zero$ for $k\le-1$). We denote the set of all linear maps from $\cI$ to $\cJ$ that have degree $n$ by $\sL_{(n)}(\cI,\cJ)$. The elements of $\sL_{(0)}(\cI,\cJ)$ are called \hl{even}.

We say, an even map $a$ from a pre-Fock space $\cI=(\bfam{H_n}_{n\in\N_0},\Om)$ to a pre-Fock space $\cJ=(\bfam{G_n}_{n\in\N_0},\Om')$ is a \hl{Fock map} if it is \hl{vacuum-preserving}, that is, if $a\Om=\Om'$.

If $\cI$ and $\cJ$ are direct sums of pre-Hilbert spaces, then $\Om$ and $\Om'$ will be required to be unit vectors. Moreover, a vacuum-preserving map $a\colon\cI\rightarrow\ol{\cJ}$ will also be called a \hl{Fock map} if $aH_n\subset\ol{G_n}$, that is, if $a$ is a Fock map into the algebraic direct sum over the completions $\ol{G_n}$.

We, usually, will use the \bf{same symbol} $\Om$ (without varying it to something like $\Om'$) for all occurring pre-Fock spaces, so that \it{vacuum-preserving} for a Fock map $a$, really, means $a\Om=\Om$.

For Section \ref{boundSEC}, only: We use analogous terminology in the category of right modules (bimodules) over a fixed unital algebra $\cB$ with the variation that $H_0=\Om\cB$ is required to be isomorphic as right module (as bimodule) to $\cB$ via $\Om\mapsto\U$.

\newpage

\section{POI-Interacting Fock spaces} \label{poiSEC}

\bdefi \label{posdefi}
Let $H$ be a pre-Hilbert space. An operator $L\in\sL(H,\ol{H})$ is \hl{positive} (writing $L\ge0$) if $\AB{x,Lx}\ge0$ for all $x\in H$.
\edefi

\begin{itemize}
\item
We might have called this \it{weakly positive}, reserving \it{positive} for operators in $\sL(H)$. We opted not to do so, and will refer to the latter situation as a \hl{positive operator on $H$}.

\item
A positive operator is weakly selfadjoint. A positive operator on $H$ is selfadjoint.

\item
Positivity induces a partial order among operators in $\sL(H,\ol{H})$ by defining $L\ge L'$ if $L-L'\ge0$.
\end{itemize}

\bdefi \label{POIdefi}
An ALV-interacting Fock space $\cI=(H,\bfam{(\bullet, \bullet)_n}_{n\in\N_0})$ (according to Definition \ref{ALVdefi} and, by our convention, with $H$ being a pre-Hilbert space) is a \hl{positive operator induced} or \hl{POI-interacting Fock space} if there is a Fock map $L\in\sL(\sF(H),\ol{\sF(H)})$ such that $(\bullet,\bullet)=\AB{\bullet,L\bullet}$.
\edefi

\begin{itemize}
\item
Recall the pre-Fock notations from the end of Section \ref{intro}: $L$ being a Fock map means that $L$ goes into the pre-Fock space $\bigoplus_{n\in\N_0}\ol{H^{\otimes n}}\subset\ol{\sF(H)}$ and that, as such, it is even and vacuum-preserving. That is, $L=\bigoplus_{n\in\N_0}L_n$ and $L_0=\id_{H_0}$.

\item
A Fock map $L$ induces a POI-interacting Fock space via  $(\bullet,\bullet):=\AB{\bullet,L\bullet}$ if and only if $L\ge0$ (that is, $L_n\ge0$ for all $n$) and $H\otimes\ker L\subset\ker L$ (that is, $H\otimes\ker L_n\subset\ker L_{n+1}$). The latter follows from $\sN=\ker L$. (Indeed, if $X\in\ker L$, then $(X,X)=\AB{X,LX}=0$, so $X\in\sN$. Conversely, if $X\in\sN$, so that $(X,X)=0$, then, by Cauchy-Schwartz inequality, $0=(Y,X)=\AB{Y,LX}$ for all $Y\in\sF(H)$. Since, $\sF(H)$ is dense, $\ol{\sF(H)}\ni LX=0$, so $X\in\ker L$.) 
\end{itemize}
Typical classes of examples are:
\bex \label{POIex}
\begin{enumerate}
\item \label{POI1}
By setting
\beqn{
L_n
\colon
x_n\otimes\ldots\otimes x_1
~\longmapsto~
\sum_{\sigma\in S_n}x_{\sigma(n)}\otimes\ldots\otimes x_{\sigma(1)}q^{\mrm{inv}(\sigma)},
}\eeqn
for $q\in\SB{-1,1}$ ($\mrm{inv}(\sigma)$ being the number of \it{inversions} of the permutation $\sigma\in S_n$), we get Bozejko's and Speicher's \cite{BoSp91a} \nbd{q}Fock space, whose creators and their adjoints satisfy the \nbd{q}commutation relations
\beqn{
a(x)a^*(y)-qa^*(y)a(x)
~=~
\AB{x,y}.
}\eeqn
The case $q=0$ (hence, $L=\id_{\sF(H)}$) is just the full Fock space. The cases $q=1$ and $q=-1$ give the \it{Boson} (or \it{symmetric}) and the \it{Fermion Fock space}, respectively. While in these extreme cases the $\frac{L_n}{n!}$ are projections and, therefore, easily established to be positive, in the general case $0<\abs{q}<1$ showing positivity is a tough problem.

\item \label{POI2}
A large class of examples, so-called \it{standard interacting Fock spaces} \cite{ALV97}, arises from $H\subset L^2(M,\mu)$ (usually, referred to as \it{test function space}) for some (\nbd{\sigma}finite) measure space $(M,\mu)$ and $L_n$ given by multiplication of the elements in $H^{\otimes n}\subset L^2(M^n,\mu^{\otimes n})$ with (measurable) positive functions on $M^n$. Standard interacting Fock spaces have been examined in particular by Lu and his coworkers. For instance, multiplying with the indicator function of the set $\CB{\alpha_n\ge\ldots\ge\alpha_1\ge 0}$ on $\R^n$, gives rise to the \it{time-ordered} or \it{chronological} or \it{monotone Fock space} examined first as interacting Fock space by Lu and Ruggieri \cite{LuRu98p}. No nontrivial symmetric Fock space is standard.
\end{enumerate}

\lf\noindent
Note that the first class has operators $L_n$ that map into $H^{\otimes n}$, while in the second class (unless in very special cases) we will need completion.
\eex

\brem
Note that Example \ref{1-mex} is a standard interacting Fock space. In fact, $\C$ is the $L^2$ of a probability measure concentrated in a single point. What we did in that example, can be generalized to standard interacting Fock spaces. So, let $L_n$ be positive measurable functions on $M^n$ acting as multiplication operators on $L^2(M^n,\mu^{\otimes n})$ in such a way that for the dense subspace $H$ of $L^2(M,\mu)$ each $H^{\otimes n}$ is in the natural domain of $L_n$. By the \it{Radon-Nikodym theorem}, the kernel condition on the $L_n$ is satisfied (if and) only if there are positive (``\nbd{L}almost surely'' unique) measurable functions $K_n$ on $M^n$ such that $L_{n+1}(t_{n+1},t_n,\ldots,t_1)=K_{n+1}(t_{n+1},t_n,\ldots,t_1)L_n(t_n,\ldots,t_1)$, almost surely. In terms of operators, this reads
\vspace{-1ex}
\beqn{
L_{n+1}
~=~
K_{n+1}(\id_H\otimes L_n),
\vspace{-2ex}
}\eeqn
so that
\vspace{-2ex}
\beqn{
L_n
~=~
K_n(\id_H\otimes K_{n-1})\ldots(\id_{H^{\otimes(n-1)}}\otimes K_1)
}\eeqn
for all $n\in\N$ (with initial condition $L_0=\id_{\Om\C}$). Also here, $\xi:=\sqrt{L}$, considered as operator $\cI\rightarrow\ol{\sF(L^2(M,\mu))}$, defines an isometry fulfilling \eqref{**}, where (modulo adjusting domain and codomain appropriately) $\vk=\sqrt{K}$.

Applying brute force linear algebra to the kernel condition $\ker L_{n+1}\subset H\otimes\ker L_n$ (see \cite[Lemma 5.4]{AcSk08}), also for a general POI-interacting Fock space there exist $K_{n+1}\in\sL(H\otimes\ol{H^{\otimes n}},\ol{H^{\otimes(n+1)}})$ such that $L_n$ is given by the preceding recursion. The recursion, yes, does capture entirely the kernel condition, by expressing the $L_n$ in terms of the $K_n$. However, if $L_{n+1}$ and $\id_H\otimes L_n$ do not commute, it leaves completely out of control the question for which $K_n$ the preceding sequence would consist of positive operators. We come back to this problem (and resolve it) in Section \ref{classSEC}.
\erem

It is natural to ask, if all ALV-interacting Fock spaces are POI (answer no), and (if not) how they are distinguished. We, first, answer the second question.

\blem \label{adj=poslem}
Let $H$ be a pre-Hilbert space. For another semiinner product $(\bullet,\bullet)$ on $H$, put $H_\sN:=H/\sN$ where $\sN:=\ker(\bullet,\bullet)$. Define the quotient map $\Lambda\colon x\mapsto x+\sN$. Then $(\bullet,\bullet)=\AB{\bullet,L\bullet}$ for some positive operator $L\in\sL(H,\ol{H})$ if and only if $\Lambda$ has an adjoint in $\sL(H_\sN,\ol{H})$.
\elem

\proof
If $\Lambda$ has a (weak!) adjoint, then $L:=\Lambda^*\Lambda$ is the positive operator we seek. Conversely, if $\AB{\Lambda x,\Lambda y}=(x,y)=\AB{x,Ly}$, then for each $z=\Lambda y\in H_\sN$ ($\Lambda$ is surjective!), the linear functional $x\mapsto\ol{\AB{\Lambda x,z}}$ on $H$ is bounded by $\norm{Ly}$, so that there is a unique element in $\ol{H}$, denoted by $\Lambda^*z$, such that $\AB{x,\Lambda^*z}=\AB{\Lambda x,z}$. The map $\Lambda^*\colon z\mapsto\Lambda^*z$ is a (weak!) adjoint of $\Lambda$.\qed

\bcor \label{Lam*ALVcor}
For an interacting Fock space $\cI$ based on $H$ the following properties are equivalent:
\begin{enumerate}
\item
The operator $\Lambda$ defined by \eqref{****} has an adjoint $\Lambda^*\in\sL(\cI,\ol{\sF(H)})$.

\item
The corresponding ALV-interacting Fock space is POI.
\end{enumerate}
\ecor

\vspace{.3ex}\noindent
We say, an interacting Fock space based on $H$ is \hl{regular} or \hl{regularly based} on $H$ if $\Lambda$ has a weak adjoint. The corollary says, then, that the POI-interacting Fock spaces obtainable from $\sF(H)$ via positive Fock maps $L$, are precisely the interacting Fock spaces regularly based on $H$.

Based on the following lemma about orthogonal(\bf{!}) dimensions, POI-interacting Fock spaces share an important property.

\blem \label{posdimlem}
Let $H$ be a pre-Hilbert space with a positive operator $L\in\sL(H,\ol{H})$. Define the semiinner product $(x,y):=\AB{x,Ly}$ and put $H_L:=H/\sN_L$. Then there exists an isometry $H_L\rightarrow\ol{H}$. Equivalently: $\dim\ol{H_L}\le\dim\ol{H}$.
\elem

\proof
By Friedrich's theorem, $L$ has a positive extension $\ol{L}\colon\cD_{\ol{L}}\rightarrow\ol{H}$ which is self-adjoint in the usual sense (that is, $\cD_{\ol{L}}=\cD_{\ol{L}^*}$ is the maximal domain in $\ol{H}$ for an adjoint of $\ol{L}$; see Section \ref{regSS}). By spectral calculus, $\ol{L}$ has a unique positive square root $\ol{\lambda}\colon\cD_{\ol{\lambda}}\rightarrow\ol{H}$, where $\cD_{\ol{\lambda}}\supset\cD_{\ol{L}}\supset H$ and $\AB{\ol{\lambda}x,\ol{\lambda}x}=\AB{x,Lx}$ for all $x\in\cD_{\ol{L}}$. By $x+\sN_L\mapsto\ol{\lambda}x$ we define an isometry $H_L\rightarrow\ol{H}$, which extends as an isometry $\ol{H_L}\rightarrow\ol{H}$.\qed

\bcor\label{POIembcor}
If the interacting Fock space $\cI$ based on $H$ is a POI-interacting Fock space, then $\cI$ is \hl{embeddable} in the sense that there exists an isometric Fock map $\xi\colon\cI\rightarrow\ol{\sF(H)}$.
\ecor

\proof
Apply Lemma \ref{posdimlem} to $\cI=\cF(H)/\sN$, component-wise.\qed

\bex\label{nonembedex}
Suppose $H$ is a separable infinite-dimensional Hilbert space and choose a Hamel basis $S$ of $H$. Let $H_1$ be a pre-Hilbert space with orthonormal Hamel basis $\bfam{e_s}_{s\in S}$. Put $H_n:=\zero$ for $n>1$. Then $\cI:=\C\Om\oplus H_1$ with $a^*(s)\Om:=e_s$ is an interacting Fock space based on $H$. But, $\dim\ol{H_1}=2^{\aleph_0}>\dim H=\aleph_0$, so that $H_1$ does not embed into $H$, so $\cI$ is not embeddable.
\eex

\it{A fortiori}, by Corollary \ref{Lam*ALVcor}, this non-embeddable $\cI$ is not regular, too. But while missing embeddability can be repaired (and after ``repairing'' the example is regular; see the discussion following Definition \ref{baselydefi}), there are examples of non-regularity that cannot be repaired. Both are the subject of the next section.

\newpage

\section{(Abstract) interacting Fock spaces} \label{absIFSEC}

The notion of embeddability of an interacting Fock space based on $H$, as defined in Corollary \ref{POIembcor}, has been recognized in Accardi and Skeide \cite{AcSk08} as a property of outstanding importance; we reconfirm this in this paper by, in particular, the results in Section \ref{classSEC}. 

When an interacting Fock space $\cI$ based on $H$ is embeddable, we also will say, $\cI$ is \hl{embeddably based} on $H$. This formulation already suggests what comes next, in that the space $\cI$ may be embeddably based on $H$ or it may not be embeddably based on $H$, depending on \it{how} we base it on $H$. This choice includes both different choices for the creator map $a^*\colon H\rightarrow\cI$ (for fixed $H$) and different choices for $H$ itself. The following new, more flexible definition of (abstract) interacting Fock space makes this precise.

\bdefi \label{IFSdefi}
Let $\bfam{H_n}_{n\in\N_0}$ be a family of pre-Hilbert spaces where $H_0=\Om\C$ for a fixed unit vector $\Om$, the \hl{vacuum}, and put $\cI:=\bigoplus_{n\in\N_0}H_n$. Let $A^*$ be a linear subspace of $\sL(\cI)$ satisfying
\beqn{
\ls A^*H_n
~=~
H_{n+1}
}\eeqn
for all $n\in\N_0$ (a condition that, clearly, replaces \eqref{***} in Definition \ref{ASdefi}). Then $\cI$ is an (\hl{abstract}) \hl{interacting Fock space} (denoted as $\cI=(\bfam{H_n}_{n\in\N_0},A^*)$). Usually, we will omit `abstract', and just say `interacting Fock space'.
\edefi

Clearly, an interacting Fock space based on $H$ (and, therefore, any other interacting Fock space in the preceding sections) is turned into an interacting Fock space by setting $A^*:=a^*(H)$. Conversely, choosing a pre-Hilbert space and a linear surjection $a^*\colon H\rightarrow A^*$, we \hl{base} an interacting Fock space on $H$. Of course, the latter is always possible by choosing an arbitrary inner product on $A^*$, turning it that way into a pre-Hilbert space denoted $H$, and choosing for $a^*$ the identification of $H$ and $A^*$. Note that the resulting interacting Fock space based on $H$ is \hl{injective} in the sense that the creator map $a^*$ is injective. So, every interacting Fock space is trivially not only \hl{baseable} but even \hl{injectively baseable}.

\bdefi \label{baselydefi}
An interacting Fock space $\cI=(\bfam{H_n}_{n\in\N_0},A^*)$ is \hl{embeddable} (\hl{regular}) if we can base it on a pre-Hilbert space $H$ \hl{embeddably} (\hl{regularly}), that is, resulting in an interacting Fock space based on $H$ that is embeddable (regular).
\edefi

We know from Example \ref{nonembedex} that there are interacting Fock spaces based on $H$ that are not embeddable, hence, not regular. In this Section we will show that all interacting Fock spaces can be embeddably based, hence, are embeddable (Theorem \ref{aIFSembthm}), while there exist (in abundance) interacting Fock spaces that cannot be regularly based, hence, are not regular (Theorem \ref{nonregthm}).

For instance, in Example \ref{nonembedex} we just have chosen a bad \hl{basing} $a^*\colon H\rightarrow A^*:=a^*(H)$. If we replace $H$ with $H_1$ (that is, if we change the inner product on the vector space $H=H_1$) and keep the same $a^*$ now considered as map $H_1\rightarrow A^*$, then $\cI$ is perfectly embeddable. In fact, it sits already as a subspace $\Om\C\oplus H_1$ in $\sF(H_1)$ and for $\xi$ we may choose the canonical embedding. Of course, $\cI$ based in this way on $H_1$ is regular. (Indeed, $L=\id_{H_0}\oplus\id_{H_1}\oplus\bigoplus_{n\ge2}0$.) So, this is not an example of an interacting Fock space that is not regular.

Example \ref{nonembedex} and the proof of Lemma \ref{posdimlem} suggest that missing orthogonal dimension of $\ol{H}$ is an obstacle to embeddability. We show now that this is essentially the only obstacle. The following lemma, making a somehow quite obvious statement with a surprisingly difficult proof, is key.

\blem\label{totdimlem}
Let $S$ be a total subset of a Hilbert space $H$. Then $\dim H\le\#S$.
\elem

\proof
Choose a well-order $\le$ on $S$. For each $s\in S$ define
\beqn{
H_s
~:=~
\bCB{s'\colon s'<s}^{\perp\perp}
}\eeqn
if $s$ is non-minimal, and define $H_s:=\zero$ if $s$ is minimal. ($H_s$ is the Hilbert subspace of $H$ generated by all $s'$ with $s'<s$. Note that $s$ may be an element of $H_s$ or not.) Define the function $f\colon S\rightarrow H$ by setting
\beqn{
f(s)
~:=~
(\id_H-p_s)s,
}\eeqn
where $p_s$ is the projection onto $H_s$. Note that $\C s+H_s=\C f(s)+H_s$, but $f(s)$ is perpendicular to $H_s$, while $s$ need not be. In particular, both spaces are closed since $\C f(s)+H_s$ is closed. Note, too, that $H_s$ can also be written as
\beqn{
H_s
~=~
\ol{\bigcup_{t<s}\CB{s'\colon s'\le t}^{\perp\perp}}.
}\eeqn

We claim $\bCB{s'\colon s'\le s}^{\perp\perp}=\bCB{f(s')\colon s'\le s}^{\perp\perp}$ for all $s\in S$. Indeed, denote by $\Sigma$ the set of all $s\in S$ for which the statement is true. For some $s\in S$ suppose that $t\in\Sigma$ for all $t<s$. This means in particular that
\beqn{
H_s
~=~
\ol{\bigcup_{t<s}\CB{s'\colon s'\le t}^{\perp\perp}}
~=~
\ol{\bigcup_{t<s}\CB{f(s')\colon s'\le t}^{\perp\perp}}.
}\eeqn
Then
\bmun{%\tag{$*$}
\CB{s'\colon s'\le s}^{\perp\perp}
~=~
\BCB{\CB{s}\cup\bigcup_{t<s}\CB{s'\colon s'\le t}}^{\perp\perp}
~=~
\ol{\C s+\ol{\bigcup_{t<s}\CB{s'\colon s'\le t}^{\perp\perp}}}
~=~
\C s+H_s
~=~
\C f(s)+H_s
\\
~=~
\ol{\C f(s)+\ol{\bigcup_{t<s}\CB{f(s')\colon s'\le t}^{\perp\perp}}}
~=~
\BCB{\CB{f(s)}\cup\bigcup_{t<s}\CB{f(s')\colon s'\le t}}^{\perp\perp}
~=~
\CB{f(s')\colon s'\le s}^{\perp\perp},
}\emun
so that also $s\in\Sigma$. By transfinite induction, $\Sigma=S$.

Define $S_0:=\CB{s\in S\colon f(s)\ne0}$. For each $s\in S_0$ put $e_s:=\frac{f(s)}{\norm{f(s)}}$ so that all $e_s$ ($s\in S_0$) are unit vectors. Put $E:=\bfam{e_s}_{s\in S_0}$. Since $\cls\CB{e_s\colon s\in S_0}=\cls\CB{f(s)\colon s\in S}$, $E$ is total. Since $e_s\perp H_s$ for all $s\in S$, the set $E$ is orthonormal. So $E$ is an ONB. Therefore, $\dim H=\#S_0\le\#S$.\qed

\bthm\label{aIFSembthm}
Every (abstract) interacting Fock space is embeddable.
\ethm

\proof
Let $\cI=(\bfam{H_n}_{n\in\N_0},\Om,A^*)$ be an interacting Fock space. Choose a vector space basis $S$ of $A^*$. Equip $A^*$ with the inner product that makes $S$ orthonormal and denote the arising pre-Hilbert space by $H$. Let $a^*$ denote the canonical identification. Then $\cI$ is an interacting Fock space based on $H$.

For each $n\in\N$, the set $\Lambda_n(S\otimes\ldots\otimes S)$ spans $H_n$. In particular, it is total for $\ol{H_n}$. By Lemma \ref{totdimlem}, $\dim\ol{H_n}\le\#S^n=\dim\ol{H^{\otimes n}}$. By Lemma \ref{posdimlem}, there exists an isometry $\xi_n\colon H_n\rightarrow\ol{H^{\otimes n}}$. Then the Fock map $\xi$ with components $\xi_n$ is the desired isometry.\qed

\lf
In a sense, this shows that it is good to look at interacting Fock spaces as abstract ones. If they come along with a basing $a^*\colon H\rightarrow A^*$ and turn out to be embeddable, this is fine. However, if an interacting Fock space based on $H$ turns out to be not embeddable, then it is better to change the basing. The results that follow from embeddability (see Section \ref{classSEC}) are too important to allow their loss by insisting in an unfortunate choice of a basing.

We thank Roland Speicher who asked, when the second author was on sabbatical leave in Kingston, if the condition of embeddability might not be automatic. The answer, it is automatic provided we choose a reasonable basing, confirms his suspect \it{cum grano salis}. The proof turned out to be much more subtle than expected. It would not have been possible without the crucial Lemma \ref{totdimlem}. Despite making a sufficiently natural and intuitive statement, we felt that its proof was particularly difficult to find.

\lf
After having shown that every interacting Fock space is embeddable, of course, we wish to know if the same is true for regularity: Is every interacting Fock space regular, that is, does every interacting Fock space arise, by basing it appropriately on a suitable pre-Hilbert space, as a POI-interacting Fock space? The following theorem answers this question in the negative sense.

Clearly, for every interacting Fock space $\cI=(\bfam{H_n}_{n\in\N_0},A^*)$, necessarily $A^*\subset\sL_{(1)}(\cI)$. On the other hand, for every sequence $\bfam{H_n}_{n\in\N_0}$ of pre-Hilbert spaces $H_n$ with $H_0=\Om\C$ for some unit vector $\Om$, the pair $\cI=(\bfam{H_n}_{n\in\N_0},\sL_{(1)}(\cI))$ is an interacting Fock space, provided that $H_n=\zero$ implies $H_{n+1}=\zero$. (This is the only, necessary and sufficient, condition that  assures that we get all of $H_{n+1}$ by applying degree one maps to elements of $H_n$.)  We call such $\cI$ the \hl{full interacting Fock space} over $\bfam{H_n}_{n\in\N_0}$. We say, an interacting Fock space is \hl{non-nilpotent} if $H_n\ne\zero$ for all $n$.

\bthm \label{nonregthm}
Every non-nilpotent full interacting Fock space is non-regular.
\ethm

\proof
Let $\cI=(\bfam{H_n}_{n\in\N_0},\sL_{(1)}(\cI))$ be a non-nilpotent full interacting Fock space. Choose a (sufficiently big) pre-Hilbert space $H$ and a surjective linear map $a^*\colon H\rightarrow\sL_{(1)}(\cI)$.

Since $H_n\ne\zero$ for all $n$, we may fix a sequence of unit vectors $\Om_n\in H_n$ (for simplicity, with $\Om_0=\Om$). For every sequence $\bfam{c_n}_{n\in\N}$ of complex numbers define the operator $c:=\sum_{n\in\N}\Om_nc_n\Om_{n-1}^*$ in $\sL_{(1)}(\cI)$. Since $a^*$ is surjective, there exists $x_c\in H$ such that $a^*(x_c)=c$. By definition,
\beqn{
\Lambda(x_c^{\otimes n})
~=~
c^n\Om
~=~
\Om_nc_n\ldots c_1,
}\eeqn
so, $\AB{\Om_n,\Lambda(x_c^{\otimes n})}=c_n\ldots c_1$.

Now, if $\Lambda$ had an adjoint $\Lambda^*\in\sL(\cI,\ol{\sF(H)})$, then
\beqn{
\abs{c_n\ldots c_1}
~=~
\abs{\AB{\Lambda^*\Om_n,x_c^{\otimes n}}}
~\le~
\norm{\Lambda^*\Om_n}\,\snorm{x_c^{\otimes n}}
~=~
\norm{\Lambda^*\Om_n}\,\norm{x_c}^n
}\eeqn
Since $\Lambda^*\Om_n\ne0$ (for instance, because $\Lambda$ is surjective, or by inserting the special choice $c_k=1$ for all $k$), we would get
\beqn{
\norm{x_c}
~\ge~
\sqrt[n]{\frac{\,\abs{c_n\ldots c_1}\,}{\norm{\Lambda^*\Om_n}}}
}\eeqn
for all $c$ and $n$. Choosing $c_n>0$ recursively by setting $c_1:=\norm{\Lambda^*\Om_1}$ and
\beqn{
c_{n+1}
~:=~
(n+1)^{n+1}\frac{\norm{\Lambda^*\Om_{n+1}}}{c_n\ldots c_1},
}\eeqn
we would get for this particular choice of $c$ that
\beqn{
\norm{x_c}
~\ge~
\sqrt[n+1]{c_{n+1}\frac{\,c_n\ldots c_1\,}{\norm{\Lambda^*\Om_{n+1}}}}
~=~
n+1
}\eeqn
for all $n$. As this is not possible, $\Lambda$ cannot have an adjoint.\qed

\lf
Note that, in particular, the interacting Fock space $\cI=\oplus_{n\in\N_0}\Om_n\C$ with $A^*=\sL_{(1)}(\cI)$ is not regular. Of course this changes entirely if we take $\cI=\sF(\C)$ with the usual creators $\ell^*(\C)=\ell^*(1)\C$ which are only a quite small subset of $\sL_{(1)}(\sF(\C))$. More generally, also the one-mode interacting Fock spaces (Example \ref{1-mex}) are regular independently of the number of direct summands. This shows how very much the structure of an interacting Fock space depends on how many creators we allow on the pre-Hilbert space $\cI$.

In Theorem \ref{aIFSembthm}, we completely settled the question of embeddability; we will not be able to do the same for regularity in this paper.
Some (non-)possibilities open up several directions for future work, and will be hinted at in Sections \ref{classSS} and \ref{2FockSS}. Examples \ref{kappaunbex} and \ref{Phinobex} present other non-regular interacting Fock space.

Full interacting Fock spaces with their operator (\nbd{*})algebras generated by $A^*=\sL_{(1)}(\cI)$ are not ``bad guys''. In fact, we will see in Section \ref{boundSEC} that these (possibly unbounded) operator (\nbd{*})algebras are analogues of tensor algebras \cite{MuSo98} (Pimsner-Toeplitz algebras \cite{Pim97}). Theorem \ref{nonregthm} just tells we might be better off, not looking at them as operator algebras of an interacting Fock space.

\newpage

\section{Squeezings: Embedded interacting Fock spaces} \label{classSEC}

In this section, we examine the consequences of having a Fock embedding $\xi\colon\cI\rightarrow\ol{\sF(H)}$ of an interacting Fock space based on $H$ into $\ol{\sF(H)}$. The formula in \eqref{**}, which expresses the images $\xi a^*(x)\xi^*$ of the creators when acting on the subspace $\xi\cI$ of $\ol{\sF(H)}$ in terms of the usual Fock creators $\ell^*(x)$ \it{squeezed} by an operator $\vk$ as $\vk\ell^*(x)$, has been observed already in \cite{AcSk08}. But in this section we go far beyond \cite[Theorem 5.5 and Corollary 5.7]{AcSk08}, and obtain a classification of interacting Fock spaces (based on $H$ or not) in terms of such \it{squeezings} $\vk$.

This is the moment to specify better when we consider two interacting Fock spaces to be ``the same''.  Recall that we have the two fundamentally different notions of interacting Fock space and interacting Fock space based on $H$, the latter being ``the same'' as ALV-interacting Fock space, while POI-interacting Fock spaces are a subspecies of ALV-interacting Fock spaces corresponding to interacting Fock spaces that are based regularly.

\bdefi \label{isodefi}
\begin{enumerate}
\item
The interacting Fock spaces $\cI=(\bfam{H_n}_{n\in\N_0},A^*)$ and $\cI'=(\bfam{H'_n}_{n\in\N_0},{A^*}')$ are \hl{isomorphic} if there exists a \hl{Fock unitary} $u=\bigoplus_{n\in\N_0}u_n$ (that is, the $u_n$ are unitaries $H_n\rightarrow H'_n$ and $u_0=\id_{\Om\C}$, where $H_0=\Om\C=H'_0$) such that
\beqn{
uA^*u^*
~=~
{A^*}'.
}\eeqn

\item
The interacting Fock spaces $\cI=(\bfam{H_n}_{n\in\N_0},a^*)$ and $\cI'=(\bfam{H'_n}_{n\in\N_0},{a^*}')$ based on (the same) $H$ are \hl{isomorphic} if
\beqn{
\Lambda X
~\longmapsto~
\Lambda' X
}\eeqn
($X\in\sF(H)$) defines a unitary $u\colon\cI\rightarrow\cI'$.
\end{enumerate}
\edefi

\noindent
There are other reasonable notions of isomorphism (see also Section \ref{autoSS}), which we postpone to future work. We collect some more or less obvious properties.

\bob \label{isoob}
Recall that for an interacting Fock space based on $H$ not only $\Lambda$ is defined in terms of $a^*$ by \eqref{****}, but that also $\Lambda$ determines $a^*$ via $a^*(x)\Lambda X=\Lambda(x\otimes X)$.

\begin{enumerate}
\item
Clearly, the unitary $u$ for isomorphic interacting Fock spaces based on $H$, is a Fock unitary. Moreover, by the preceding reminder, $ua^*(x)u^*={a^*}'(x)$ for all $x\in H$. Therefore, isomorphic interacting Fock spaces based on $H$ are also isomorphic as interacting Fock spaces.

Conversely, suppose we have two interacting Fock spaces that are isomorphic via $u$. If we base the first one on $H$ via $a^*\colon H\rightarrow A^*$, then by ${a'}^*\colon x\mapsto ua^*(x)u^*$ we turn the second one into an isomorphic interacting Fock space based on $H$. Moreover, by $\Lambda X\mapsto\Lambda'X$ we recover the $u$ we started with.

\item
By the discussion following Definition \ref{ASdefi}, we know:  Every ALV-interacting Fock space $\cI=(H,\bfam{(\bullet, \bullet)_n}_{n\in\N_0})$ \bf{is} (understood as) an interacting Fock space based on $H$ with the canonical basing $a^*\colon x\mapsto a^*(x)$ (where $a^*(x)$ are the creators with which an ALV-interacting Fock space comes along); every interacting Fock space $\cI=(\bfam{H_n}_{n\in\N_0},a^*)$ based on $H$ is (canonically) \bf{isomorphic} to the ALV-interacting Fock space coming from the semiinner product $(\bullet,\bullet):=\AB{\Lambda\bullet,\Lambda\bullet}$ on $\sF(H)$. Moreover, if $\cI'=(H,\bfam{(\bullet, \bullet)'_n}_{n\in\N_0})$ is another ALV-interacting Fock space isomorphic to the interacting Fock space $\cI$ based on $H$, then $(\bullet, \bullet)'=\AB{\Lambda'\bullet,\Lambda'\bullet}=\AB{\Lambda\bullet,\Lambda\bullet}=(\bullet, \bullet)$, that is, as ALV-Fock interacting space it is identical to the ALV-interacting Fock space arising from $\cI$. \it{A fortiori} there is one and only one POI-interacting Fock space isomorphic to a given interacting Fock space regularly based on $H$. 
\end{enumerate}
\eob

\noindent
We now fix an interacting Fock space $\cI$ based on $H$ and assume it is embedded via a fixed $\xi\colon\cI\rightarrow\ol{\sF(H)}$. In this situation, we say $\cI$ is an \hl{embedded} interacting Fock space, and it is understood that an interacting Fock space to be embedded has to be based.

By assuming that $\cI$ is embeddably based (always possible by Theorem \ref{aIFSembthm}), we do not lose any interacting Fock space. (After all, choosing a basing does not change the interacting Fock space.) Clearly, the Fock isometry $\xi$ may be viewed as a Fock unitary $u_\xi$ onto $\xi\cI\subset\ol{\sF(H)}$. Clearly, defining $a^*_\xi\colon x\mapsto u_\xi a^*(x)u_\xi^*$ turns $\xi\cI$ into an interacting Fock space based on $H$ isomorphic to $\cI$. (And if we started with another embedding $\xi'$, then the interacting Fock spaces $\xi\cI$ and $\xi'\cI$ based on $H$ are isomorphic via $u_{\xi'}u_\xi^*$.) So, starting with an interacting Fock space embeddably based on $H$, actually embedding it, we stay in the same isomorphism class of interacting Fock spaces based on $H$.

We distinguished carefully between the unitary $u_\xi$ onto $\xi\cI$ and the isometry $\xi$ into $\ol{\sF(H)}$. We, tacitly, used already that a unitary $u$ between pre-Hilbert spaces always has an adjoint, namely, $u^*=u^{-1}$. This need not be so, for an isometry. (One may show that an isometry has an adjoint if and only if its range is complemented in its codomain; see, for instance, Skeide \cite[Proposition 1.5.13]{Ske01}.) Fortunately, our isometry $\xi$ goes into a Hilbert(!) space and, like every isometry from a pre-Hilbert space into a Hilbert space, it has a densely defined, surjective adjoint $\xi^*\colon\cD_{\xi^*}:=\xi\cI\oplus(\xi\cI)^\perp\rightarrow\cI$, determined by $\xi^*(\xi x)=x$ and $\xi^*y=0$ for $y\in(\xi\cI)^\perp$. (The complement $(\xi\cI)^\perp$ is taken in the Hilbert space $\ol{\sF(H)}$, and since $\cD_{\xi^*}$ has zero-complement in this Hilbert space, it is dense; the last conclusion may fail for subspaces of pre-Hilbert spaces.) It follows that
 \beqn{
 a
 ~\longmapsto~
 \xi a\xi^*
 }\eeqn
defines an algebra monomorphism from the algebra $\sL(\cI)$ onto the corner $\sL(\xi\cI)\subset \sL(\cD_{\xi^*})=\sL\rtMatrix{\xi\cI\\(\xi\cI)^\perp}=\rtMatrix{\sL(\xi\cI)&\sL((\xi\cI)^\perp,\xi\cI)\\\sL(\xi\cI,(\xi\cI)^\perp)&\sL((\xi\cI)^\perp)}$. If $a$ has an adjoint, $a^*$, then $\xi a^*\xi^*$ is, clearly, an adjoint of $\xi a\xi^*$. So, $\xi\bullet\xi^*$, when restricted to $\sL^a(\cI)$ is actually a \nbd{*}monomorphism. Moreover, since $\xi$ respects the vacuum, $\xi\bullet\xi^*$ also respects the \hl{vacuum expectation} $\AB{\Om,\bullet\Om}$. Since $\xi$ is even, also the degree of $a$ is preserved, that is, the monomorphism itself is an even map.

So, via $\xi$, we have identified $\cI$ as a subspace $\xi\cI$ of $\cD_{\xi^*}\subset\ol{\sF(H)}$ and we have identified $A^*$ (and the algebras generated by it) as a subspace $\xi A^*\xi^*$ (and the algebras generated by it) of $\sL(\cD_{\xi^*})$. The map corresponding to $\Lambda$ for this interacting Fock space $\xi\cI$ is $\Lambda_\xi=u_\xi\Lambda$, when considered as map onto $\xi\cI$, as it has to be by definition. However, we prefer to consider it as map $\lambda:=\xi\Lambda\colon\sF(H)\rightarrow\ol{\sF(H)}$, taking also into account that its range is actually $\xi\cI\subset\cD_{\xi^*}\subset\ol{\sF(H)}$. Recall that $\lambda$ depends on $\xi$. But for reasons of readability in formulae with indices, we dispense with the idea, calling it $\lambda_\xi$. (We will rather write $\lambda'$ to indicate when it is originating in another $\xi'$. It is clear that $\lambda'=\xi'\xi^*\lambda$.)

We are now almost ready to formulate and prove \eqref{**} in this general context. The only question that remains to be made precise in order to make sense out of $\vk\ell^*(x)$, is the question of the appropriate domain and codomain of $\vk$. As we wish that $\vk\ell^*(x)=\xi a^*(x)\xi^*$, the codomain should coincide with domain $\cD_{\xi^*}$ of $\xi^*$. The domain should contain what $\ell^*(x)$ generates out of $\cD_{\xi^*}$. We just mention that $\sF(H)=(H\otimes\sF(H))\oplus\Om\C$ in the obvious way; consequently, for every subspace $\cD$ of $\ol{\sF(H)}$ we get the subspace $(H\otimes\cD)\oplus\Om\C$  of $\ol{\sF(H)}$,  and the latter is dense if (and only if) the former is dense.

\bthm \label{**thm}
Let $\cI$ be an embedded (via $\xi$) interacting Fock space (based on $H$). There exists a unique vacuum-preserving map (necessarily also a Fock map) $\vk\in\sL((H\otimes\cD_{\xi^*})~\oplus~\C\Om,\cD_{\xi^*})$ such that
\beqn{
\vk\ell^*(x)
~=~
\xi a^*(x)\xi^*
}\eeqn
(that is, Equation \eqref{**}).Therefore, the algebra monomorphism $a\mapsto\xi a\xi^*$ sends $a^*(x)$ to $\vk\ell^*(x)$.

If $a^*(x)$ has an adjoint $a(x)\in\sL(\cI)$, then $\xi a(x)\xi^*=(\vk\ell^*(x))^*$ (though, $\vk$ need not be adjointable).

In either case, the (\nbd{*})monomorphism respects the vacuum state $\AB{\Om,\bullet\Om}$.

Moreover, $\lambda$ can be recovered as the unique Fock map satisfying the equation
\beq{ \label{lkeq}
\lambda
~=~
\vk((\id_H\otimes\lambda)+\id_{\Om\C}),
}\eeq
that is, as the unique $\lambda$ whose components satisfy the recursion
\beqn{
\lambda_{n+1}
~=~
\vk_{n+1}(\id_H\otimes\lambda_n)
\text{~~~and~~~}
\lambda_0
~=~
\id_{\C\Om},
}\eeqn
that is,
\beqn{
\lambda_n
~=~
\vk_n(\id_H\otimes\vk_{n-1})\ldots(\id_H^{\otimes(n-1)}\otimes\vk_1)~~~(n\geq 1).
}\eeqn
\ethm

\proof[Notes on the proof.]
Why `notes on the proof'? Well, for adjointable interacting Fock spaces and without the uniqueness statement, this theorem is \cite[Theorem 5.5 and Corollary 5.7]{AcSk08}. The proof in \cite{AcSk08} does not depend on adjointability, and once we have $\vk$ satisfying \eqref{**}, it was just an omission in \cite{AcSk08} not to have noticed uniqueness. However, the proof in \cite{AcSk08} went by first proving (by brute-force linear algebra \cite[Lemma 5.4]{AcSk08}) existence of $\vk$ satisfying the recursion with $\lambda$. And while \eqref{**} fixes $\vk$, the recursion alone does not. (It may be considered a sort of ``lucky punch'' that the freedom in choosing $\vk$ for satisfying the recursion has been used ``wisely'' to also satisfy \eqref{**}.) Starting from \eqref{**} and uniqueness, straightens up and simplifies the proof considerably, so we sketch this briefly.

$\vk$ is determined uniquely by $\vk\ell^*(x)=\xi a^*(x)\xi^*$ on the span of the ranges of all $\ell^*(x)$, that is, on $H\otimes\cD_{\xi^*}$. The remaining uncertainty is taken away by the requirement that $\vk$ is vacuum-preserving.

For existence of $\vk$, we simply put $\vk\Om:=\Om$ and define it on $H\otimes\cD_{\xi^*}=H\otimes(\xi\cI\oplus(\xi\cI)^\perp)$ as \eqref{**} suggests: Necessarily, $\vk(x\otimes Y)=\xi a^*(x)\xi^*Y=0$ for $Y\in(\xi\cI)^\perp$. And for $\xi X\in\xi\cI$ we obtain $\vk(x\otimes\xi X)=\xi a^*(x)\xi^*\xi X=\xi a^*(x)X$. Since $\xi$ is an isometry, this map $\vk$ is well defined.

By definition, this $\vk$ satisfies \eqref{**}. And it is routine (using how $\Lambda$ and $a^*$ determine each other as explained in the beginning of Observation \ref{isoob} and the interplay between $\Lambda$ and $\lambda$ via $\xi$) to verify \eqref{lkeq}.\qed

\lf
Let us sum up again what we achieved. From an interacting Fock space $\cI$ based on $H$ and embedded via $\xi$, we extracted the pre-Fock subspace $\xi\cI$ of $\ol{\sF(H)}$ and the operator $\vk$ from the (dense, pre-Fock) subspace $(H\otimes\cD_{\xi})\oplus\Om\C$ to the (dense, pre-Fock) subspace $\cD_{\xi^*}:=\xi\cI\oplus(\xi\cI)^\perp$. From $\vk$ we reconstruct $\lambda$ via the recursion encoded in \eqref{lkeq}, and from $\lambda$ we reconstruct $\xi a^*(x)\xi^*$. (Or, better, from $\Lambda_\xi$, the surjective corestriction of $\lambda$, we reconstruct $u_\xi a^*(x)u_\xi^*\in\sL(\xi\cI)$ as explained in the beginning of Observation \ref{isoob}, which, when embedded into $\sL(\cD_{\xi^*})$, becomes $\xi a^*(x)\xi^*$.) ~That is, we have encoded the entire information about the embedded interacting Fock space $\cI$, and  up to isomorphism about the interacting Fock space $\cI$ based on $H$, in the operator $\vk$ (including, of course, how its domain and codomain are made up out of $\xi\cI$); and $\vk$, on the other hand, is uniquely determined by $\cI$ and $\xi$, that is, by the embedded interacting Fock space $\cI$. Moreover, if we started from another embedding, $\xi'$, then everything is under control via the partial isometry $\xi'\xi^*$ in the sense that $\vk'=\xi'\xi^*\vk((\id_H\otimes\xi\xi'^*)\oplus\id_{\Om\C})$ and the corresponding $u_{\xi'}u_\xi^*$ is an isomorphism between the interacting Fock spaces $\xi\cI$ and $\xi'\cI$ based on $H$.

Additionally, let us observe that $\vk$ fulfills the following two properties: Firstly, $\vk$ is onto $\xi\cI$ (simply because $\lambda$ is onto $\xi\cI)$). Secondly, $\vk$ is $0$ on the subspace $H\otimes(\xi\cI)^\perp$ (as computed in the proof of Theorem \ref{**thm}).

We now show that these two conditions are the only conditions a Fock-map $\vk$ has to satisfy in order to be the $\vk$ of an embedded interacting Fock space. To that goal, we now free the discussion from the embedding $\xi$.

\bdefi
Let $\cI$ be a pre-Fock subspace of $\ol{\sF(H)}$, and define the dense, pre-Fock subspace $\cD_\cI:=\cI\oplus\cI^\perp$ of $\ol{\sF(H)}$. A Fock map $\vk\colon(H\otimes\cD_\cI)\oplus\Om\C\rightarrow\cD_\cI$ is called a \hl{squeezing} (\hl{relative to $\cI$}) if $\vk$ is onto $\cI$ and vanishes on $H\otimes\cI^\perp$.
\edefi

Observe that the \it{squeezed creators} $\vk\ell^*(x)$ (co)restrict to maps $\cI\rightarrow\cI$, which we denote by $a^*_\vk(x)$. This gives rise to the linear map $a^*_\vk\colon H\rightarrow\sL(\cI)$. Occasionally, we leave out the subscript $\vk$ when there is no danger of confusion.

\bthm \label{kIFSthm}
If $\vk$ is a squeezing relative to $\cI=\bigoplus_{n\in\N_0}H_n\subset\ol{\sF(H)}$, then $\cI_\vk:=(\bfam{H_n}_{n\in\N_0},a^*_\vk)$ is an interacting Fock space based on $H$. Moreover, the (unique) $\vk_{\xi_\vk}$ constructed by Theorem \ref{**thm} from the canonical embedding $\xi_\vk\colon\cI\rightarrow\cI\subset\ol{\sF(H)}$ is $\vk$.
\ethm

\bdefi \label{kdefi}
We call $\cI_\vk$ a \hl{\nbd{\vk}interacting Fock space}, and denote it by $\cI_\vk=(H,\vk)$ (also here leaving occasionally out the subscript).
\edefi

\proof[Proof of Theorem \ref{kIFSthm}.~]
There is not really much to prove. $\vk$ being a squeezing, by surjectivity of $\vk$ it follows that $\vk\ell^*(x)$ maps  $\ol{H^{\otimes n}}$ onto $\cI\cap\ol{H^{\otimes(n+1)}}=H_{n+1}$ and by $\vk$ being $0$ on $H\otimes\cI^\perp$ it follows that to exhaust the range it is sufficient to restrict to $H\otimes(\cI\cap\ol{H^{\otimes n}})=H\otimes H_n$ (and $\Om\C$). Therefore, $\ls a^*_\vk(H)H_n=H_{n+1}$. Clearly, $\vk$ does satisfy \eqref{**} for the canonical embedding $\xi_\vk$, so by the uniqueness statement in Theorem \ref{**thm}, $\vk$ coincides with $\vk_{\xi_\vk}$.\qed

\bcor
We, thus, established a one-to-one correspondence between embedded interacting Fock spaces and squeezings.
\ecor

The following theorem is a mere corollary of Theorems \ref{aIFSembthm} and \ref{**thm}.

\bthm \label{kIFSisothm}
Every interacting Fock space is isomorphic to a \nbd{\vk}interacting Fock space (suitably varying $H$, $\cI\subset\ol{\sF(H)}$, and $\vk$ relative to $\cI$).

Every interacting Fock space based embeddably on $H$ is isomorphic to a \nbd{\vk}interacting Fock space for a squeezing $\vk$ relative to a pre-Fock subspace $\cI$ of $\ol{\sF(H)}$.
\ethm

We have already discussed the influence of different choices $\xi$ how to embed into $\ol{\sF(H)}$ a given interacting Fock space based on $H$. Maybe a bit surprisingly, the answer is the same if we vary also $H$, that is, if we vary also the basing. Without the obvious proof, we state the following:

\bprop \label{kpartprop}
Let $\vk$ and $\vk'$ be squeezings relative to pre-Fock subspaces $\cI\subset\ol{\sF(H)}$ and $\cI'\subset\ol{\sF(H')}$, respectively. Then the interacting Fock spaces $\cI_\vk$ and $\cI'_{\vk'}$ are isomorphic (as interacting Fock spaces) if and only if there is a partial Fock isometry $v\in\sB(\cD_\cI,\cD_{\cI'})$ with $v^*v=p_\cI$ and $vv^*=p_{\cI'}$ ($p_{\cI^{(')}}$ the projection in $\sB(\cD_{\cI^{(')}})$ onto $\cI^{(')}$) such that
\beqn{
\vk'
~=~
v\vk((\id_H\otimes v^*)\oplus\id_{\Om\C}).
}\eeqn
\eprop

We see for getting an interacting Fock space as a \nbd{\vk}interacting Fock space, it does not only not matter (via an extremely obvious relation among different $\vk$) how we embed it, but it does not even depend (up to the same obvious relation) on how we based it, as long as we based it embeddably. Also: $\vk$ is bounded if and only $\vk'$ is bounded; that is, no change of basing makes an interacting Fock space that has an unbounded $\vk$ into one that has a bounded $\vk$. Since we have examples of bounded interacting Fock spaces with unbounded $\vk$, restricting to bounded $\vk$ will lose us examples of bounded interacting Fock spaces.

Recall that \nbd{\vk}interacting Fock spaces \bf{are} embedded interacting Fock spaces and, therefore, based. Some properties of an interacting Fock space (for instance, boundedness of the set $A^*$) are intrinsic; other properties (for instance, regularity of a basing) depend on the basing. This raises several question how these properties can be seen by looking only at $\vk$, or by guaranteeing existence of certain good choices for $\vk$. Regarding regularity -- a property with reference to a given basing --, we close this section by stating the quite obvious result that regularity does not depend on the representative within the same isomorphism class of interacting Fock spaces based on the same pre-Hilbert space $H$.

\bprop \label{lLadprop}
If $\cI$ and $\cI'$ are isomorphic interacting Fock space based on $H$, then $\cI$ is regular if and only if $\cI'$ is regular.
\eprop

\proof
Let $u$ be the isomorphism. Then if $\Lambda^*$ exists, $\Lambda^*u^*$ is an adjoint of $\Lambda'$, and \it{vice versa}.\qed

\bcor \label{lLcor}
Suppose $\xi$ is a Fock embedding into $\ol{\sF(H)}$. Then $\cI$ is regular if and only if $\xi\cI$ is regular, that is, if $\lambda$ has an adjoint.
\ecor

In the following section we address questions of boundedness. More general questions require more refined notions of isomorphism and more reasonable choices for our basings. As with this we run into problems that do not allow for a single solution but split into subclasses, we postpone the discussion, indicating some future work in Section \ref{olSEC}.

\newpage

\section{Boundedness: Cuntz-Pimsner-Toeplitz algebras} \label{boundSEC}

As already noticed in Accardi and Skeide \cite[Section 4]{AcSk08}, if $\cI$ is an adjointable(\bf{!}) interacting Fock space (in \cite{AcSk08} based on $H$, but that is irrelevant), then we may define the \hl{full Fock module}
\beqn{
\sF(\sL_{(1)}^a(\cI))
~:=~
\sL^a_{(0)}(\cI)
\,\oplus\,
\bigoplus_{n\in\N}\ls\bfam{\,\underbrace{\!\sL_{(1)}^a(\cI)\ldots\sL_{(1)}^a(\cI)\!}_{n\text{ times}}\,}
}\eeqn
($\sL_{(n)}^a$ denoting the adjointable part of $\sL_{(n)}$) on which the elements of $A^*$ act by operator multiplication. How is this a Fock module? Well, $\sL^a_{(0)}(\cI)$ is a \nbd{*}algebra of operators in $\sL^a(\cI)$ and for each $n$ ($n=0$ included), $\sL^a_{(n)}(\cI)$ is a bimodule over $\sL^a_{(0)}(\cI)$ with an \it{inner product} $\AB{X_n,Y_n}:=X_n^*Y_n$ in $\sL^a_{(0)}(\cI)$. Moreover, the tensor product $\sL^a_{(n)}(\cI)\odot\sL^a_{(m)}(\cI)$ over $\sL^a_{(0)}(\cI)$ sits naturally as ~$\ls\sL^a_{(n)}(\cI)\sL^a_{(m)}(\cI)$~ in $\sL^a_{(n+m)}(\cI)$. We do not explain in detail how to make this more precise.%
\footnote{
It occupies the whole lengthy \cite[Section 3]{AcSk08} (see also Skeide \cite[Appendix C]{Ske01}) to develop a notion of positivity for general \nbd{*}algebras that is sufficiently general for applications (for instance, the \it{square of white noise} Fock module in \cite{AcSk00a}) and still allows to control positivity in tensor products, before the Fock module of an interacting Fock space $\cI$ can be defined in \cite[Section 4]{AcSk08}.
}
Here, we are interested in the case when $A^*$ consists of \bf{bounded} operators. We shall say, $\cI$ is a \hl{bounded} interacting Fock space. In this case, we really get a (completed) full Fock module and embed the operators and algebras into Cuntz-Pimsner-Toeplitz algebras. In the end, we free this from the unnecessary hypothesis that the elements of $A^*$ are adjointable. Criteria that show how boundedness of $A^*$ is reflected by other ways to describe interacting Fock spaces ($\vk$, $\lambda$, $L$, ...), are postponed to Section \ref{bcritSEC}.

Since in this section we put emphasis on $A^*$ and do not  consider $\cI$ to be based (Example \ref{CPTex} being the only exception), $a^*$ stands for a typical element of $A^*$, and not for a basing.

Clearly, if $A^*\subset\sB^a(\cI)$, then we restrict everything to the bounded portions, and define
\beqn{
\sF(\sB_{(1)}^a(\cI))
~:=~
\sB^a_{(0)}(\cI)
\,\oplus\,
\bigoplus_{n\in\N}\ls\bfam{\,\underbrace{\!\sB_{(1)}^a(\cI)\ldots\sB_{(1)}^a(\cI)\!}_{n\text{ times}}\,},
}\eeqn
on which, again, the elements of $A^*$ act by operator multiplication. Now, $\sB^a_{(0)}(\cI)$ is a pre-\nbd{C^*}algebra of operators in $\sB^a(\cI)$ and $\sB^a_{(n)}(\cI)$ is a pre-correspondence (that is like a correspondence but not necessarily complete and possibly over a pre--\nbd{C^*}algebra with \hl{contractive} left action) over $\sB^a_{(0)}(\cI)$. (Even if all $H_n$ are Hilbert spaces, $\cI$, and with $\cI$ also $\sB^a_{(0)}(\cI)$ and $\sB^a_{(n)}(\cI)$, will not be complete, unless $\cI$ is nilpotent.) We may complete, and obtain
\beq{ \label{Ba1Fock}
\sF(\ol{\sB_{(1)}^a(\cI)})
~:=~
\ol{\sB^a_{(0)}(\cI)}
\,\oplus\,
\coplus_{n\in\N}\ls\bfam{\,\underbrace{\!\sB_{(1)}^a(\cI)\ldots\sB_{(1)}^a(\cI)\!}_{n\text{ times}}\,}
~=~
\coplus_{n\in\N_0}\sB_{(1)}^a(\cI)^{\odot n}.
}\eeq

\brem
Still, while $\ol{\sB_{(n)}^a(\cI)\odot\sB_{(m)}^a(\cI)}$ is contained in $\ol{\sB_{(n+m)}^a(\cI)}$, it is usually only a proper subset. If we insist in equality, we have to pass to the von Neumann objects $\sB_{(0)}\bfam{\,\ol{\cI}\,}=\ol{\sB_{(0)}^a(\cI)}^s$ and $\sB_{(1)}\bfam{\,\ol{\cI}\,}=\ol{\sB_{(1)}^a(\cI)}^s$ (strong closure in $\sB(\ol{\cI})$). We ignore this ramification in this paper.
\erem

Let us briefly recall a couple of general facts about full Fock modules and inducing representations. (See below and see Footnote \ref{corrFN} for \it{correspondence}.)

Firstly, if $E$ is a correspondence over a \nbd{C^*}algebra $\cB$, then the \hl{full Fock module} over $E$ is the correspondence $\sF(E):=\coplus_{n\in\N_0}E^{\odot n}$. (Here $E^{\odot 0}:=\cB$. But if $\cB$ is unital, then we will write it as $E^{\odot 0}:=\om\cB$, with the central unit vector $\om:=\U\in\cB$.) For each $x\in E$, the \hl{creator} $\ell^*(x)\colon X\mapsto x\odot X$ is an adjointable operator on $\sF(E)$, denoted $\ell^*(x)\in\sB^a(\sF(E))$. Since $\sF(E)$ is a correspondence and since $\cB$ acts faithfully from the left on the direct summand $E^{\odot 0}=\cB$, also $\cB$ sits as a \nbd{C^*}subalgebra in $\sB^a(\sF(E))$. The \hl{tensor algebra} over $E$ is the Banach subalgebra of  $\sB^a(\sF(E))$ generated by $\ell^*(E)$ and $\cB$ (Muhly and Solel \cite{MuSo98}). The \hl{Cuntz-Pimsner-Toeplitz algebra} over $E$ is the \nbd{C^*}subalgebra of  $\sB^a(\sF(E))$ generated by $\ell^*(E)$ and $\cB$ (Pimsner \cite{Pim97}).

Secondly, if $E$ is a Hilbert \nbd{\cB}module and if $G$ is a correspondence from $\cB$ to $\C$ (that is, $G$ is a Hilbert space with a nondegenerate representation of $\cB$), then $\sB^a(E)$ acts (nondegenerately) on the Hilbert space $E\odot G$ via $\sB^a(E)\ni a\mapsto a\odot\id_G\in\sB(E\odot G)$. If the correspondence $G$ is \hl{faithful} (that is, if the left action defines a faithful representation of $\cB$), then also the action of $\sB^a(E)$ on $\sB(E\odot G)$ is faithful. (In our applications to the Fock module $\sF(E)$, $G$ will be ``very non-faithful'' and we have to work to show by hand that the action of $\sB^a(E)$ for that $G$ is, nevertheless, faithful.) If $E$ is a \hl{correspondence} from $\cA$ to $\cB$ (that is, the left action of $\cA$ on the Hilbert \nbd{\cB}module $E$ defines a nondegenerate homomorphism), then the canonical homomorphism $\cA\rightarrow\sB^a(E)\rightarrow\sB(E\odot G)$ defines a nondegenerate representation of $\cA$ on $E\odot G$ (turning $E\odot G$ into a correspondence from $\cA$ to $\C$), the representation \hl{induced} from (the representation on) $G$ by $E$.%
\footnote{ \label{corrFN}
There are several definitions of \it{\nbd{C^*}correspondence} around. Despite the possibility to construct (tensor products and) the full Fock module also over Hilbert \nbd{\cB}modules with a degenerate left action by $\cB$, in several places in the theory, to our taste, degeneracy of the left action is not acceptable. (Just one instance: The algebra should act as ``identity correspondence'' under tensor product.) So, we insist that a correspondence, to merit the name, has nondegenerate left action, by definition. On the other hand, while many authors allow for degenerate left action, in the construction of the full Fock module they insist in that the correspondence should be full, which we do not.
}

\lf
After these reminders, we return to the beginning. The $\sF(\ol{\sB_{(1)}^a(\cI)})$ defined in \eqref{Ba1Fock} is, indeed, the full Fock module $\sF(E)$ for the correspondence $E:=\ol{\sB_{(1)}^a(\cI)}$ over the (unital!) \nbd{C^*}algebra $\cB:=\ol{\sB_{(0)}^a(\cI)}$. We wish to identify the \nbd{C^*}algebra $\sB^a(\sF(E))$ as a subalgebra of $\sB\bfam{\,\ol{\!\cI}\hspace{.1ex}}$; and we wish to do it in such a way that the creators $\ell^*(a_1)\in\sB^a(\sF(E))\subset\sB\bfam{\,\ol{\!\cI}\hspace{.1ex}}$ act like the operators $a_1\in E\subset\sB^a(\cI)\subset\sB\bfam{\,\ol{\!\cI}\hspace{.1ex}}$ act on $\cI\subset\,\ol{\!\cI}$. For that goal, we tensor $\sF(E)$ with the representation space $G:=H_0=\Om\C\subset\cI$ of $\cB$, which is left invariant by $\cB$ because all elements of $\cB$ are even. (Tensoring with $\cI$ would, yes, guarantee faithfulness of the representation on $\sF(E)\sodot\cI\cong\coplus_{n\in\N_0}H_n^{n+1}$, but this space would be much too big, and it also would be quite tedious to invent a good notation for how $\ell^*(a_1)\odot\id_\cI$ acts between the several direct summands.) The following proposition is obvious. (It also triggers a gapless proof for \cite[Theorem 4.1]{AcSk08}.)

\lf
\bprop \label{TCPTprop}
The map
\beqn{
\bfam{a^*_n\odot\ldots\odot a^*_1}\odot\Om
~\longmapsto~
a^*_n\ldots a^*_1\Om
~~~~~~
(a^*_i\in A^*\subset E)
}\eeqn
defines a unitary $\sF(E)\sodot H_0\rightarrow\,\ol{\!\cI}$ and, under this isomorphism, $\ell^*(a_1)\odot\id_{H_0}=a_1$ for all $a_1\in E\subset\sB\bfam{\,\ol{\!\cI}\hspace{.1ex}}$. Therefore, the map $a^*\mapsto\ell^*(a^*)$ $(a^*\in A^*)$ extends to a completely isometric isomorphism from the (\nbd{*})algebra $\cA^{(*)}$ generated by $A^*$ onto the (\nbd{*})subalgebra of the tensor algebra (the Cuntz-Pimsner-Toeplitz algebra) of $E$ generated by $\ell^*(A^*)$.
\eprop

\brem \label{BFErem}
Note that also the representation of $\cB\subset\sB^a(\sF(E))$ on $\sF(E)\odot  H_0$, under the isomorphism with $\,\ol{\!\cI}$, is just the identity representation. This is enough to show that the representation $\sB^a(\sF(E))\rightarrow\sB^a(\sF(E))\odot\id_{H_0}\subset\sB\bfam{\,\ol{\!\cI}\hspace{.1ex}}$ of $\sB^a(\sF(E))$ (containing the Cuntz-Pimsner-Toeplitz algebra of $E$, containing the tensor algebra of $E$) on $\,\ol{\!\cI}$ is faithful. (Indeed, first of all  for $0\ne a\in\sB^a(\sF(E))$ there exist $k,m,n$ and $X_n\in E^{\sodot n}$, $Y_m\in E^{\sodot m}$, $Z_k,Z'_k\in E^{\sodot k}$ such that $\AB{(Z_k\odot\Om),(\AB{X_n,aY_m}\odot\id_{H_0})(Z'_k\odot\Om)}\ne0$. (Recall that $\AB{X_n,aY_m}\in\cB$ is even, and that if $\AB{X_n,aY_m}\ne0$, then also $\AB{X_n,aY_m}\odot\id_{H_0}\ne0$.) By
\beqn{
0
~\ne~
\AB{(Z_k\odot\Om),(\AB{X_n,aY_m}\odot\id_{H_0})(Z'_k\odot\Om)}
~=~
\BAB{((X_n\odot Z_k)\odot\Om)\,,\,(a\odot\id_{H_0})\,((Y_n\odot Z'_k)\odot\Om)},
}\eeqn
we see $a\odot\id_{H_0}\ne0$.) We do not really need that result. Nevertheless, it is surely worthwhile mentioning it.
\erem

$\cB$ and $E$, as defined above, are rather big. (If we passed to the von Neumann case, that is, taking strong closures everywhere, we would end up with the type I von Neumann algebras $\ol{\cB}^s$ and $\ol{\sB^a(\sF(E))}^s$ which have isomorphic atomic centers $\ell^\infty$.) In view of our interest in the Banach (\nbd{*})algebra generated by $A^*$, we had better try and keep the tensor algebra (the Cuntz-Pimsner-Toeplitz algebra) into which we embed as small as reasonably possible. More precisely, instead of $E$ and $\cB$ we had better pass to a subspace $F\subset E$ and to a \nbd{C^*}subalgebra $\cC\subset\cB$ such that $F$ still contains $A^*$ and such that $F$ is a correspondence over $\cC$ with respect to the inner product and bimodule operations inherited from $\sB\bfam{\,\ol{\!\cI}\hspace{.1ex}}\supset F,\cC$.

\bcor \label{smallcor}
Under these conditions, Proposition \ref{TCPTprop} remains true. That is, $\sF(F)\sodot H_0\cong\,\ol{\!\cI}$ via the same isomorphism, and $a^*\mapsto\ell^*(a^*)$ $(a^*\in A^*)$ extends to (completely isometric) embeddings of the tensor algebra and the Cuntz-Pimsner-Toeplitz algebra of $F$ into $\sB^a(\sF(F))$. (Also Remark \ref{BFErem} remains true.)
\ecor

Note that even for fixed $F$, the Fock module $\sF(F)$ and $\sB^a(\sF(E))$ and its tensor and Cuntz-Pimsner-Toeplitz subalgebras still depend on the choice of $\cC$. The corollary is, of course, true for all possible choices.

The condition that $F$ be a Hilbert module over some \nbd{C^*}subalgebra $\cC$ of $\cB$, means that $F$ is a closed subspace of $E$ invariant under the ternary product $(x,y,z)\mapsto x\AB{y,z}$; the minimal choice for $\cC$ is $\cC_F:=\cls\AB{F,F}$ (in which case $F$ is full) and every other choice must contain $\cC_F$ as an ideal. (See, for instance, the lemma in Skeide \cite[Section 0]{Ske18p1}.) It is easy to see that the smallest choice containing $A^*$, the closed ternary subspace generated by $A^*$, is
\baln{
\vspace{-1ex}
E_{A^*}
&
~:=~
\cls\bigcup_{n\in\N_0}A^*((A^*)^*A^*)^n;
&
\cB_{A^*}
&
~:=~
\cls\bigcup_{n\in\N}((A^*)^*A^*)^n.
}\ealn

\vspace{-1.5ex}\noindent
No smaller choice for $F$ and $\cC$ fitting the assumptions of Corollary \ref{smallcor} is possible. But is $E_{A^*}$ a correspondence over $\cB_{A^*}$? Or, more generally, if we have a closed ternary subspace $F$ of $E$ containing $A^*$ and a \nbd{C^*}subalgebra $\cC$ of $\cB$ containing $\cC_F$ as an ideal (so that $F$ is a Hilbert \nbd{\cC}module), is $F$ a \nbd{\cC}correspondence? This means actually two questions regarding the left action of $\cC$:
\begin{enumerate}
\item
Is $F$ invariant under $\cC$, that is, is $\cC F\subset F$?

\item
Does $\cC$ act nondegenerately on $F$, that is, is $\cls \cC F\supset F$?
\end{enumerate}
Both questions together may be united in the single question whether $\cls \cC F=F$; but we prefer to keep the two questions separate.

As far as the second question is concerned, this problem can be resolved once for all by passing to the unitalization $\wt{\cC}$ of $\cC$, provided the answer to the first question is affirmative. (Recall that $\cB$ is unital, so if $\U_\cB\notin\cC$, then by identifying the new unit $\wt{\U}_\cC$ with $\U_\cB$, $\wt{\cC}$ may be naturally identified as a unital subalgebra of $\cB$. This is independent on whether $\cC$ has its own unit $\U_\cC\ne\U_\cB$ or not.) Note that if we do so, then even if $F$ was a full Hilbert \nbd{\cC}module, it is now a definitely non-full Hilbert \nbd{\wt{\cC}}module. But, as explained in Footnote \ref{corrFN}, for us this is not a problem. (This also explains as simply as possible how and why, as claimed in Footnote \ref{corrFN}, the construction of $\sF(F)$ for degenerate left actions of $\cC$ works, too. Simply pass to $\wt{\cC}$ and construct $\sF(F)$ for the \nbd{\wt{\cC}}correspondence $F$. Then pass to $\sF(F)\sodot\cC=\cls\sF(F)\cC$, which removes from $\sF(F)$ the only (one-dimensional subspace spanned by the) element $\U_\cB\in\wt{\cC}=F^{\sodot 0}\subset\sF(F)$ that has inner products outside $\cC$. Corollary \ref{smallcor} remains true for $\sF(F)\sodot\cC$.) A case where nondegeneracy is clear, is when $\cC\ni\U_\cB$. It is easy to see that for interacting Fock spaces coming from subproduct systems (to be discussed in Section \ref{piSEC}) $\cB_{A^*}$ acts non-degenerately on $E_{A^*}$ if and only if the subproduct system is actually a product system (in which case the interacting Fock space is actually a full Fock space $\sF(H)$ and we really recover $\cB_{A^*}=\C$ and $E_{A^*}=\ol{H}$). Also if $\cI\ne H_0$ is nilpotent, then $\cB_{A^*}$ necessarily acts degenerately on $E_{A^*}$. (Indeed, since $H_{N+1}=\zero$, $A^*$ annihilates $H_N\ne\zero$, so none of the (even!) elements in $\cB_{A^*}$ can reach $H_N\backslash\zero$.)

So, after we have resolved (in an uncomplicated, pragmatic way) the second question (nondegeneracy), we are left with the first question (invariance). For $E_{A^*}$ and $\cB_{A^*}$, the only answer we can give is ``\it{rather no than yes}''; it depends highly on the interacting Fock space in question. Typical elements of $E_{A^*}$ are products or \hl{words} of elements or \hl{letters} that come aternatingly from $A^*$ and from $(A^*)^*$, starting and ending with a letter from $A^*$. The typical elements of $\cB_{A^*}$ are similar alternating words, but the first letter is from $(A^*)^*$ instead of $A^*$ (while the last one is still from $A^*$). If we multiply a word of  $E_{A^*}$ from the left with a word of $\cB_{A^*}$, then the  last letter of the latter (an element of $A^*$) meets the first letter of the former (also an element of $A^*$). So the resulting product word is no longer alternating. Whether or not it can be written as the limit of linear combinations of alternating words is totally unclear.
%\OW[Hast Du irgendeine Idee? Ansonsten laß ich's als Frage oder ich streich's.]{Answer for SubPS possible?}

With the notation $\ve_i=\pm1$ and, for $a^*\in A^*$, putting $a^1:=a^*$, $a^{-1}:=(a^*)^*$, one choice (usually) smaller than $\cB$, $E$ is
\begin{subequations} \label{BEcI}
\bal{ \label{BcI}
\cB_\cI
&
~:=~
\textstyle
\cls\BCB{a_n^{\ve_n}\ldots\,a_1^{\ve_1}\colon n\in\N,a_i^*\in A^*,\sum_{i=1}^n\ve_i=0}
\intertext{and} \label{EcI}
E_\cI
&
~:=~
\textstyle
\cls\BCB{a_n^{\ve_n}\ldots\,a_1^{\ve_1}\colon n\in\N,a_i^*\in A^*,\sum_{i=1}^n\ve_i=1}.
 }\eal
\end{subequations}
Clearly, $E_\cI$ is a full Hilbert \nbd{\cB_\cI}module. It is unclear if $\cB_\cI$ acts nondegenerately, but, clearly, it leaves $E_\cI$ invariant. If $\cB_\cI$ should act degenerately on $E_\cI$, then we would pass to $\wt{\cB}_\cI$ by adding to the generating set in \eqref{BcI} the term $a_n^{\ve_n}\ldots a_1^{\ve_1}=\U_\cB$ for $n=0$. (Modulo completion, this is the choice that has been discussed in \cite[Theorem 4.6]{AcSk08}.) Then $E_\cI$ is considered a (definitely non-full) correspondence over $\wt{\cB}_\cI$.

An even smaller choice, not discussed before, is
\begin{subequations} \label{BENC}
\bal{ \label{BNC}
\cB_\cI^{NC}
&
~:=~
\textstyle
\cls\BCB{a_n^{\ve_n}\ldots\,a_1^{\ve_1}\colon n\in\N,a_i^*\in A^*,\sum_{i=1}^k\ve_i\ge0\forall k\le n,\sum_{i=1}^n\ve_i=0}
\intertext{and} \label{ENC}
E_\cI^{NC}
&
~:=~
\textstyle
\cls\BCB{a_n^{\ve_n}\ldots\,a_1^{\ve_1}\colon n\in\N,a_i^*\in A^*,\sum_{i=1}^k\ve_i\ge0\forall k\le n,\sum_{i=1}^n\ve_i=1}.
}\eal
\end{subequations}
(\hl{NC} is referring to the fact that the difference of tuples occurring in \eqref{BcI} and \eqref{BNC} resembles the difference between \it{pair partitions} and \it{non-crossing pair partitions} of the set $\CB{1,\ldots,n}$ for even $n$.) Clearly, $\cB_\cI^{NC}$ is an algebra and $E_\cI^{NC}$ is invariant under left and right multiplication by elements of $\cB_\cI^{NC}$.

\bprop
$\cB_\cI^{NC}$ is a \nbd{C^*}algebra and the restriction of the inner product of $E$ turns $E_\cI^{NC}$ into a full Hilbert \nbd{\cB_\cI^{NC}}module.
\eprop

\proof
Suppose we have a word $a_n^{\ve_n}\ldots\,a_1^{\ve_1}$ from the generating set in \eqref{BNC}, that is, $\sum_{i=1}^k\ve_i\ge0\forall k\le n$ and $\sum_{i=1}^n\ve_i=0$. Then
\beqn{ \textstyle
\sum_{i=1}^k(-\ve_{n+1-i})
~=~
-\sum_{i=n-k+1}^n\ve_i
~=~
-(0-\sum_{i=1}^{n-k}\ve_i)
~\ge~0
}\eeqn
for all $k\le n$. Therefore, the word $(a_n^{\ve_n}\ldots\,a_1^{\ve_1})^*=a_1^{-\ve_1}\ldots\,a_n^{-\ve_n}$ is from the generating set, too. So, the Banach subalgebra $\cB_\cI^{NC}$ of $\cB$ is a \nbd{C^*}algebra.

In a similar way, one shows that $x,y\in E_\cI^{NC}$ implies $\AB{x,y}\in \cB_\cI^{NC}$. So, $E_\cI^{NC}$ is a Hilbert \nbd{\cB_\cI^{NC}}module.

Since every generating word $a_n^{\ve_n}\ldots\,a_1^{\ve_1}$ of $\cB_\cI^{NC}$ contains a factor of the form $a_{i+1}^-a_i^+$, the Hilbert \nbd{\cB_\cI^{NC}}module $E_\cI^{NC}$ is full.\qed

 \lf
Again, if $\cB_\cI^{NC}$ should act degenerately on $E_\cI^{NC}$, we may pass to the unitalization $\wt{\cB}_\cI^{NC}\ni\U_\cB$.

Summing up, we have presented three (usually) different ways to embed the Banach (\nbd{C^*})al\-ge\-bra $\ol{\cA^{(*)}}$ generated by $A^*$ into a tensor (Cuntz-Pimsner-Toeplitz) algebra. It is noteworthy that the latter (two) have no choice but containing $\cB_{A^*}$, which coincides with the Banach algebra generated by the set $(A^*)^*A^*$ and is a \nbd{C^*}algebra. It is usually not contained in $\cA$, so the containing tensor algebras will usually be bigger than $\ol{\cA}$.

\bex \label{CPTex}
Let $\cI=\Om\C\oplus H\oplus\Om_2\C$ for a pre-Hilbert space $H$ ($\dim H\ge2$) and some unit vector $\Om_2$, and assume $H$ has an anti-unitary involution $x\mapsto\bar{x}$. Turn $\cI$ into an interacting Fock space based on $H$ by defining $a^*(x)$ as
\baln{
\Om
&
~\longmapsto~
x,
&
y
&
~\longmapsto~
\Om_2\AB{\bar{x},y},
&
\Om_2
&
~\longmapsto~
0.
}\ealn
(The involution serves to assure that $x\mapsto a^*(x)$ is linear.) One easily checks that the adjoint $a(x)$ of $a^*(x)$ acts as
\baln{
\Om
&
~\longmapsto~
0,
&
y
&
~\longmapsto~
\Om\AB{x,y},
&
\Om_2
&
~\longmapsto~
\bar{x}.
}\ealn
We prefer to write these as finite-rank operators, getting $a^*(x)=x\Om^*+\Om_2\bar{x}^*$ and, consequently, $a(x)=\Om x^*+\bar{x}\Om_2^*$. Clearly, $a(x)$ leaves $\cI$ invariant, so $\cI$ with $A^*:=a^*(H)$ is an adjointable interacting Fock space with bounded creators.

For simplicity (in particular, notationally), we assume $H$ is a Hilbert space. We find $\cB=\C\oplus\sB(H)\oplus\C=\rtMatrix{\C&&\\&\sB(H)&\\&&\C}$ and $E=\sB_{(1)}(\cI)=\rtMatrix{&&\\H\Om^*&&\\&\Om_2H^*&{~~~}}$. From $a(x)a^*(y)=\Om\AB{x,y}\Om^*+\bar{x}\,\bar{y}^*$, we see that $\cB_{A^*}\subset\C\oplus\sK(H)\oplus 0$. Choosing $x\ne0\ne y$ perpendicular, we see that $\cB_{A^*}\ni a(x)a^*(y)=\bar{x}\bar{y}^*$; multiplying with other $a(x')a^*(y')$, it follows that $\cB_{A^*}$ contains all rank-one operators on $H$. Therefore, $\cB_{A^*}=\C\oplus\sK(H)\oplus 0$ and, consequently, $E_{A^*}=E$. From $E_{A^*}\subset F\subset E$ for any possible choice fulfilling the hypotheses of Corollary \ref{smallcor}, we find $E_\cI^{NC}=E_\cI=E$. In $\cB_\cI$ we find the word $a(x)^*a(x')^*a(y')a(y)=\Om_2\AB{\bar{x},x'}\AB{y',\bar{y}}\Om_2^*$, so $\cB_\cI=\C\oplus\sK(H)\oplus\C$. Now, the elements of $\cB_\cI^{NC}$ vanish on $\Om_2$ , so $\cB_{A^*}\subset\cB_\cI^{NC}\ne\cB_\cI$; we conclude that $\cB_\cI^{NC}=\cB_{A^*}$.

So, $E_{A^*}$, $E_\cI^{NC}$, and $E_\cI$ all coincide with $E=\sB_{(1)}(\cI)$ and $\cB_\cI^{NC}$ coincides with $\cB_{A^*}$, but the inclusions $\cB_{A^*}\subset\cB_\cI\subset\cB=\sB_{(0)}(\cI)$ are strict. Since $E$ is invariant under $\cB$, it is invariant under any subalgebra of $\cB$; and $\cB$ and $\cB_\cI$ act nondegenerately. However, since $\cB_{A^*}E=\rtMatrix{&&\\[-1ex]H\Om^*&&\\[-1ex]&\!0&{~~~~}}\ne E$, the action of $\cB_{A^*}$ is degenerate. So, we have to pass to the unitalization $\wt{\cB}_{A^*}=\cB_{A^*}+\id_\cI\C$. (Note that this does not coincide with $\C\oplus\wt{\sK(H)}\oplus\C=\wt{\cB}_\cI$; indeed, the latter contains $\id_H\in\wt{\sK(H)}$, while the former does not.) So, the tensor (Cuntz-Pimsner-Toeplitz) algebras into which we embed $\cA^{(*)}$ differ only by how much $\cB$ differs from $\cB_\cI$ differs from $\wt{\cB}_{A^*}$ (respectively, from $\cB_{A^*}$ if we do not insist in nondegenerate left actions) and the latter two are not contained in one another.

Going one step further to $\cI=\Om\C\oplus H_1\oplus H_2\oplus\Om_3\C$ with various choices for $A^*$, allows to produce more distinctive examples. (See also Open Problem 13 in Section \ref{2FockSS}.)
\eex

\lf
So far, we assumed an interacting Fock space $\cI$ with bounded $A^*$ that is adjointable. We briefly show how to free the preceding discussion and results from the hypothesis of adjointability.

So, we  now only assume that all elements of $A^*$ are bounded, but \bf{not necessarily adjointable}. (Of course, they are all weakly adjointable.) We may complete all pre-Hilbert spaces $\cI$ and $H_n$ and extend every element $a^*$ of $A^*$ to a (now adjointable) operator in $\sB\bfam{\,\ol{\cI}\,}$, which we continue denoting by $a^*$. (We do not assume that $\cI$ is based. In fact, completing $H$, wishing to extend also the map $H\rightarrow A^*$ involves unavoidably to change also $A^*$.) Clearly, $\cls A^*\ol{H}_n=\ol{H}_{n+1}$.

We also may immediately start with a family $\bfam{H_n}_{n\in\N_0}$ of Hilbert spaces where $H_0=\Om\C$, the Hilbert space $\cI=\coplus_{n\in\N_0}H_n$, and with a subset $A^*\subset\sB(\cI)$ such that
\beq{ \label{cls***}
\cls A^*H_n
~=~
H_{n+1}.
}\eeq
In this case, we may define the pre-Hilbert subspaces
\beqn{
\ul{H\!}\,_n
~:=~
\ls{A^*}^n\Om
}\eeqn
of $H_n$. Since all elements of $A^*$ are bounded, we may show by induction that $\ul{H\!}\,_n$ is dense in $H_n$ for all $n$. Clearly, elements of $A^*$ send $\ul{H\!}\,_n$ into $\ul{H\!}\,_{n+1}$. Therefore, the $H_n$ and $\cI$ may be obtained by the completion procedure described above, from the interacting Fock space $\ul{\cI\!}\,$ obtained from the $\ul{H\!}\,_n$ with the set $\ul{A}^*$ of all (co)restrictions of the elements of $A^*$ to operators on $\ul{\cI\!}\,$.

So, it does not really matter if we complete an interacting Fock space with bounded (but not necessarily adjointable) creators, or if we start start with a Hilbert-space-version of interacting Fock space where the axiom corresponding to \eqref{***} is replaced with the weaker condition in \eqref{cls***}. But, once we have Hilbert spaces, the elements of $A^*$ \bf{are} adjointable. It is clear that everything about $E$, $E_{A^*}$, $E_\cI^{NC}$, and $E_\cI$ (with the corresponding versions of $\cB$) goes through exactly, as before. We do not give details.

\brem
 We preferred not to mess up this section, which puts the application of \cite[Sections 3 and 4]{AcSk08} to the bounded case on firm ground, with too many references to Sections \ref{piSEC} and \ref{olSEC}. At least, we wish to emphasize again that this section together with Section \ref{piSEC} and its relation to the works \cite{ShaSo09,DRS11,KaSha15p} (explained in Section \ref{boundSS}) motivated this paper to large extent.
\erem

\newpage

\section{Boundedness: Criteria} \label{bcritSEC}

In the preceding section we have seen the nice consequences when $A^*$ has only bounded elements; in this section we wish to examine when the latter happens. Well, if we just have an (abstract) interacting Fock space, then we cannot do much more than just look at $A^*$ and check if its elements are bounded. What we mean is that in this section we will assume that $\cI$ is based on $H$ via the creator map $a^*\colon H\rightarrow A^*$ (so that there is $\Lambda$) or even embeddably based (so that there is $\vk$) or that it is regularly based (so that there is $L$). Recall that the first two things can be done for every interacting Fock space, while the last is limited to regular ones. We wish to understand boundedness of the creators in $A^*:=a^*(H)$ in terms of $\Lambda$, $\vk$, or $L$.

The question of boundedness has several layers. First of all, note that $a^*(x)$ is bounded if and only if all restrictions to the \nbd{n}particle sectors $H_n$ have finite norms $\norm{a^*(x)}_n:=\norm{a^*(x)\upharpoonright H_n}$ and if $\sup_n\norm{a^*(x)}_n$ ($=\norm{a^*(x)}$) is finite. (The same is true for $\Lambda$, $\vk$, $L$ ...) For being unbounded it is sufficient to show that $\norm{a^*(x)}_n=\infty$ for one $n$. On the other hand, if all $\norm{a^*(x)}_n$ are finite and $a^*(x)$ is unbounded just because the supremum is not finite, then this unboundedness is of a much nicer type. (For instance, the symmetric Fock space, that is, Example \ref{POIex}\eqref{POI1} for $q=1$, has creators of that type.) Such operators, clearly, have weak adjoints; their unboundedness is technically not more complicated than that of a selfadjoint operator with discrete spectrum. $a^*(x)$ that are unbounded on an \nbd{n}particle sector, may be arbitrarily irregular. All the criteria for boundedness in this section have (more or less obvious) versions for boundedness on each \nbd{n}particle sector (but not necessarily global), but we dispense with formulating them.

On the other hand, apart from the question whether $a^*(x)$ is bounded for every $x$, we may ask whether the creator map $a^*$ itself is bounded or not. This question, we will address.

Let us start with an example illustrating that even for a POI-interacting Fock space boundedness of the operator $L$ (or its \it{square root} $\Lambda$) does not guarantee boundedness of the creators $a^*(x)$.

\bex \label{bLunbex}
Let $H=L^2\SB{0,1}$ (as functions of $t\in\SB{0,1}$). For $L_1$ choose multiplication by $t$, for $L_2$ choose $\id_{H\otimes H}$, and put $L_n=0$ for $n\ge3$, so that $L$ is bounded. Then for $y_n=\I_{\SB{0,\frac{1}{n}}}$ we find $\norm{y_n}_\cI=\sqrt{\int_0^\frac{1}{n}t\,dt}=\frac{1}{\sqrt{2}n}$. For any $x\in H$ we find $\norm{a^*(x)y_n}_\cI=\norm{x}\sqrt{\frac{1}{n}}$, so,
\beqn{
\frac{\norm{a^*(x)y_n}_\cI}{\norm{y_n}_\cI}
~=~
\norm{x}\sqrt{2n},
}\eeqn
that is, despite $L$ is bounded, the operator $a^*(x)$ is unbounded whenever $x\ne0$.
\eex

We see, looking directly at boundedness of the operators $L$ or $\Lambda$ is not promising. So, in the sense of concluding from boundedness of `something' boundedness of all $a^*(x)$, the following obvious theorem in terms of $\vk$ is the best we can do.

\bthm \label{kappabthm}
Let $\cI=(H,\vk)$ be a \nbd{\vk}interacting Fock space. If $\vk$ is bounded, then the creator map $a^*$ is bounded by $\norm{a^*}\le\norm{\vk}$.
\ethm

\proof
$\norm{a^*(x)}=\norm{\vk\ell^*(x)}\le\norm{\vk}\,\norm{x}$.\qed

\lf
The condition that $\vk$ be bounded is not necessary. (See, however, Theorem \ref{prodSkthm}.)

\bex \label{kappaunbex}
Returning to Example \ref{CPTex}, we consider the interacting Fock space $\cI=\Om\C\oplus H\oplus\Om_2\C$ based on $H$ as embedded by choosing for $\Om_2$ a unit vector in $H\otimes H$. The norm of $a^*(x)$ is the norm of $x$, so the creator map $a^*$ is an isometry.

$\lambda_1(x)=a^*(x)\Om=x$, so $\vk_1=\lambda_1=\id_H$. For $\vk_2$ we compute $\lambda_2(x\otimes y)=a^*(x)a^*(y)\Om=\AB{\bar{x},y}\Om_2$, so
\beqn{
\vk_2(x\otimes y)
~=~
\vk_2(x\otimes\vk_1y)
~=~
\lambda_2(x\otimes y)
~=~
\Om_2\AB{\bar{x},y}.
}\eeqn
If $\dim H\ge\infty$, we may choose a self-adjoint orthonormal sequence $e_n=\bar{e}_n$. Since $\norm{\sum_{n=1}^N\frac{e_n\otimes e_n}{n}}^2$ $=\sum_{n=1}^N\frac{1}{n^2}$ converges, but $\sum_{n=1}^N\frac{\AB{e_n,e_n}}{n}=\sum_{n=1}^N\frac{1}{n}$ diverges, the map $\vk_2$, hence, $\vk$, is unbounded.

Note that $\lambda_2$ is not weakly adjointable. (The linear functional $\AB{\Om_2,\lambda_2\bullet}$ is unbounded, so there is no vector $Z=\lambda_2^*\Om_2\in\ol{H\otimes H}$ generating it as $\AB{Z,\bullet}$.) That is, $\cI$ is not regular. Note, too, that there is no difference if we assume $H$ is a Hilbert space. In Example \ref{bA*unbLex}, we will see a regular example.
\eex

The preceding example is based on (and an example for) the fact that the tensor product of Hilbert spaces does not share the usual universal property of tensor products: Not every bounded bilinear map $j\colon H\times H\rightarrow\C$ gives rise to a bounded linear map $\breve{j}\colon H\sbar{\otimes}H\rightarrow\C$ satisfying $\breve{j}(x\otimes y)=j(x,y)$. This gives the right idea. For boundedness of $a^*(x)$ or $a^*\colon x\mapsto a^*(x)$ not boundedness of $\vk$ is the relevant question, but boundedness of the bilinear map $(x,X)\mapsto\vk(x\otimes X)$. (We could replace the pre-Hilbert norm on $H\otimes\sF(H)$ with the projective norm on the tensor product, which has the universal property. But it would not give any better insight, so we dispense with this idea.) Keeping this in mind, the following improvement of Theorem \ref{kappabthm} is immediate.

\bthm\label{a*iffthm}
Let $\cI=(H,\vk)$ be a \nbd{\vk}interacting Fock space. Then:
\begin{enumerate}
\item
$a^*(x)$ is bounded if and only if there exists a constant $M_x$ such that $\norm{\vk(x\otimes X)}\le M_x\norm{X}$ for all $X\in\cD_\cI$.

\item
$a^*$ is bounded if and only if there exists a constant $M$ such that $\norm{\vk(x\otimes X)}\le M\norm{x}\norm{X}$ for all $x\in H$ and $X\in\cD_\cI$.
\end{enumerate}
\ethm

\noindent
Recalling the properties of $\vk$ and the interrelation of $\vk$ with $\lambda$, we observe that $\norm{\vk(x\otimes X)}$, for fixed $x$, takes its supremum varying over vectors of the form $\lambda X$ $(X\in\sF(H))$. The first condition transforms into
\beqn{
\norm{\lambda\ell^*(x)X}
~=~
\norm{\lambda(x\otimes X)}
~=~
\norm{\vk(x\otimes\lambda X)}
~\le~
M_x\norm{\lambda X},
}\eeqn
and analogously for the second condition. Recalling that the $\Lambda$ of an interacting Fock space $\cI$ embeddably based on $H$ is related to the $\lambda$, when we actually identify $\cI$ as a \nbd{\vk}interacting Fock via the embedding $\xi$, by $\lambda=\xi\Lambda$, we obtain the following criterion in terms of $\Lambda$, which is independent of how we actually embedded $\cI$. The nice thing is that (as the equation $\Lambda(x\otimes X)=a^*(x)\Lambda X$, which we used already so many times and which holds for arbitrary interacting Fock spaces based on $H$, shows) the inequalities expressed in terms of $\Lambda$ hold independently on whether $\cI$ is based embeddably or non-embeddably.

\bcor
Let $\cI$ be an interacting Fock space based on $H$. Then:
\begin{enumerate}
\item
$a^*(x)$ is bounded if and only if there exists a constant $M_x$ such that $\norm{\Lambda\ell^*(x)X}\le M_x\norm{\Lambda X}$ for all $X\in\sF(H)$.

\item
$a^*$ is bounded if and only if there exists a constant $M$ such that $\norm{\Lambda\ell^*(x)X}\le M\norm{x}\norm{\Lambda X}$ for all $x\in H$ and $X\in\sF(H)$.
\end{enumerate}
\ecor

\noindent
Now suppose $\cI$ is regular, that is, $\Lambda$ has a weak adjoint so that $L:=\Lambda^*\Lambda\ge0$ induces $\cI$ as POI-interacting Fock space. Then
\baln{
\norm{\Lambda\ell^*(x)X}^2
&
~=~
\AB{X,(\ell(x)L\ell^*(x))X},
&
\norm{\Lambda X}^2
&
~=~
\AB{X,LX}.
}\ealn
This allows, finally, to answer the long standing question, when a POI-interacting Fock space has bounded creators, in terms of operator inequalities.

\bthm \label{a*Lthm}
Let $\cI$ be a POI-interacting Fock space induced by the positive Fock operator $L\in\sL(\sF(H),\ol{\sF(H)})$. Then:
\begin{enumerate}
\item
$a^*(x)$ is bounded if and only if there exists a constant $M_x$ such that
\beqn{
\ell(x)L\ell^*(x)
~\le~
M_x^2L.
}\eeqn

\item
$a^*$ is bounded if and only if there exists a constant $M$ such that
\beqn{
\ell(x)L\ell^*(x)
~\le~
M^2\norm{x}^2L.
}\eeqn
\end{enumerate}
\ethm

\noindent
It is noteworthy that for the components $L_n$ of $L$, the inequalities read
\baln{ %\label{Lineq}
\ell(x)L_{n+1}\ell^*(x)
&
~\le~
M_x^2L_n,
&
\ell(x)L_{n+1}\ell^*(x)
&
~\le~
M^2\norm{x}^2L_n
}\ealn
(with $M_x$ and $M$, respectively, independent of $n$; in fact, if the constants exist, but depends on $n$, then this means the restrictions of $a^*(x)$ and $a^*$, respectively, to $H_n$ are bounded).

We know from Example \ref{bLunbex} that boundedness of $L$ is not sufficient for $L$ to fulfill the conditions in Theorem \ref{a*Lthm}. The following example shows that boundedness of $L$ is also not necessary.

\bex \label{bA*unbLex}
The construction of a counter example is based on the following computation. Denote by $e_1,\ldots,e_n$ the standard ONB of $\C^n$, and define the unit vector $e^n:=\sum_i\frac{e_i\otimes e_i}{\sqrt{n}}\in \C^n\otimes\C^n$. Then $\AB{e^n,x\otimes y}=\frac{1}{\sqrt{n}}\sum_ix_iy_i$. With the projection $p_n:=e^n{e^n}^*$, it follows that
\beqn{
\AB{x^n\otimes y^n,(np_n)(x^n\otimes y^n)}
~\le~
\norm{x^n}^2\norm{y^n}^2,
\text{~~~so,~~~}
(x^n\otimes\id_{\C^n})^*(np_n)(x^n\otimes\id_{\C^n})
~\le~
\norm{x^n}^2\id_{\C^n}.
}\eeqn
Consequently, if we define $H:=\bigoplus_{n\in\N}\C^n$ and the unbounded operator $L_2:=\bigoplus_{m,n\in\N}\delta_{m,n}np_n$ on $H\otimes H$, then for $x=\bigoplus_{n\in\N}x^n\in H$ we get
\bmun{
(x\otimes\id_H)^*L_2(x\otimes\id_H)
~=~
\bigoplus_{n\in\N}(x^n\otimes\id_{\C^n})^*(np_n)(x^n\otimes\id_{\C^n})
\\
~\le~
\bigoplus_{n\in\N}\norm{x^n}^2\id_{\C^n}
~\le~
\sup_{n\in\N}\norm{x^n}^2\bigoplus_{n\in\N}\id_{\C^n}
~\le~
\norm{x}^2\id_H.
}\emun
Therefore, putting $L_1:=\id_H$ and $L_n=0$ for $n\ge3$, we get a POI-interacting Fock space with bounded creator map but unbounded $L_2\le L$.
\eex

\lf
Let us collect the (non)implications we have in a diagram.
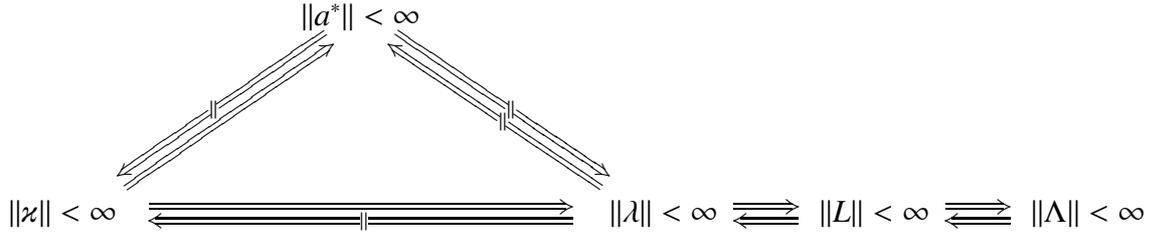
\beqn{
\xymatrix{
&&
\norm{a^*}<\infty	\ar@{=>}@<-.5ex>[ddll]|{||}	\ar@{=>}@<.5ex>[ddrr]|{||}
&&
\\\\
\norm{\vk}<\infty~~~\Big.		\ar@{=>}@<.5ex>[rrrr]	\ar@{=>}@<-.5ex>[uurr]
&&&&
\Big.~~~\norm{\lambda}<\infty~		\ar@{=>}@<.5ex>[llll]|{||}	\ar@{=>}@<.5ex>[uull]|{||}		\ar@{=>}@<.5ex>[r]
&~\norm{L}<\infty~		\ar@{=>}@<.5ex>[r]		\ar@{=>}@<.5ex>[l]
&~\norm{\Lambda}<\infty		\ar@{=>}@<.5ex>[l]
}
}\eeqn
The tail that starts from $\norm{\lambda}<\infty$ to the right, needs a comment. Clearly, a bounded $\Lambda$ is weakly adjointable, so there exists $L=\Lambda^*\Lambda$ and, necessarily, is bounded, too. And if $L$ exists, so that $\cI$ is embeddable, then also $\lambda$ exists (and is bounded, if $L$ is). If $\lambda$ exists (because we started with an embeddable interacting Fock space based on $H$), then $\lambda$ is just the $\Lambda$ for an isomorphic interacting Fock space based on $H$; again $\lambda$ bounded implies existence of $L$, which is bounded, too. So, bounded $\lambda$ and bounded $\Lambda$ are ``the same'', but only the situtation with $\Lambda$ is one that does not come along with an explicitly chosen embedding; and if $\Lambda$ is not bounded, then the situation is more general in that $\cI$ need not be embeddably based. So, it would add to the diagram if we made the same non-arrows which are there between $\lambda$ and $a^*$ also between $\Lambda$ and $a^*$. Last but not least, also the the non-arrow from $\norm{\lambda}<\infty$ to $\norm{\vk}<\infty$ requires a word; indeed if $\norm{\lambda}<\infty$ implied $\norm{\vk}<\infty$, then together with the arrow from $\norm{\vk}<\infty$ to $\norm{a^*}<\infty$ we would get the arrow from $\norm{\lambda}<\infty$ to $\norm{a^*}<\infty$, which , as we know, is not true.

\newpage

\section{Subproduct systems: A class of examples} \label{piSEC}

A class of operator algebras ($*$ or not) generated by creators on Fock type spaces arises from so-called \it{subproduct systems}. Subproduct systems (even of correspondences) have been introduced by Shalit and Solel \cite{ShaSo09} and, independently, (under the name of \it{inclusion systems} and limited to Hilbert spaces) by Bhat and Mukherjee \cite{BhMu10}. The operator algebras of our interest in this paper, have been introduced by Shalit and Solel \cite{ShaSo09} and led to several forthcoming papers by Shalit and his collaborators; see also Section \ref{boundSS}. During the 2011 Spring School and Conference on ``Product Systems and Independence in Quantum Dynamics'' in Greifswald, when listening to Shalit's talk, several participants noted instantaneously, that the Fock type spaces of subproduct systems are interacting Fock spaces; this also includes the same set of creators in a canonical basing.

The scope of this section is to examine the structure of these interacting Fock spaces arising from subproduct systems (namely, \nbd{\vk}interacting Fock spaces, where $\vk=\pi$ is a projection that, apart from being a squeezing, fulfills an extra condition). On the fly, we examine the general structure of \nbd{\vk}interacting Fock spaces, where $\vk=\pi$ is a projection.

\lf
A (\hl{discrete}) \hl{subproduct system} (of Hilbert spaces) is a family $H^\botimes=\bfam{H_n}_{n\in\N_0}$ of Hilbert spaces $H_n$ with isometric \hl{coproduct} maps $w_{m,n}\colon H_{m+n}\rightarrow H_m\sbar{\otimes}H_n$ iterating coassociatively, and with $H_0=\C$ such that the marginal maps $v_{n,0},v_{0,n}$ become the canonical identifications $H_n\otimes\C\cong H_n\cong\C\otimes H_n$. (In several places, there occurred also \it{superproduct systems}, replacing the isometries with coisometries. A far reaching generalization of both (arising in the dilation theory of multi-parameter CP-semigroups in Shalit and Skeide \cite{ShaSk10p}) we discuss in Section \ref{pivesSS}.

For our purposes, it is better to pass to the \hl{product} maps $v_{m,n}:=w_{m,n}^*\colon H_m\sbar{\otimes}H_n\rightarrow H_{m+n}$, which are coisometries. The associativity condition, then, really means that the \hl{product} $(x_m,y_n)\mapsto x_my_n:=v_{m,n}(x_m\otimes y_n)$ is associative.

If $H^\botimes$ is a subproduct system, then the \hl{Fock space} over $H^\botimes$ is $\sF(H^\botimes):=\coplus_{n\in\N_0}H_n$. For each $x\in H_1$, we define the creator $a^*(x)\in\sB(\sF(H^\botimes))$ by $a^*(x)X_n:=xX_n$ for all $n,X_n\in H_n$; see, for instance, \cite{ShaSo09}.

Since $v_{1,n}$ is a coisometry, it is surjective. More precisely, it maps the Hilbert space $H_1\sbar{\otimes}H_n$ onto the Hilbert space $H_{n+1}$. If we take only the algebraic tensor product $H_1\otimes H_n$, then (as soon as $H_{n+1}$ is not finite-dimensional) it is no longer surjective, but only with dense range. So, thinking of $\sF(H^\botimes)$ as an interacting Fock space (writing also $H_0=\Om\C$ with $\Om=1\in\C=H_0$), we are in the situation sketched in the end of Section \ref{boundSEC}, where \eqref{***} is replaced by the weaker \eqref{cls***}. As explained there, we know how to pass to the proper interacting Fock space $\ul{\sF(H^\botimes)}:=\bigoplus_{n\in\N_0}\ul{H\!}\,_n$ determined by the family of dense pre-Hilbert subspaces
\beqn{
\ul{H\!}\,_n
~:=~
\ls{a^*(H_1)}^n\Om
~\subset~
H_n.
}\eeqn

Roughly, we started with a subproduct system (that is, by definition) of Hilbert spaces and obtained the topological version of interacting Fock space as discussed in the end of Section \ref{boundSEC}. The reduction, there, to a proper interacting Fock space (with only bounded creators) can be interpreted, in the context of subproduct systems, as the passage to the \hl{algebraic subproduct system} of (dense) pre-Hilbert (sub)spaces and their algebraic tensor products \hl{generated} by $H_1$. The fact that this, actually, \bf{is} the algebraic subproduct system generated by $H_1$, follows clearly from writing the structure with (coisometric) product maps. Indeed, the $n$th pre-Hilbert space is just what is spanned by \nbd{n}fold products of elements from $H_1$; it is clear by construction that the iterated products $v_{m,n}$ leave these algebraic domains invariant. If we insisted to work with the (isometric) coproduct maps $w_{m,n}$, then it would not at all be clear if we could find \bf{dense} pre-Hilbert subspaces so that the restriction of $w_{m,n}$ would map into their algebraic tensor product. (This is \it{a priori} not even clear for $w_{1,1}$. But, while for the products $v_{m,n}$ the problem is solved inductively, here, for the coproducts $w_{m,n}$ no inductive solution is possible, because with each new level $N+1$, the possible solution for $n,m\le N$ will be affected; this situation is ``anti-inductive''.) For this reason, the following observation, which tells that by the co/isometric property we actually do obtain an algebraic subproduct system $\bfam{\ul{H\!}\,_n}_{n\in\N_0}$ with respect to the (co)restricted coproduct maps $w_{m,n}$, is quite remarkable:

\bob \label{preob}
Suppose we have (pre-)Hilbert spaces $H\supset H'$ and $G\supset G'$, and suppose we have a (necessarily adjointable) coisometry $w\colon H\rightarrow G$ that (co)restricts to a surjective map $w'\colon H'\rightarrow G'$. Then the adjoint $w^*$ of $w$ (co)restricts, too, to a map $G'\rightarrow H'$, necessarily the adjoint of $w'$. (Indeed, by replacing $H$ with the range of the projection $w^*w$ (so that, in particular, surely $w^*$ maps $G$ into that space no matter how small or big the subspace $G'$ is), we may assume that $w$ is actually unitary. Then, like for every invertible map, the restriction of the inverse map $w^*$ to the image $G'$ of a restriction of the map $w$ to $H'$, sends $G'$ into (hence, onto) $H'$. If we add again what we cut away to make $w$ unitary, we see that $w^*$ maps $G'$ onto $H'\cap(w^*wH)$. Of course, $ww^*$ (co)restricts to $\id_{G'}$; the only question was if the first map $w^*$ of the product $ww^*$ does lead or does not lead out of $H'$.) Consequently, the (coisometric!) product maps of the algebraic subproduct system $\bfam{\ul{H\!}\,_n}_{n\in\N_0}$ have (isometric) adjoints for the algebraic (co)domains. Therefore, while in the general case considered in the end of Section \ref{boundSEC} the restrictions of the creators to dense interacting Fock space need not be adjointable, in our case here the (co)restrictions of the $a^*(x)$ remain adjointable. (Indeed, $a^*(x)$, on the algebraic domain, is adjointable if and only if each $a^*(x)\upharpoonright\ul{H\!}\,_n$ (considered as map into $\ul{H\!}\,_{n+1}$) is adjointable, and $a^*(x)\upharpoonright\ul{H\!}\,_n=w_{1,n}(x\otimes\id_{\ul{H\!}\,_n})$ has an adjoint, namely, $(x\otimes\id_{\ul{H\!}\,_n})^*v_{1,n}$.) Therefore, the (proper) interacting Fock space of a subproduct system $H^\botimes$ is adjointable.
\eob

We now wish to understand the structure of interacting Fock spaces derived from subproduct systems. More precisely, we wish to understand them as \nbd{\vk}interacting Fock spaces, and distinguish those $\vk$ that lead to interacting Fock spaces coming from subproduct systems. The following result shows that not only the basing is embeddably, but that there is actually a very canonical embedding into $\sF(H_1)$.

Here and in the sequel, we denote by $v_{n_1,\ldots,n_k}\colon H_{n_1}\sbar{\otimes}\ldots\sbar{\otimes}H_{n_k}\rightarrow H_{n_1+\ldots+n_k}$ the iterated product of $k$ factors (which, by associativity, does not depend on how we iterate), and we denote the special case of $n$ factors from $H_1$ as $v_{(n)}:=v_{1,\ldots,1}$.

\bthm \label{SubPsubSthm}
\begin{enumerate}
\item \label{SPSS1}
Suppose $H$ is a Hilbert space and $\pi_n$ are projections in $\sB(H^{\bar{\otimes}n})$ (with $\pi_0=\id_\C$). Then the maps $v_{m,n}\colon(\pi_mX_m)\otimes(\pi_nY_n)\mapsto \pi_{m+n}(X_m\otimes Y_n)$ turn the family $\bfam{\pi_nH^{\bar{\otimes}n}}_{n\in\N_0}$ into a subproduct system if and only if the projections $\pi_n$ satisfy
\beq{ \label{pnSubPS}
\id_H\otimes \pi_n
~\ge~
\pi_{n+1}
~\le~
\pi_n\otimes\id_H
}\eeq
for all $n\in\N$.

\item \label{SPSS2}
Suppose $H^\botimes$ is a subproduct system, and put $\pi_n:=v_{(n)}^*v_{(n)}\in\sB(H_1^{\sbar{\otimes}n})$. Then the $\pi_n$ fulfill \eqref{pnSubPS} and
\beqn{
X_n
~\longmapsto~
v_{(n)}^*X_n
}\eeqn
is an isomorphism of subproduct systems from $H^\botimes$ to $\bfam{\pi_nH_1^{\bar{\otimes}n}}_{n\in\N_0}$.
\end{enumerate}
\ethm

\proof
\ref{SPSS1}. Associativity is manifest, once the $v_{m,n}$ are well-defined. It is clear that $v_{m,n}$ is well-defined if and only if the kernel of $\pi_m\otimes \pi_n$ is contained in the kernel of $\pi_{m+n}$, that is, if and only if
\beq{ \label{pmnSubPS}
\pi_m\otimes \pi_n
~\ge~
\pi_{m+n}.
}\eeq
What remains is to show that the necessary conditions in \eqref{pnSubPS} (they form a subset of the conditions in \eqref{pmnSubPS}) are also sufficient. Note that \eqref{pnSubPS} may also be written as $(\id_H\otimes \pi_n)\pi_{n+1}=\pi_{n+1}=\pi_{n+1}(\pi_n\otimes\id_H)$. We find
\bmun{
\pi_{m+n}
~=~
(\id_H\otimes \pi_{m-1+n})\pi_{m+n}
~=~
(\id_{H^{\bar{\otimes}2}}\otimes \pi_{m-2+n})(\id_H\otimes \pi_{m-1+n})\pi_{m+n}
\\
~=~
(\id_{H^{\bar{\otimes}2}}\otimes \pi_{m-2+n})\pi_{m+n}
~=~
\ldots
~=~
(\id_{H^{\bar{\otimes}m}}\otimes \pi_n)\pi_{m+n},
}\emun
that is, $\id_{H^{\bar{\otimes}m}}\otimes \pi_n\ge \pi_{m+n}$, and, similarly, $\pi_{m+n}=\pi_{m+n}(\pi_m\otimes\id_{H^{\bar{\otimes}n}})$, that is, $\pi_m\otimes\id_{H^{\bar{\otimes}n}}\ge \pi_{m+n}$. Both together give \eqref{pmnSubPS}.

\ref{SPSS2}. Clearly, $v_{(n)}^*$, being an isometry, defines a unitary onto $v_{(n)}^*H_n=\pi_n H_1^{\bar{\otimes}n}$. By the family $v_{(n)}^*$ of unitaries, the product maps $v_{m,n}$ lift to a product on the family $\bfam{\pi_nH_1^{\bar{\otimes}n}}_{n\in\N_0}$ as
\bmun{
(\pi_mX_m)\otimes(\pi_nY_n)
~\longmapsto~
v_{(m)}(\pi_mX_m)\otimes v_{(n)}(\pi_nY_n)
~=~
v_{(m)}X_m\otimes v_{(n)}Y_n
\\
~\longmapsto~
v_{m,n}(v_{(m)}X_m\otimes v_{(n)}Y_n)
~=~
v_{(m+n)}(X_m\otimes Y_n)
\\
~\longmapsto~
v_{(m+n)}^*v_{(m+n)}(X_m\otimes Y_n)
~=~
\pi_{(m+n)}(X_m\otimes Y_n)
}\emun
(first sending the elements $\pi_mX_m$ and $\pi_nY_n$ of the family $\bfam{\pi_nH_1^{\bar{\otimes}n}}_{n\in\N_0}$ to the family $H^\botimes$ where, then, $v_{m,n}$ is applied to send, in the end, the result $v_{(m+n)}(X_m\otimes Y_n)$ back to $\bfam{\pi_nH_1^{\bar{\otimes}n}}_{n\in\N_0}$). This is not only precisely the action we wish to define in Part \ref{SPSS1}. It also establishes the latter, being an image of the subproduct system structure of $H^\botimes$, as a properly defined operation of a subproduct system, therefore, necessarily satisfying \eqref{pnSubPS}. By construction, the family of unitaries $v_{(n)}^*$ is an isomorphism of subproduct systems.\qed

\brem
Using the conditions in \eqref{pmnSubPS}, this is just a suitably reformulated version of \cite[Lemma 6.1]{ShaSo09}, referring to the family $\bfam{\pi_nH^{\bar{\otimes}n}}_{n\in\N_0}$ as a \hl{standard subproduct system}. That the weaker conditions in \eqref{pnSubPS} already suffice, is new. These conditions are modeled after and motivated by an analogous set of combinatorial conditions in the combinatorics of words systems and their associated subproduct systems, discussed in Gerhold and Skeide \cite{GeSk14p}.
\erem

Recall that by Observation \ref{preob}, $v_{(n)}^*$ maps $\ul{H\!}\,_n$ really into the algebraic tensor power $H_1^{\otimes n}$. Therefore, $\pi_n$ (co)restricts to a projection in $\sB^a(H_1^{\otimes n})$, which we continue denoting $\pi_n$. Their direct sum $\pi$ is a Fock projection in $\sB^a(\sF(H_1))$. If we define $\xi:=\bigoplus_{n\in\N_0}v_{(n)}^*$, then we embed the interacting Fock space $\cI:=\ul{\sF(H^\botimes)}$ onto
\beqn{
\xi\cI
~=~
\pi\sF(H_1)
~\subset~
\pi\ol{\sF(H_1)}
~\subset~
\ol{\sF(H_1)}.
}\eeqn
By definition $\xi\cI$ is a subspace of the completion $\ol{\sF(H_1)}$ and the complement $(\xi\cI)^\perp$ is relative to that Hilbert space. But thanks to being the range of the projection $\pi\in\sB^a(\sF(H_1))$, the subspace $\xi\cI$ is complemented also in $\sF(H_1)$. (The complement in this space is just the intersection of the topological complement $(\xi\cI)^\perp$ with $\sF(H_1)$.) Then, $\pi$ is literally everything we can know about that embedded interacting Fock space: $\pi=L=\lambda=\vk$. (Indeed, clearly, $\lambda_n=\pi_n$, so $L_n=\lambda_n^*\lambda_n=\pi_n$. Clearly, inserting $\pi_n$ as candidate for $\vk_n$ into the recursion for $\lambda_n$, we recover $\lambda_n=\pi_n=\vk_n$. For being the (uniquely determined) squeezing $\vk$, the resulting Fock projection $\pi$, with which we wish to identify $\vk$, has to be a squeezing. But, also this is true, because clearly $\pi_{n+1}$ is surjective, and since $\pi_{n+1}\le\id_H\otimes\pi_n$, we get that $\pi_{n+1}$ is $0$ on $H\otimes H_n^\perp$. We see, how nicely the algebraic invariance properties discussed in Observation \ref{preob} in the case of interacting Fock spaces from subproduct systems work together with the more topological definitions of \nbd{\pi}interacting Fock space.)

The squeezing $\pi$ is a projection. We ask what other properties a squeezing has to satisfy to be the one that comes from a subproduct system as described. This question requires also to understand which Fock projections are squeezings. Actually, we first need a sufficiently flexible notion of projection. We say, a map $\pi$ from a pre-Hilbert space $H$ into its completion $\ol{H}$ is a \hl{weak} projection if $\AB{x,\pi y}=\AB{\pi x,\pi y}$ for all $x,y\in H$.

A weak projection extends uniquely to a projection in $\sB(\ol{H})$, and every restriction of a projection in $\sB(\ol{H})$ to $H$ is a weak projection. A weak projection need not be a projection.

\bex \label{wpnonpex}
Consider the (completed) one-mode full Fock space $\ol{\sF(\C)}$ with the \hl{exponential vectors} $\ee(z):=\sum_{n\in\N_0}\frac{z^ne_n}{\sqrt{n!}}$ for all $z\in\C$. Put $H:=\ls\CB{\ee(z)\colon z\ne0}$. Since the set of all exponential vectors (including $\ee(0)=\Om$) is linearly independent, the projection $\Om\Om^*\in\sB(\ol{H})$, when restricted to $H$ does not leave $H$ invariant. It is, therefore, a weak projection that is not a projection.
\eex

For simplicity, in the following theorem we assume $\pi_1=\id_H$, identifying this way $H_1$ with $H$ (otherwise being only a subspace of $H$). One can show that we always may replace $H$ with $H_1:=\pi_1 H$.

\lf
\bthm
Let $\cI=(H,\pi)$ be a \nbd{\pi}interacting Fock space where the squeezing $\pi\colon(H\otimes\cI)\oplus\Om\C\rightarrow\cI\subset\ol{\cI}$ is a weak projection with $\pi_1=\id_H$. Then
\beq{ \label{pirec}
\pi_{n+1}
~\le~
\id_H\otimes\pi_n.
}\eeq
Conversely, if $H$ is a pre-Hilbert space and $\pi\in\sB(\ol{\sF(H)})$ a weak Fock projection such that the components $\pi_n$ fulfill \eqref{pirec} (and $\pi_1=\id_H$), then $\cI:=\pi\sF(H)$ is a \nbd{\pi}interacting Fock space.

Moreover, in either case among the summands $\pi_nH^{\otimes n}$ there exist coisometries $\pi_mH^{\otimes m}\otimes\pi_n H^{\otimes n}\rightarrow\pi_{m+n}H^{\otimes m+n}$ satisfying $\pi_mX_m\otimes\pi_nY_n\mapsto\pi_{m+n}(X_m\otimes Y_n)$ (so that the $\pi_n\ol{H^{\otimes n}}$ form a subproduct system and $\cI$ is its associated interacting Fock space) if and only the $\pi_n$ also fulfill
\beq{ \label{SPSker}
\pi_{n+1}
~\le~
\pi_n\otimes\id_H.
}\eeq
\ethm

\proof
As discussed two paragraphs before the theorem, if $\pi$ is a squeezing, then the condition \eqref{pirec} is fulfilled. On the other hand, if a $\pi$ is Fock projection in $\sB(\ol{\sF(H)})$, then by definition $\pi$ sends $\sF(H)$ surjectively onto $\cI$, and if $\pi$ fulfills \eqref{pirec}, then $\pi_{n+1}$ is $0$ on $H\otimes(\pi_nH^{\otimes n})^\perp$, so $\pi$ is a squeezing. We argued already that the last statement is true.\qed

\lf
It is noteworthy that, despite Example \ref{wpnonpex}, the two inequalities together imply the algebraic invariance discussed in Observation \ref{preob}.

\bex
There are \nbd{\pi}interacting Fock spaces that do not come from a subproduct system. Let $H$ be a pre-Hilbert space with an orthonormal Hamel basis $\bfam{e_n}_{n\in\N}$ and put $p_n=e_ne_n^*$. Then $\pi_n=p_n\otimes\ldots\otimes p_1$ define a squeezing $\pi$ that does not satisfy \eqref{SPSker}.
\eex

\bob
By \eqref{pirec} and Theorem \ref{a*Lthm}, a \nbd{\pi}interacting Fock space is a POI-interacting Fock space with bounded creator map $a^*$.
\eob

\newpage

\section{Outlook} \label{olSEC}

\renewcommand{\thesubsection}{\thesection.\Alph{subsection}}

As outlined in the introduction, this section presents results, (counter) examples, and considerations, pointing to future work. While we present a considerable number of problems we leave open, the results serve to be able to formulate the problems or to provide (counter) examples, and the (counter) examples serve to motivate problems or to illustrate why the answer to a problem is too involved to be included in the present paper. Especially the counter examples, also have the scope to prevent ourselves and the reader to formulate problems that are tempting but turn out to be useless to consider.

We split the discussion into six parts, \it{regularity} \eqref{regSS}, \it{bounded creators} \eqref{boundSS}, \it{productive systems} \eqref{pivesSS}, \it{classification} \eqref{classSS}, a \it{case study} \ref{2FockSS}, and \it{automorphism groups} \ref{autoSS}. The subdivision is rather \it{ad hoc}, the order may appear a bit arbitrary. In the end, all these subsections go into the direction of questions about classification. But while the first three go by limiting to subclasses (\ref{boundSS} and \ref{pivesSS} opening up relations to other areas), \ref{classSS} goes into the direction of how to tackle, meaningfully, the question of general classification. In \ref{autoSS}, we address questions about several automorphism groups of interacting Fock spaces.

% In \ref{regSS} we ask, specifically, about regularity of \nbd{\vk}interacting Fock spaces. We show that regularity may not related directly to a ``nice'' property of $\vk$, leaving as an open question, which are the squeezings $\vk$ that lead to regular \nbd{\vk}interacting Fock spaces. Not even the hypothesis of bounded creators helps. This leads to the question, which bounded interacting Fock spaces are regular. Asking this question for abstract bounded interacting Fock spaces, leads us to the situation in \ref{boundSS}. There, we analyze the question about the classification of bounded interacting Fock spaces in terms of the associated operator algebras. This leads directly to \ref{pivesSS}, where we ask for the interacting Fock spaces that come from \it{nondegenerate productive systems} introduced recently by Shalit and Skeide \cite{ShaSk10p}. \it{Productive systems}, generalizing both subproduct systems and superproduct systems, when \it{nondegenerate}, lead to a well specified subclass of \nbd{\vk}interacting Fock spaces, containing those coming from subproduct systems, which we have discussed on Section \ref{piSEC}. It is natural to ask the questions from \ref{boundSS} restricted to this subclass.

% \newpage

\subsection{About regularity} \label{regSS}

Well, the title of this part would be better \it{about regularity in terms of $\vk$}. Can we determine in terms of $\vk$ regularity of an embedded interacting Fock space? The answer is in form of no-go-theorems (basically, Corollary \ref{padcor}).

Every interacting Fock space can be embedded. More precisely, every interacting Fock space can be embeddably based (Theorem \ref{aIFSembthm}) and, then (after having it based embeddably), embedded. By Theorem \ref{**thm}, every embedded interacting Fock space can be viewed as a \nbd{\vk}interacting Fock space in a unique way. By Corollary \ref{lLcor}, this \nbd{\vk}interacting Fock space is regularly based if and only $\lambda$ is weakly adjointable. The question we tackle here, is if we may hope that adjointability of $\lambda$ is related in some useful way to adjointability of $\vk$.

$\lambda$ and $\vk$ are related by \eqref{lkeq} which amounts, equivalently, to the recursion
\beqn{
\lambda_{n+1}
~=~
\vk_{n+1}(\id_H\otimes\lambda_n)
\text{~~~and~~~}
\lambda_0
~=~
\id_{\C\Om},
}\eeqn
which we repeat here for convenience. So, if $\lambda_n$ and $\lambda_{n+1}$ have an adjoint (happening for all $n$ if and only if $\lambda$ has an adjoint), is this enough to force that $\vk_{n+1}$ has an adjoint (happening for all $n$ if and only if $\vk$ has an adjoint)? If $\vk_{n+1}$ and $\lambda_n$ have an adjoint, is this enough to force that also $\lambda_{n+1}$ has an adjoint? (Since $\lambda_n$ is computed recursively, this just means whether or not $\vk$ adjointable implies $\lambda$ adjointable.) We are, roughly, concerned with the following situation.

Let $G$, $H$, and $K$ be pre-Hilbert spaces and let $a\colon G\rightarrow H$ and $b\colon H\rightarrow K$ be linear operators. Put $c:=ba\colon G\rightarrow K$. Does weak adjointability of two of them imply weak adjointability of the third one?

Let us recall that a densely defined operator $a\colon\ol{G}\supset\cD_a\rightarrow\ol{H}$ is \hl{closeable} if the closure of its graph $\cG_a:=\CB{(g,ag)\colon g\in\cD_a}$ in $\ol{G}\oplus\ol{H}$ is the graph of a (densely defined, since $\cD_{\ol{a}}\supset\cD_a$) operator $\ol{a}\colon\cD_{\ol{a}}\rightarrow\ol{H}$, the \hl{closure} of $a$. (This happens if and only if for each sequence $g_n$ in $\cD_G$ with $g_n\to 0$ we have that $ag_n\to h\in\ol{H}$ implies $h=0$.) In case $a=\ol{a}$, we say $a$ is \hl{closed}. A closeable operator has a densely defined adjoint $\cD_{a^*}\rightarrow\ol{G}$ (namely, the operator whose graph is $\f(\cG_{-a}^\perp)$ where $\f\colon(g,h)\mapsto(h,g)$ is the \hl{flip}). But it is not said that the (maximal!) domain $\cD_{a^*}\subset\ol{H}$ contains $H$. However, any weakly adjointable operator is closeable. So, weak adjointability is stronger a property than closeability.

\bex \label{padex}
Put $G:=H:=\ls\CB{e_n}$ for some orthonormal family $\bfam{e_n}_{n\in\N}$, and put $K:=\C$.
\begin{enumerate}
\item \label{pex1}
$a,b$ closeable/adjointable $\not\Rightarrow$ $c$ closeable/adjointable.

Let $ae_n:=ne_n$ and $be_n:=\frac{1}{n}$. Then $a$ is adjointable (in fact, $a$ is selfadjoint) and $b$ is adjointable (in fact, $b$ is bounded), but $c\colon e_n\mapsto 1$ is unbounded, hence, not closeable, \it{a fortiori} not adjointable. (A densely defined adjoint $\C\supset\cD_{c^*}\rightarrow\ol{H}$, being an operator with finite-dimensional domain, is necessarily bounded, which implied that $c$ itself had to be bounded.)

\item \label{pex2}
$a,c$ closeable/adjointable $\not\Rightarrow$ $b$ closeable/adjointable.

Let $ae_n:= \frac{e_n}{n}$ and $be_n:=1$. Then $a$ is adjointable (in fact, $a$ is bounded) and $c\colon e_n\mapsto\frac{1}{n}$ is adjointable (in fact, $c$ is bounded), but (like $c$ in Number \ref{pex1}) $b$ is not closeable, \it{a fortiori} not adjointable.

\item \label{pex3}
We add also the last case, $b,c$ closeable/adjointable $\not\Rightarrow$ $a$ closeable/adjointable (which we do not need).

Let $ae_n:= e_n-ne_1$ and $be_n:=\frac{1-\delta_{1,n}}{n}$. Then $c=b$ are bounded, hence, adjointable, but $a$ is not closeable, hence, not adjointable. (Indeed, the sequence $\frac{e_n}{n}$ converges to $0$, while $a\frac{e_n}{n}=\frac{e_n}{n}-e_1$ converges to $-e_1\ne0$, showing $a$ is not closeable.)
\end{enumerate}
\eex

\bcor \label{padcor}
Neither does weak adjointability of $\vk$ imply regularity, nor does regularity imply weak adjointability of $\vk$.
\ecor

\proof
For the overall setting as in Example \ref{padex}, define $\vk_1:=a$, $\vk_2:=\id_H\otimes e_1b$, and $\vk_n=0$ for all $n\ge3$. Then (no matter which of the possibilities for $a$ and $b$ we choose) $\vk$ with components $\vk_n$ is a squeezing relative to $\cI=\Om\C\oplus H\oplus(H\otimes e_1)\oplus 0\ldots\subset\sF(H)$. Moreover, $\lambda_1=\vk_1=a$ and $\lambda_2=\vk_2(\id_H\otimes\lambda_1)=\id_H\otimes e_1c$. The first statement of the corollary follows from the choice in \ref{padex}\eqref{pex1}, the second statement follows from the choice in \ref{padex}\eqref{pex2}.\qed

\bOP
What are the squeezings $\vk$ that lead to regular \nbd{\vk}interacting Fock spaces? By Corollary \ref{Lam*ALVcor}, the regular interacting Fock spaces are exactly the POI-interacting Fock spaces. So one might try to approach that problem, starting directly from POI-interacting Fock spaces and see if it is possible to specify the special properties of their $\vk$. (Recall that they are not only embeddable by Corollary \ref{POIembcor}, but that the embedding constructed for that goal in the proof of Lemma \ref{posdimlem} is actually quite canonical.)
\eOP

It might be worthwhile to look at other properties the \nbd{\vk}interacting Fock spaces constructed in the proof of Corollary \ref{padcor} from (all three cases of) Example \ref{padex} have. (For instance, the properties of the creators $a^*(x)=\vk\ell^*(x)$ depend entirely on the corresponding properties of $b$.) We omit this, but we ask:

\bOP
What are the squeezings $\vk$ that lead to adjointable \nbd{\vk}interacting Fock spaces, or at least to \nbd{\vk}interacting Fock spaces with closeable creators?
\eOP

We have already characterized the squeezings $\vk$ that lead to bounded creators or even bounded creator maps in Theorem \ref{a*iffthm}. Let us recall that by the discussion in the end of Section \ref{bcritSEC}, bounded $\vk$ implies that everything else is bounded (and, therefore, weakly adjointable). Finally, recall that for regularity, boundedness is neither sufficient (see Example \ref{kappaunbex}) nor necessary (see symmetric Fock space in Example \ref{POIex}\eqref{POI1} for $q=1$).

With the last question of this subsection, we leave the situation of a given \nbd{\vk}interacting Fock space, or even of based interacting Fock spaces, but pass to abstract ones. This points straight at \ref{boundSS} and \ref{classSS}.

\bOP
Which bounded interacting Fock spaces $\cI=(\bfam{H_n}_{n\in\N_0},A^*)$ are regular?
\eOP

\newpage

\subsection{About bounded creators} \label{boundSS}

We can say that the situation of a bounded interacting Fock space $\cI=(\bfam{H_n}_{n\in\N_0},A^*)$, with the consequences in Section \ref{boundSEC} and the results by Davidson, Ramsey, and Shalit \cite{DRS11}, were what motivated this paper. With a bounded interacting Fock space $\cI$ we associate the Banach (\nbd{C^*})algebra $\ol{\cA^{(*)}}$ generated by $A^*$. We ask:

\bOP
To what extent bounded interacting Fock spaces are classified by their associated operator algebras $\ol{\cA^{(*)}}$?
\eOP

Davidson, Ramsey, and Shalit \cite{DRS11} have analyzed the operator algebras $\ol{\cA^{(*)}}$ for interacting Fock spaces coming from a subclass of the finite-dimensional subproduct systems. Since $H_1$ is \nbd{d}dimensional and since $H_n$ may be thought of as image of a projection in $\sB(H_1^{\sbar{\otimes}n})$, as already observed by Shalit and Solel \cite{ShaSo09}, one may think of $H^\botimes$ with its product as an algebra generated by $d$ indeterminates subject to homogeneous relations. \cite[Theorem 8.4]{ShaSo09} showed that $\ol{\cA}$ is the universal operator algebra generated by a row contraction of $d$ operators subject to the same relations. \cite{DRS11} showed that the classification of the arising (non-selfadjoint) operator algebras $\ol{\cA}$ is the same as the classification of the associated subproduct systems is the same as the classification of the homogeneous relations up to permutations of the indeterminates, while the classification by the associated self-adjoint operator algebra $\ol{\cA^*}$ may be coarser. If the relations contain commutativity of the product, we are in the realm of polynomials in $d$ variables. \cite{DRS11} show that among subproduct systems that correspond to quotients by radical ideals, the classification is already done on the level of algebraic isomorphism of $\ol{\cA}$. Kakariadis and Shalit \cite{KaSha15p} do a similar program for the case of noncommuting \nbd{d}tuples. We ask:

\bOP
What are nice classes of bounded interacting Fock spaces that \bf{are} classified by their associated operator algebras $\ol{\cA^{(*)}}$?
\eOP

Natural suggestions for subclasses are all interacting Fock spaces of \bf{all} subproduct systems, not only of all finite-dimensional ones (for which \cite{KaSha15p} give the answer). Recall that subproduct systems lead to \nbd{\vk}interacting Fock spaces, where $\vk$ is a projection fulfilling an extra condition. This may be generalized to just any (squeezing) projection, or any bounded squeezing. In \ref{pivesSS}, we propose another class that arises from \it{nondegenerate productive systems}. Finally, we ask:

\bOP
How is the classification in terms of the tensor (Cuntz-Pimsner-Toeplitz) algebras into which $\ol{\cA^{(*)}}$ embeds? (Recall that there are different choices.) How, further, under the quotient to the Cuntz-Pimsner algebras?
\eOP

Kakariadis and Shalit \cite{KaSha15p} address some of these question in their framework.

\newpage

\subsection{About productive systems} \label{pivesSS}

Following the definition in Shalit and Skeide \cite[Section 6]{ShaSk10p}, a (discrete one-parameter) \hl{productive system} of Hilbert spaces is a family $H^<=\bfam{H_n}_{n\in\N_0}$ of Hilbert spaces $H_n$ with product maps $v_{m,n}\in\sB(H_m\sbar{\otimes}H_n,H_{m+n})$ fulfilling all requirements for the product maps of a subproduct system, except that they are not required coisometric. (In \cite{ShaSk10p}, the definition is for correspondences instead of just Hilbert spaces, and the indexing monoid can be arbitrary. There is also a \it{coproductive system}; even for Hilbert spaces, the two categories are different for the suggested morphisms, but we do not need this sophistication. A super(sub)product system is a (co)isometric productive system.) A productive system is \hl{nondegenerate} if $v_{m,n}(H_m\sbar{\otimes}H_n)$ is dense in  $H_{m+n}$. (A subproduct system is a nondegenerate productive system, while a superproduct system is nondegenerate if and only if it is a product system.) Sometimes, we will like that a productive system be \hl{contractive} (all $v_{m,n}$ are contractions) or \hl{bounded} (their norms are bounded uniformly).

Exactly as in the beginning of Section \ref{piSEC} (with the same product notation), for every $x\in H_1$ we define the creators $a^*(x)X_n:=xX_n$ ($X_n\in H_n$). If $H^<$ is bounded (contractive), then all creators are bounded (by $\norm{x}$) and  the creator map is bounded (contractive). (Only if the creators are bounded, they can be defined everywhere on the Fock space $\sF(H^<)$ of the productive system.) Similarly, we extract pre-Hilbert spaces
\beqn{
\ul{H\!}\,_n
~:=~
\ls{a^*(H_1)}^n\Om
~\subset~
H_n.
}\eeqn
The creators fulfill \eqref{cls***} if and only if $H^<$ is nondegenerate, in which case each $\ul{H\!}\,_n$ is dense in $H_n$. In any case, the family $\ul{H\!}\,_n$ gives rise to an interacting Fock space $\cI$ based on $H_1$ with creators $a^*(x)$ (co)restricted to $\cI$ and still denoted by $a^*(x)$. If $H^<$ is nondegenerate, then $\cI$ is dense in $\sF(H^<)$. But even if $H^<$ is not nondegenerate, then the Hilbert subspaces $\ol{\ul{H\!}\,_n}\subset H_n$ form a productive subsystem of $H^<$, which, now, is nondegenerate.

Obviously, $\Lambda$ is just given by the iterated products $v_{(n)}$ as
\beqn{
\Lambda_n(x_n\otimes\ldots\otimes x_1)
~=~
a^*(x_n)\ldots a^*(x_1)\Om
~=~
v_{(n)}(x_n\otimes\ldots\otimes x_1).
}\eeqn
It follows that $\Lambda_n$ is bounded. Therefore, $\Lambda$ is weakly adjointable, so $\cI$ is regular, hence, embeddable. If $H^<$ is contractive, then $\Lambda$ is a contraction.

Whatever the embedding $\xi$ and $\lambda=\xi\Lambda$ are, $\vk$ is the unique vacuum preserving operator that is $0$ on $H_1\otimes(\xi\cI)^\perp\subset\ol{\sF(H_1)}$ and that sends $x\otimes\lambda_nX_n$ to $\lambda_{n+1}(x\otimes X_n)$ for $X_n\in H_1^{\otimes n}$ ($n\ge0$). That is, the norm of $\vk_{n+1}$ it the same as the norm of
\beqn{
v_{1,n}
\colon
x\otimes\Lambda_nX_n
~=~
(\id_{H_1}\otimes v_{(n)})(x\otimes X_n)
~\longmapsto~
v_{1,n}(\id_{H_1}\otimes v_{(n)})(x\otimes X_n)
~=~
v_{(n+1)}(x\otimes X_n)
~=~
\Lambda_{n+1}(x\otimes X_n).
}\eeqn
Consequently, $\vk$ is bounded (a contraction), if $H^<$ is bounded (contractive).

So, analyzing the interacting Fock spaces of bounded productive systems, we have found the first instance (apart from subproduct systems, where we identified $\vk$ explicitly and it was bounded) of a class that have necessarily bounded $\vk$.

What else does it need for that a bounded squeezing $\vk$ determines an interacting Fock space $\cI$ that comes from a bounded productive system? Well, if $\cI$ comes from a productive system, then the product is recovered from $\Lambda$ as $v_{m,n}\colon\Lambda_mX_m\otimes\Lambda_nY_n\mapsto\Lambda_{m+n}(X_m\otimes Y_n)$. Recalling that we have a \nbd{\vk}interacting Fock space $\cI\subset\ol{\sF(H)}$ ($H$, for convenience, immediately assumed to be a Hilbert space), where $\lambda=\Lambda$, this reads
\beqn{
\lambda_mX_m\otimes\lambda_nY_n
~\longmapsto~
\lambda_{m+n}(X_m\otimes Y_n).
}\eeqn
If these maps are well-defined, then the corresponding product is manifestly associative. Recall from linear algebra that $\ker(\lambda_m\otimes\lambda_n)=\ls(\ker\lambda_m\otimes H_n)\cup(H_m\otimes\ker\lambda_n)$. From the recursion for $\lambda$ in terms of $\vk$, it follows that $\lambda_{m+n}$ vanishes on $H_m\otimes\ker\lambda_n$. If we find the analogous recursion $\lambda_{n+1}=\vk'_{n+1}(\lambda_n\otimes\id_H)$ for some $\vk'$  (for the same (\bf{!}) $\lambda_n$) from the other side, then $\lambda_{m+n}$ also vanishes on $\ker\lambda_m\otimes H_n$, and our product is well-defined. We tell why existence of $\vk'$ is necessary, in between Open Problems 7 and 8, below. Now, since obviously $\vk$ is, \it{cum grano salis} (that is, up to questions of (co)domain), $\bigoplus_{n\in\N_0}v_{1,n}$, we see that all $v_{1,n}$, hence, all $v_{m,n}$ are contractions if $\vk$ is. We have proved the following:

\vspace{-1ex}
\bthm \label{prodSkthm}
The interacting Fock space of any contractive productive system $H^<$ is isomorphic to a \nbd{\vk}interacting Fock space based on $H_1$ for a contractive squeezing $\vk$. 

Conversely, if $\cI$ is a \nbd{\vk}interacting Fock space based on a Hilbert space $H=H_1$, then $\cI$ is the interacting Fock space of a nondegenerate contractive productive system if and only if $\vk$ is a contraction and the there exits another (contraction) $\vk'$ such that the $\lambda$ constructed from $\vk$ fulfills $\lambda=\vk'((\lambda\otimes\id_H)\oplus\id_{H_0})$.
\ethm

\vspace{-2ex}
\bOP
What are the (contractive) squeezings $\vk$ that belong to bounded (contractive) productive systems? How can they be classified? (Of course, also the questions about the associated operator algebras are meaningful  for this subclass of bounded interacting Fock spaces.)
\eOP

\vspace{-4ex}
Obviously, one may model ALV-interacting Fock spaces and the whole theory that follows also for \hl{right interacting Fock spaces} starting in Definition \ref{ALVdefi} not from the left creators $\ell^*(x)$ but from the right creators $r^*(x)\colon X_n\mapsto X_n\otimes x$. (That would lead to call our interacting Fock spaces \hl{left}.) Interacting Fock spaces of productive systems appear, then, to be left \bf{and} right. (This entirely explains origin and properties of $\vk'$ in the preceding theorem.)

\bOP
Elaborate the precise relationship between productive systems and interacting Fock spaces that are left and right.
\eOP

\vspace{-3ex}
\bOP
(Entirely speculatively.) Is there a notion \it{dual} to interacting Fock space relating to the notion of nondegenerate \it{coproductive system} (generalizing superproduct system)?
\eOP

% \newpage
\subsection{About classification} \label{classSS}

One of the basic open classification problems, that of regularity, we described already in \ref{regSS}, showing that $\vk$ will, in general, not give a good answer. The present part is rather directed to point out possible strategies to give positive answers. A good strategy is to try to answer the question for subclasses; so, we basically propose more subclasses, ornamented with some preliminary insights.

As we know from Theorem \ref{nonregthm}, the non-nilpotent full interacting Fock spaces are all not regular: No matter how we base them, they will never be regularly based. But unbounded creators were essential in the proof. We asked already in Open Problem 3, which bounded interacting Fock spaces are regular.  Varying the hypothesis of Theorem \ref{nonregthm} by adding boundedness, this leads to the following questions:

\bOP
A \hl{bounded full interacting Fock space} is an interacting Fock space of the form $\cI=(\bfam{H_n}_{n\in\N_0},\sB_{(1)}(\cI))$. Is $\cI$ regular? Is $\cI$ regular, if it is also adjointable (so $A^*=\sB^a_{(1)}(\cI)$)? (Weak adjointability does not add anything, because all bounded operators are weakly adjointable.)
\eOP

We know from Example \ref{kappaunbex} that a bounded interacting Fock space based on $H$ (even with bounded creator map) need not be regularly based. But can it be regularly based by choosing a better basing? We address this question later, in the more general context of \ref{2FockSS}. Here we are interested in a property that, independently of the basing, the interacting Fock space in Example \ref{kappaunbex} possesses, but, that full interacting Fock spaces (bounded or not) do not possess. Likewise the chosen basing in Example \ref{kappaunbex} fulfills certain, apparently desirable, conditions we would like to add (in various combinations) as hypotheses to the question of regularity.

Almost all interacting Fock spaces we considered in this paper (and that have been considered elsewhere) are \hl{vacuum separated} in the sense that $a^*\Om=0$ implies $a^*=0$ for all $a^*\in A^*$. (This is a general property that does not refer to any basing.) The full interacting Fock spaces, however, are not vacuum separated (unless $H_2=\zero$). We will rather say \hl{\nbd{0}separated}, because we think that \hl{\nbd{n}separated}, meaning that $a^*H_n=\zero$ implies $a^*=0$, (and also the property to be \nbd{n}separated for all $n$ or for $n\le N$) will play a role in future discussion.

If $\cI$ is \nbd{0}separated, then the map $a^*\mapsto a^*\Om$ from $A^*$ to $H_1$ (which, we know, is surjective) is a bijection. Therefore, $\cI$ may be \hl{iso-based}, that is, based on $H_1$ in such a (unique, if possible) way such that $a^*(x)\Om=x$. Of course, any iso-based interacting Fock space is \nbd{0}separated. We summarize:

\bprop
The \nbd{0}separated interacting Fock spaces are exactly those that can be iso-based, and if an interacting Fock space can be iso-based then the \hl{iso-basing} $a^*\colon H_1\rightarrow A^*$ is unique.
\eprop

\bob
The \nbd{q}Fock spaces (Example \ref{POIex}\eqref{POI1}) and the interacting Fock spaces of productive systems are iso-based, hence, \nbd{0}separated. ALV-(POI-)interacting Fock spaces are iso-based if an only if $(\bullet,\bullet)_1=\AB{\bullet,\bullet}_1$ ($L_1=\id_{H_1}$). They are \nbd{0}separated (that is, can be iso-based) if and only if $(\bullet,\bullet)_1$ is an inner product ($L_1$ is injective).
\eob

Iso-based interacting Fock spaces fulfill another property, we mentioned already in Section \ref{absIFSEC}. They are \hl{injectively based} in the sense that the basing $a^*$ is an injective, hence, bijective map. That is, for an interacting Fock space that is injectively based on $H$, we may identify $H$ with $A^*$ as vector spaces, and $H$ induces an inner product on $A^*$. Conversely, starting with an interacting Fock space, we base it injectively on a pre-Hilbert space $H$ by defining just any inner product $(\bullet,\bullet)$ on $A^*$ and call the resulting pre-Hilbert space $H:=(A^*,(\bullet,\bullet))$. In Theorem \ref{aIFSembthm}, we even showed that the choice of $(\bullet,\bullet)$ can be done such that the basing is embeddable. Let us repeat that regarding regularity, we know this is not always possible. In fact one of the major questions in the background of this section is whether or not for a given (class of) interacting Fock space(s) the basing may be chosen regularly embeddable.

\it{A priori} we do not know, if the restriction to injective basing does influence the answer to the question. We show, it does not. In fact, Proposition \ref{reg-inregprop} below, is a simple consequence of the following obvious result.

\blem
Suppose $\Lambda\colon H\rightarrow I$ is an operator between pre-Hilbert spaces $H$ and $I$, and $\Lambda^*\colon I\rightarrow\ol{H}$ a weak adjoint. Then $p\Lambda^*\colon I\rightarrow\ol{G}$, with $p$ the projection from $\ol{H}$ onto $\ol{G}$, is a weak adjoint of $\Lambda\upharpoonright G$.
\elem

\bprop \label{reg-inregprop}
If an interacting Fock space is regular, then there is also a regular basing that is injective.
\eprop

\proof
Suppose we have an interacting Fock space $\cI$ regularly based on $H$ and denote by $\Lambda^*\colon\cI\rightarrow\ol{\sF(H)}$ the adjoint of $\Lambda$. Choose any subspace $G$ of $H$ such that the restriction of $a^*$ to $G$ is bijective, and apply the lemma to $\Lambda$ and the subspace $\sF(G)$ of the domain of $\Lambda$.\qed

\lf
So, as far as regularity is concerned, we do not lose anything restricting our attention to injective basings. The question for existence of a regular basing is, therefore, equivalent to the question of existence of an inner product on $A^*$ such that the corresponding injective basing on $H:=(A^*,(\bullet,\bullet))$ is regular. This is paired with the question if for a \nbd{0}separated interacting Fock space $\cI$ its iso-basing on $H_1$ is regular; if the answer is yes, we say $\cI$ is \hl{iso-regular}.

For the balance of this part, we examine examples of \nbd{0}separated interacting Fock spaces, their iso-basings, and the possibilities for other injective basings. They illustrate that the natural questions to be asked, do not have general answers, but depend on the case.

Regarding embeddability, we have seen that there are ``unfortunate'' choices for a basing; but for embeddability, Theorem \ref{aIFSembthm} tells us that the choice may be fixed (even injectively, but by Theorem \ref{nonregthm}, not necessarily regularly). A very simple class of examples tells that also regularity does depend on the choice of the basing: To show that a certain interacting Fock space is not regular, it is not sufficient to find just one non-regular basing.

\bex
(See also Example \ref{nonembedex}.) We look at full interacting Fock spaces of the form $\cI_1=\Om\C\oplus H_1$, so that $A^*:=\sL^a_{(1)}(H_0,H_1)=\sB^a_{(1)}(H_0,H_1)=H_1\Om^*$. Apart from being bounded, $\cI$ is also adjointable and \nbd{0}separated. Obviously, the iso-basing $a^*(x):=x\Om^*$ is regular, and we may consider $\cI_1$ as embedded via the canonical identification $\cI_1\subset\sF(H_1)$. (Then $\Lambda$ is actually the adjoint of $\xi$, $\lambda$ is the projection onto $\Om\C\oplus H_1\subset\sF(H_1)$. Up to possible missing completeness, $\cI$ is the interacting Fock space of the simplest nontrivial subproduct systems possible, and $\vk=\lambda$ the projection identifying it as \nbd{\vk}interacting Fock space.) But we may choose ``less fortunate'' basings.

In fact, any other basing $a'^*\colon H\rightarrow A^*$ factors through $a^*$ via a surjective linear map $T\colon H\rightarrow H_1$ in the sense that $a'^*=a^*T$. In fact, for the basing $a'^*$, we find $\Lambda'=\id_{H_0}\oplus T$. So, the basing $a'^*$ is regular if and only $T$ is weakly adjointable. Already when $H_1=\Om_1\C$ is one-dimensional, we may choose $T=\Om_1\vp$, where $\vp$ is unbounded. There is the entirely justified objection that $H$ has to be infinite-dimensional, so that $a'^*$ for this $T$ is not at all injective. However, if $H_1$ is not finite-dimensional, for instance,if we put $H:=H_1$, then there are invertible maps $T$ on $H_1$ that are not weakly adjointable. Here is an example:

Put $H_1=\ls\CB{e_n\colon n\in\N}$ for some orthonormal family $\bfam{e_n}_{n\in\N}$. Then $T\colon e_n\mapsto e_n+ne_1$ is, clearly, a bijection. But from $\AB{e_1,Te_n}=n$, it follows that there is no vector $x=T^*e_1$ fulfilling $\AB{x,e_n}=n$.
\eex

With this example and the following case study in mind, which shows that all three classes are different, we pose:

\bOP
Classify \nbd{0}separated interacting Fock spaces. Among those, classify the regular ones. Among those, classify the iso-regular ones. ``Classification'' means any sense of classification discussed in Section \ref{olSEC} or still to be uncovered elsewhere, when restricted to this hierarchy of subclasses.
\eOP

\newpage

\subsection{A case study: Interacting Fock 2--spaces} \label{2FockSS}

Before examining in detail the one step more complicated situation, $\cI_2^0:=\Om\C\oplus H_1\oplus\Om_2\C$, let us refine our notations regarding a general \nbd{0}separated interacting Fock space $\cI$ and its iso-basing $a^*\colon H_1\rightarrow A^*$. Denote by $a_n\in\sL(H_n,H_{n+1})$ the (co)restriction of $a\in\sL_{(1)}(\cI)$ to the occurring subspaces of $\cI$. Then
\beqn{
a
~=~
a_0\oplus a_1\oplus\ldots,
}\eeqn
when considering $a$ as map into $H_1\oplus H_2\oplus\ldots=\cI\ominus H_0$. Applying this to a creator $a^*(x)$ (in iso-basing!), we get
\beqn{
a^*(x)
~=~
x\Om^*+a^*_1(x)+a^*_2(x)+\ldots,
}\eeqn
where, in a sense, the (unique!) linear maps $a^*_n\colon x\mapsto a^*_n(x)\in\sL(H_n,H_{n+1})\subset\sL_{(1)}(\cI)$ capture the entire structure of $\cI$ in iso-basing.

Now let us fix $\cI_2^0:=\Om\C\oplus H_1\oplus\Om_2\C$ for some pre-Hilbert space $H_1$, unit vector $\Om_2\in H_2$ such that $H_2=\Om_2\C$ and a linear map $0\ne a^*_1\in\sL(H_1,\sL(H_1,H_2))$ so that $\cI_2^0$ with $a^*(x):=x\Om^*+a_1^*(x)$ is a \nbd{0}separated interacting Fock space in iso-basing. We call $\cI_2^0$ an \hl{interacting Fock \nbd{2}space}. (An \hl{interacting Fock \nbd{n}space} would be a \nbd{0}separated interacting Fock space in iso-basing satisfying $H_n=\Om_n\C$ and $H_{n+1}=\zero$.) We shall write $a^*_1=\Om_2\vp\colon x\mapsto\Om_2\vp_x$ where $\vp\colon x\mapsto\vp_x$ is a linear map from $H_1$ into the linear functionals $\sL(H_1,\C)$ on $H_1$, so that $a^*_1(x)y=\Om_2(\vp_xy)$.

The map $\vp$, which characterizes $\cI_2^0$, may equally well be described by the linear functional $\Phi\in\sL(H_1\otimes H_1,\C)$ induced via the universal property of the tensor product $H_1\otimes H_1$ by the bilinear map $(x,y)\mapsto\vp_xy$. The only condition to be satisfied is that either of them is nonzero. For $\Lambda_1$ we find $\Lambda_1(x)=a^*(x)\Om=x$, so $\Lambda_1=\id_{H_1}$. For $\Lambda_2$ we find
\beqn{
\Lambda_2(x\otimes y)
~=~
a^*(x)a^*(y)\Om
~=~
a^*_1(x)y
~=~
\Om_2(\vp_xy)
~=~
\Om_2\Phi(x\otimes y),
}\eeqn
so $\Lambda_2=\Om_2\Phi$. Therefore:

\bprop \label{20isoregprop}
$\cI_2^0$ is iso-regular if and only if $\Phi$ is bounded.

Consequently, there exist iso-based \nbd{0}separated interacting Fock spaces that are non iso-regular.
\eprop

\brem
Note that if $\Phi$ is bounded and $H_1$ complete, then also $\cI_2^0$ is the interacting Fock space of a productive system with the only non-obvious product map given by $v_{1,1}=\Lambda_2$.
\erem

\brem \label{isoembrem}
Note, too, that choosing for $\Om_2$ a unit vector in $H_1\otimes H_1$, we identify $H_2=\Om_2\C$ as a subspace of $H_1\otimes H_1$. For the interacting Fock space $\cI_2^0$ based on $H_1$ and embedded this way in $\sF(H_1)$, we find $\lambda=\Lambda$ (coextended as map into $\sF(H_1)$), and $\vk=\lambda$.
\erem

Of course, as discussed in Example \ref{kappaunbex}, we knew the latter statement of Proposition \ref{20isoregprop} already for the special case in Examples \ref{CPTex} and \ref{kappaunbex}, which, in fact, is an interacting Fock \nbd{2}space. Here, we recover that statement as a part of a more general situation.

We now switch out attention to the question if $\cI_2^0$ is regular, that is, if, changing that basing, we can turn $\cI_2^0$ into an interacting Fock space that is regularly based. By Proposition \ref{reg-inregprop}, it is sufficient to look at injective basings, only. So, given $\Phi$, can we change the inner product of $H_1$ in such a way that $\Phi$ becomes bounded?

We prefer to formulate the question in a slightly more abstract way: Given a vector space $H$ and a linear functional $\Phi\colon H\otimes H\rightarrow\C$, does there exist an inner product on $H$ such that $\Phi$ is bounded? The answer -- no in general, but yes if $H$ has a countable basis -- is provided by the following (counter) Example \ref{Phinobex} and Theorem \ref{Phicbthm}.

\bex \label{Phinobex}
As frequently with spaces $H$ that may be viewed as a space of functions on $\SB{0,1}$ (or other subsets of $\R_+$ with accumulation points), if problems can be caused in $H\otimes H$, then they arise by looking at the \hl{diagonal} $D:=\CB{(t,t)\colon t\in\SB{0,1}}$ ~of~ $\SB{0,1}\times\SB{0,1}$.

Choose $H$ to be a vector space with basis $\bfam{e_t}_{t\in\SB{0,1}}$. Define $\Phi\in\sL(H\otimes H,\C)$ by setting
\beqn{
\Phi(e_s\otimes e_t)
~:=~
\begin{cases}
\frac{1}{s-t}&s\ne t,
\\
0&\text{otherwise.}
\end{cases}
}\eeqn
Then, for whatever inner product we might choose on $H$, the functional $\Phi$ is unbounded. In the following lemma, we prove a more general statement, which might be useful also for general tensor products of general normed spaces.
\eex

\blem
For whatever norm $\norm{\bullet}$ we choose on $H$, there is no (sub)cross norm on $H\otimes H$ that made $\Phi$ bounded.
\elem

\proof
Define $S_n:=\CB{t\in\SB{0,1}\colon\norm{e_t}\le n}$. Since $\SB{0,1}$ is the countable(!) union of all $S_n$, from a certain $N$ on all $S_n$ $(n\ge N$) are uncountable. The infinite set $S_N$ has an accumulation point, say, $t_0$. For each $\ve>0$, the intersection $(U_{\frac{\ve}{2}}(t_0)\times U_{\frac{\ve}{2}}(t_0))\cap(S_N\times S_N\backslash D)$ is nonempty, so there are $s\ne t\in S_N$ such that $\abs{\Phi(e_s\otimes e_t)}>\frac{1}{\ve}$. Therefore,
\beqn{
\norm{\Phi}
~=~
\sup_{0\ne X\in H\otimes H}\frac{\abs{\Phi(X)}}{\norm{X}}
~\ge~
\sup_{s,t\in S_N}\frac{\abs{\Phi(e_s\otimes e_t)}}{\norm{e_s\otimes e_t}}
~\ge~
\sup_{s,t\in S_N}\frac{\abs{\Phi(e_s\otimes e_t)}}{\norm{e_s}\norm{e_t}}
~\ge~
\sup_{s,t\in S_N}\frac{\abs{\Phi(e_s\otimes e_t)}}{N^2}
~\ge~
\frac{1}{\ve N^2}.
}\eeqn
(The second ``$\ge$'' follows from sub-cross; if it was cross, as for pre-Hilbert spaces, then it would be ``$=$''.) Since $\ve>0$ was arbitrary, $\norm{\Phi}=\infty$.\qed

\lf
Note that $H$ may be separable. (Just take the one-mode symmetric Fock space with the (dense) subspace spanned by the exponential vectors $\ee(t)$ ($t\in\SB{0,1}$); see Example \ref{wpnonpex}.) That the index set $\SB{0,1}$ of the Hamel basis is uncountable, is crucial for the proof. (Otherwise, we cannot show existence of an $S_N$ with an accumulation point.) In fact:

\bthm \label{Phicbthm}
Let $H$ be a vector space with a countable basis $\bfam{e_n}_{n\in\N}$. Then for every linear functional $\Phi\in\sL(H\otimes H,\C)$ on $H\otimes H$, there exists an inner product on $H$ such that $\Phi$ is bounded.
\ethm

\bcor
~Every interacting Fock \nbd{2}space with a countable Hamel basis is regular (though, not necessarily iso-regular).
\ecor

(This result has some similarity with \cite[Theorem 5.3]{AcSk08}, which asserts that an interacting Fock space based on pre-Hilbert space with countable Hamel basis is even algebraically embeddable.)

For the proof of Theorem \ref{Phicbthm}, we need preparation.

\bob
For any function $F\colon\N\times\N\rightarrow\R_+$ define the function $f\colon\N\rightarrow\R_+$ by
\beqn{
f(n)
~:=~
\max\CB{1,F(i,j)\colon i,j\le n}.
}\eeqn
Then
\beqn{
F(i,j)
~\le~
f(i)f(j)
}\eeqn
for all $i,j$. (Indeed,
$F(i,j)
~\le~
f(\max\CB{i,j})
~\le~
f(\max\CB{i,j})\,f(\min\CB{i,j})
~=~
f(i)f(j)
$.)\eob

\proof[Proof of Theorem \ref{Phicbthm}.~]
(We thank Uwe Franz for assistance.)

For $F(i,j):=\sabs{\Phi(e_i\otimes e_j)}$, choose $f$ as in the observation and put $c_n:=2^nf(n)\ne0$. On $H$ define an inner product by setting $\AB{e_i,e_j}:=\delta_{i,j}c_i^2$. For $v=\sum_{i,j}e_i\otimes e_j\lambda_{i,j}\in H\otimes H$ (so, $\lambda_{i,j}\ne0$ only for finitely many $i,j$), we find
\beqn{
\norm{v}^2
~=~
\sum_{i,j}\sabs{\lambda_{i,j}}^2c_i^2c_j^2.
}\eeqn
Therefore,
\bmun{
\abs{\Phi(v)}
~\le~
\sum_{i,j}\sabs{\lambda_{i,j}}\sabs{\Phi(e_i\otimes e_j)}
~\le~
\sum_{i,j}\sabs{\lambda_{i,j}}f(i)f(j)
~=~
\sum_{i,j}\sabs{\lambda_{i,j}}c_ic_j\textstyle\frac{1}{2^i2^j}
\\
~\le~
\textstyle\sqrt{\Big.\sum_{i,j}\bfam{\sabs{\lambda_{i,j}}c_ic_j}^2}\,\sqrt{\Big.\sum_{i,j}\bfam{\frac{1}{2^i2^j}}^2}
~=~
\norm{v}\frac{1}{3},
}\emun
where the step from the first to the second line is Cauchy-Schwartz inequality for $\ell^2(\N\times\N)$.\qed

\bOP
The interacting Fock space in Examples \ref{CPTex} and \ref{kappaunbex} is $\cI_2^0$ with $\Phi(x\otimes y):=\AB{\bar{x},y}$ for some anti-unitary involution on $H$. By Theorem \ref{Phicbthm}, it is regular if $H$ has a countable Hamel basis. Is it always regular?
\eOP

The simple characterization of iso-regular interacting Fock \nbd{2}spaces as those with bounded $\Phi$, is thanks to finite-dimensionality of $H_2$. Already in the \nbd{3}space with both $H_1$ and $H_2$ infinite-dimensional we do not know what conditions for $a^*_1$ will pop up (while $a^*_2$ still has to be bounded for the same reason).

\bOP
How much of this case study goes through for interacting Fock \nbd{n}spaces for $n=3$ or bigger? (For instance, embeddability for the iso-basing, as discussed for $n=2$ in Remark \ref{isoembrem}, for $n=3$ may easily fail.)
\eOP

\newpage

\subsection{About automorphism groups} \label{autoSS}

Apart from phrasing the ``natural'' questions about classifying objects by looking at their automorphism groups, we add to the two notions of isomorphim that we defined already in Definition \ref{isodefi}, two more in Definition \ref{misodefi}. The first one, \it{quasi-isomorphism}, adds to isomorphisms of interacting Fock spaces based on pre-Hilbert spaces, a ``unitary freedom'' in choosing the ``parameter space'' $H$. The second one, \it{vague isomorphism}, aims at incorporating different choices of injective basings. Each of the four notions of isomorphism (plus, possibly, others that have not yet been invented) lead to a different notion of automorphism group. Therefore, all our ``natural'' questions, actually, have four versions of them.

\bOP
What are the automorphism groups of some concrete (classes of) interacting Fock spaces?
\eOP

\bOP
To what extent are (classes of) interacting Fock spaces determined by their automorphism groups?
\eOP

\bOP
How is this classification in terms of automorphism groups reflected by the (several!) associated operator algebras?
\eOP

\bdefi \label{misodefi}
The interacting Fock space $\cI=(\bfam{H_n}_{n\in\N_0},a^*)$ based on $H$ and the interacting Fock space $\cI'=(\bfam{H'_n}_{n\in\N_0},{a^*}')$ based on $H'$ are \hl{quasi-isomorphic} if there exists a unitary $U\colon H\rightarrow H'$ such that $\cI$ and $(\bfam{H'_n}_{n\in\N_0},{a^*}'\circ U)$ are isomorphic interacting Fock spaces based on $H$.

$\cI$ and $\cI'$ are \hl{vaguely isomorphic} if we can find an invertible linear map $T\colon H\rightarrow H'$ such that $\cI$ and $(\bfam{H'_n}_{n\in\N_0},{a^*}'\circ T)$ are isomorphic interacting Fock spaces based on $H$.
\edefi

We dispense, for now, with the idea to make a list of more or less obvious properties in the style of Observation \ref{isoob}. Anyway, this had to be repeated when we go \it{in medias res} with this program.

% \listofOWs
\newpage

\newpage

\lf\noindent
\bf{Acknowledgments.~}
We wish to thank Orr Shalit for very useful discussions about \cite{ShaSo09,DRS11,KaSha15p}. We wish to thank Roland Speicher for the hint that the condition in \cite{AcSk08} for an interacting Fock space to be \it{embeddable} might be superfluous, which we answer here by ``\it{yes, if}''. This made the results from \cite{AcSk08} applicable without any limitation. We wish to thank Uwe Franz for assistance in the proof of Theorem \ref{Phicbthm}. A big thank you to the referees who produced the two reports.

MG acknowledges funding from the German Research Foundation (DFG) through the proj\-ect "Non-Commutative Stochastic Independence: Algebraic and Analytic Aspects", project number 397960675.

\setlength{\baselineskip}{2.5ex}

% \bibliography{mybib}
\newcommand{\Swap}[2]{#2#1}\newcommand{\Sort}[1]{}
\providecommand{\bysame}{\leavevmode\hbox to3em{\hrulefill}\thinspace}
\providecommand{\MR}{\relax\ifhmode\unskip\space\fi MR }
% \MRhref is called by the amsart/book/proc definition of \MR.
\providecommand{\MRhref}[2]{%
  \href{http://www.ams.org/mathscinet-getitem?mr=#1}{#2}
}
\providecommand{\href}[2]{#2}

\lf\noindent
Malte Gerhold:
{\small\itshape Institut für Mathematik und Informatik}, \\
{\small\itshape Universität Greifswald,
17487 Greifswald, Germany, \\E-mail: \href{mailto:mgerhold@uni-greifswald.de}{\tt{mgerhold@uni-greifswald.de}}}

\noindent
and

\noindent
{\small\itshape Faculty of Mathematics}, \\
{\small\itshape Technion Israel Institute of Technology,
Haifa 3200003, Israel, \\E-mail: \href{mailto:maltegerhold@campus.technion.ac.il}{\tt{maltegerhold@campus.technion.ac.il}}}

\lf\noindent
Michael Skeide:
{\small\itshape Dipartimento di Economia, Universit\`a\ degli Studi del Molise, Via de Sanctis, 86100 Campobasso, Italy, E-mail: \href{mailto:skeide@unimol.it}{\tt{skeide@unimol.it}}}\\
{\small{\itshape Homepage: \url{http://web.unimol.it/skeide/}}}

% \newpage
% \listofOWs

\end{document}